\documentclass[a4paper,12pt,leqno, twoside]{article}
\usepackage{amssymb}
\usepackage{amsthm}
\usepackage{amstext}
\usepackage{amsbsy}
\usepackage{amscd}
\usepackage{amsopn}
\usepackage{amsmath}
\usepackage{ stmaryrd } 
\usepackage{color}
\usepackage{mathtools}
\usepackage{mathdots}
\usepackage[mathscr]{eucal}
\usepackage{graphicx}
\usepackage{makeidx}
\makeindex
\usepackage{epic,eepic}
\usepackage{epsfig}
\usepackage[all]{xy}
\usepackage[french]{babel}
\usepackage[latin1]{inputenc}
\usepackage[nottoc]{tocbibind}

    \def\AM{{\mathbb{A}}}
\def\BG{{\mathfrak B}}    
    \def\CM{{\mathbb{C}}}

    \def\FM{{\mathbb{F}}}
    \def\GM{{\mathbb{G}}}
    \def\HM{{\mathbb{H}}}

    \def\NM{{\mathbb{N}}}
    
    \def\PM{{\mathbb{P}}}
    \def\QM{{\mathbb{Q}}}
    \def\RM{{\mathbb{R}}}
\def\SG{{\mathfrak S}}

\def\XG{{\mathfrak X}}    
    
    \def\ZM{{\mathbb{Z}}}

\def\Bb{{\mathbf B}}    \def\BC{{\mathcal{B}}}
    \def\CC{{\mathcal{C}}}

    \def\FC{{\mathcal{F}}}
\def\Gb{{\mathbf G}}    
\def\Hb{{\mathbf H}}    \def\HC{{\mathcal{H}}}
    
\def\Jb{{\mathbf J}}    
    
\def\Lb{{\mathbf L}}    \def\LC{{\mathcal{L}}}
    \def\MC{{\mathcal{M}}}
    
    \def\OC{{\mathcal{O}}}
\def\Pb{{\mathbf P}}

\def\Tb{{\mathbf T}}    \def\TC{{\mathcal{T}}}
\def\Ub{{\mathbf U}}    
\def\Vb{{\mathbf V}}    
  \def\wb{{\mathbf w}}

          \def\wdo{{\dot{w}}}

\def\a{\alpha}

\def\G{\Gamma}
\def\d{\delta}
\def\D{\Delta}
\def\e{\varepsilon}
\def\ph{\varphi}
\def\l{\lambda}
\def\L{\Lambda}

\def\o{\omega}
\def\O{\Omega}

\def\s{\sigma}
\def\Sig{\Sigma}
\def\th{\theta}
\def\Th{\Theta}

\def\z{\zeta}

\def\Lamb{{\boldsymbol{\Lambda}}}

\DeclareMathOperator{\Aut}{{\mathrm{Aut}}}

\DeclareMathOperator{\Cosoc}{{\mathrm{Cosoc}}}

\DeclareMathOperator{\diag}{{\mathrm{diag}}}

\DeclareMathOperator{\End}{{\mathrm{End}}}

\DeclareMathOperator{\Fr}{{\mathrm{Fr}}}
\DeclareMathOperator{\Gal}{{\mathrm{Gal}}}
\DeclareMathOperator{\Hom}{{\mathrm{Hom}}}

\DeclareMathOperator{\im}{{\mathrm{Im}}}
\DeclareMathOperator{\ind}{{\mathrm{ind}}}
\DeclareMathOperator{\Ind}{{\mathrm{Ind}}}

\DeclareMathOperator{\Ker}{{\mathrm{Ker}}}

\DeclareMathOperator{\res}{{\mathrm{res}}}
\DeclareMathOperator{\Res}{{\mathrm{Res}}}

\DeclareMathOperator{\Soc}{{\mathrm{Soc}}}

\DeclareMathOperator{\Spec}{{\mathrm{Spec}}}

\DeclareMathOperator{\tors}{{\mathrm{tor}}}

\DeclareMathOperator{\Trace}{{\mathrm{Trace}}}

\DeclareMathOperator{\val}{{\mathrm{val}}}

\def\oQl{\o\QM_{\ell}}
\def\oZl{\o\ZM_{\ell}}

\def\oFq{\o{\FM}_q}

\def\ssi{si et seulement si }
\def\para{sous-groupe parabolique }

\def\levi{sous-groupe de Levi }
\def\levis{sous-groupes de Levi }

\def\borel{sous-groupe de Borel }

\def\Art {\mathrm {Art}}
\def\Nr {\mathrm {Nr}}

\def\Gal{{\rm Gal}}

\def\Frob{{\rm Frob}}

\def\Isom{{\rm Isom}}

\def\spe{{\rm sp}}

\def\Rep{{\rm Rep}}

\def\simto{\buildrel\hbox{$\sim$}\over\longrightarrow}

\def\leq{\leqslant}
\def\geq{\geqslant}
\def\into{\hookrightarrow}
\def\onto{\twoheadrightarrow}
\def\id{\mathop{\mathrm{Id}}\nolimits}

\def\ba{\backslash}
\def\wt{\widetilde}
\def\wh{\widehat}
\def\o#1{\overline{#1}}

\def\To#1{\buildrel\hbox{\tiny{$#1$}}\over\longrightarrow}

\def\to{\rightarrow}

\def\Hom{\mathop{\hbox{\rm Hom}}\nolimits}

\def\dim{\mathop{\mbox{\rm dim}}\nolimits}

\def \limi#1{\lim\limits_{\displaystyle\longrightarrow\atop {#1}}}

\def\ini{\setcounter{equation}{\value{subsubsection}}\addtocounter{subsubsection}{1}}

\makeatletter
\renewcommand{\subsubsection}{\@startsection{subsubsection}{3}{\parindent}{-\baselineskip}{-0.01\baselineskip}{\bf}}
\renewcommand*{\@seccntformat}[1]{%
  \csname the#1\endcsname\
} \makeatother

\def\ali{\subsubsection{}\setcounter{equation}{0}}

\newtheoremstyle{th}
  {\baselineskip}{.5\baselineskip}{\itshape}
  {\parindent}{\bf}
  {--}{.5em}{}

\newtheoremstyle{def}
  {\baselineskip}{\baselineskip}{}
  {\parindent}{\bf}
  {--}{.5em}{}

\newtheoremstyle{th*}
  {.5\baselineskip}{.5\baselineskip}{\itshape}
  {\parindent}{\bf}
  {--}{.5em}{}

\newtheoremstyle{remark*}
  {.5\baselineskip}{.5\baselineskip}{}
  {\parindent}{\bf}
  {--}{.5em}{}

\newtheoremstyle{remark}
  {.5\baselineskip}{.5\baselineskip}{}
  {\parindent}{\bf}
  {--}{.5em}{}

\swapnumbers \theoremstyle{th}
\newtheorem{theo}[subsubsection]{\sc Th{\'e}or{\`e}me.\bf}
\newtheorem{lemme}[subsubsection]{\sc Lemme.\bf}
\newtheorem{prop}[subsubsection]{\sc Proposition.\bf}
\newtheorem{coro}[subsubsection]{\sc Corollaire.\bf}

\theoremstyle{def}
\newtheorem{fact}[subsubsection]{\sc Fait\bf}

\newtheorem{DEf}[subsubsection]{\sc D{\'e}finition.\bf}
\newtheorem*{Def}{\sc D\'efinition.\bf}

\theoremstyle{remark}

\theoremstyle{th*}

\newtheorem*{lem}{\sc Lemme.}

\theoremstyle{remark*}
\newtheorem*{rem}{\sc Remarque.}
\newtheorem*{fait}{\sc Fait.}

\newtheorem*{exe}{\sc Exemple.}

\newcommand{\findem}{\hfill$\Box$\par\medskip}
\newcommand{\dem}{\indent {\it Preuve :} \rm }

\newenvironment{preuve}{\dem}{\findem}

\def\LJ{{\rm LJ}}
\def\JL{{\rm JL}}

\def\GL{{\rm GL}}
\def\Spf{{\rm Spf}}
\def\DL{{\rm DL}}
\def\PGL{{\rm PGL}}
\def\Coef{{\rm Coef}}

\theoremstyle{plain}

\begin{document}

\title{L'espace symm\'etrique de Drinfeld et correspondance de Langlands locale II}
\author{Haoran Wang}
\date{}

\maketitle


\begin{center}
{\bf {Abstract}}
\end{center}

We study the geometry and the cohomology of the tamely ramified cover of Drinfeld's $p$-adic symmetric space over a $p$-adic field $K.$ For this tame level, we prove, in a purely local way, most of a conjecture of Harris on the form of the $\ell$-adic cohomologies, as well as the torsion freeness of the integral cohomology. In this paper, we also compute the $\ell$-adic cohomology of Coxeter Deligne-Lusztig variety associated to $\GL_d,$ and some results of independent interest on the coefficient systems over the Bruhat-Tits building associated to $\GL_d(K)$ have been established.

\vspace*{6pt}\tableofcontents  

\renewcommand{\proofname}{\indent Preuve}

%
%


\section{Introduction}

Soit $K$ un corps local de caract\'eristique r\'esiduelle $p,$ d'anneau des entiers $\OC$ et de corps r\'esiduel $\FM_q.$ Pour $d\geq2,$ Drinfeld a introduit dans \cite{drinfeld-ell} son fameux ``espace sym\'etrique $p$-adique'' $\O^{d-1}$ ($\PM^{d-1}_K$ priv\'e des hyperplans $K$-rationnels), et il a trouv\'e que $\O^{d-1}$ poss\`ede un syst\`eme projectif de rev\^etements \'etales des espaces rigide-analytiques au sens de Raynaud-Berkovich $\{\Sig_n\}$ (on l'appelle la tour de Drinfeld).

La cohomologie de la tour de Drinfeld est li\'ee \`a la correspondance de Langlands locale. Elle a \'et\'e conjecturalememnt d\'ecrite par Harris dans un travail non publi\'e bas\'e sur le calcul de Schneider-Stuhler dans \cite{SS-invent} et son article \cite{Harris-Carayol}. Cette conjecture est maintenant connue, de mani\`ere indirecte, en combinant le th\'eor\`eme de Faltings-Fargues qui permet de passer \`a la tour de Lubin-Tate, et le travail de Boyer \cite{Boyer-Harris} qui d\'ecrit la cohomologie de la tour de Lubin-Tate par une approche de nature globale. En consid\'erant le complexe de cohomologie de la tour de Drinfeld comme un objet de la cat\'egorie d\'eriv\'ee des repr\'esentations lisses, Dat \cite{Dat-elliptic} en a tir\'e un raffinement spectaculaire: ce complexe de cohomologie r\'ealise \`a la fois la partie elliptique de la correspondance de Langlands locale et une forme de la correspondance de Jacquet-Langlands locale.

Cet article fait suit \`a \cite{Wang-Sigma1}. On \'etudie, de mani\`ere purement {\em locale}, la partie {\em non-supercuspidale} de la cohomologie $\ell$-adique ainsi que la cohomologie \`a coefficients entiers du niveau mod\'er\'e $\Sig_1$ de la tour ou plut\^ot d'une variante $\MC_{Dr,1}/\varpi^\ZM$ qui est en fait une r\'eunion de $d$ copies de $\Sig_1.$

Soient  $W_K$ le groupe de Weil de $K$ et $D$ l'alg\`ebre \`a division centrale d'invariant $1/d$ sur $K.$ La cohomologie de $\MC_{Dr,1}/\varpi^\ZM$ est munie d'une action des trois groupes $\GL_d(K),$ $W_K$ et $D^\times.$ Pour la d\'ecrire, on rappelle que si $\rho$ est une repr\'esentation irr\'eductible de $D^\times,$ alors la correspondance de Jacquet-Langlands lui associe une repr\'esentation $\JL(\rho)$ de $\GL_d(K).$ De plus, la correspondante de Langlands semi-simplifi\'ee $L(\JL(\rho))^{ss}$ est de la forme $\s_\rho\oplus\s_\rho(-1)\oplus\cdots\oplus \s_\rho(1-e)$ pour un diviseur $e$ de $d.$ On appelle alors ``repr\'esentation elliptique de type $\rho$'' toute repr\'esentation irr\'eductible $\pi$ de $\GL_d(K)$ telle que $L(\pi)^{ss}=L(\JL(\rho))^{ss}.$ Ces repr\'esentations sont param\'etr\'ees $I\mapsto \pi^I_\rho$ par les sous-ensembles de $\{1,\ldots,e\},$ {\em cf.} \cite{Dat-elliptic}.

Notre r\'esultat est le suivant:

\medskip
\noindent{\bf Th\'eor\`eme A.}
(\ref{Cor::2} \ref{Cor::3} \ref{4.3Prop:1} \ref{4.3Prop:2} \ref{4.3Cor1} \ref{Prop::1})
{\it
\begin{description}
\item [(a)] Pour tout $q\in\NM,$ le $\oZl$-module $H^{q}_c(\MC_{Dr,1}/\varpi^\ZM,\oZl)$ est admissible en tant que $\GL_d(K)$-module; il est sans-torsion et non-divisible.

\item [(b)] Soit $\rho$ une repr\'esentation irr\'eductible de niveau z\'ero de $D^\times$ de caract\`ere central trivial. Alors,
      $$
\Hom_{D^\times}\big(\rho,H^{d-1+i}_c(\MC_{Dr,1}/\varpi^\ZM,\oQl)\big)\cong\left\{
                                       \begin{array}{ll}
                                         \pi_\rho^{I_i}\otimes \s_\rho(-i), & \hbox{si $i\in\{0,\ldots,e-1\}$;} \\
                                         0, & \hbox{sinon,}
                                       \end{array}
                                     \right.
$$
o\`u $I_i=\{1,\ldots,i\}$ ou $\{e-i,\ldots,e-1\}.$
\end{description}}

\begin{rem}
(i) On a \'etudi\'e la partie supercuspidale de la cohomologie $\ell$-adique dans \cite{Wang-Sigma1}.

(ii) Gr\^ace \`a la m\'ethode globale (Boyer \cite{Boyer-Harris} et Faltings-Fargues \cite{FGL}), on sait que $I_i=\{1,\ldots,i\}.$

(iii) Boyer \cite{Boyer-torsionfree} a r\'ecemment annonc\'e une preuve de l'absence de torsion dans la cohomologie enti\`ere de la tour de Lubin-Tate.
\end{rem}
\medskip

Notre m\'ethode repose sur l'existence d'une suite spectrale associ\'ee \`a certain recouvrement ouvert de $\Sig_1.$ Ceci nous permet de calculer la cohomologie $H^q_c(\Sig_1,\L)$ ($\L=\oZl$ ou $\oQl$) via un {\em syst\`eme de coefficients} sur l'immeuble de Bruhat-Tits $\BC\TC$ associ\'e \`a $G:=\GL_d(K).$ Plus pr\'ecis\'ement, rappelons qu'il existe une application $\tau$ de $\O^{d-1}$ vers la r\'ealisation g\'eom\'etrique $|\BC\TC|$ de $\BC\TC.$
En prenant la composition avec la transition $p:\Sig_1\to\O^{d-1},$ on obtient un morphisme $G^\circ$-\'equivariant $\nu:\Sig_1\To{}|\BC\TC|,$ o\`u $G^\circ:=\{g\in G~|~\det(g)\in \OC^\times\}.$ Consid\'erons alors le recouvrement admissible $\{\nu^{-1}(|\s|^*)\}_{\s\in\BC\TC}$ de $\Sig_1,$ o\`u $|\s|^*$ est la r\'eunion de toutes les $|\s'|$ de la r\'ealisation g\'eom\'etrique de $\s'$ avec $\s'$ contenant $\s.$

\begin{fait}
Notons $\BC\TC_k$ l'ensemble des simplexes de dimension $k.$ Le complexe de C\v{e}ch associ\'e au recouvrement ci-dessus nous fournit une suite spectrale $G^\circ$-\'equivariante
\ini\begin{equation}\label{Intro1}
E^{pq}_1=\bigoplus_{\s\in\BC\TC_{-p}} H^q_c(\nu^{-1}(|\s|^*),\Lambda)\Longrightarrow H^{p+q}_c(\Sigma_1,\Lambda), \footnote{On triche un peu ici; il faut consid\'erer l'orientation de $\s.$}
\end{equation}
dont la diff\'erentielle $d^{pq}_1$ est celle du complexe de cha\^ines du syst\`eme de coefficients $\s\mapsto H^q_c(\nu^{-1}(|\s|^*),\Lambda).$
\end{fait}

On renvoie les lecteurs \`a 2.1 pour la notion de syst\`eme de coefficients.

Dans \cite{Wang-Sigma1}, nous avons calcul\'e les $H^q_c(\nu^{-1}(|s|^*),\Lambda)$ pour les sommets $s\in\BC\TC$ ainsi que les $H^q_c(\nu^{-1}(|\s|^*),\Lambda)$ en terme de $H^q_c(\nu^{-1}(|s|^*),\Lambda)$ avec $\s$ contenant $s$ (voir les rappels dans \ref{3.1::T1}). En particulier, on a montr\'e que le syst\`eme de coefficients $\s\mapsto V_\s:=H^q_c(\nu^{-1}(|\s|^*),\Lambda)$ satisfait la propri\'et\'e suivante: si $\s$ est un simplexe contenant un sommet $s,$ alors le morphisme $V_\s\to V_s$ induit un isomorphisme $V_\s\cong V_s^{G_\s^+},$ o\`u $G_\s^+$ d\'esigne le pro-$p$-radical du fixateur de $\s.$

Nous \'etudions les syst\`emes de coefficients ayant la propri\'et\'e sus-mentionn\'ee dans la section 2. Nous d\'emontrons le r\'esultat suivant qui implique entre autres la d\'eg\'en\'erescence en $E_1$ de la suite spectrale \ref{Intro1}.

\medskip

\noindent{\bf Th\'eor\`eme B.} (Th\'eor\`eme \ref{thm2})
{\it Soit $\{V_\s\}$ un syst\`eme de coefficients tel que le morphisme $\varphi^\s_s:V_\s\to V_s$ induit un isomorphisme $V_\s\simto V_s^{G^+_\s}.$ Alors, le complexe de cha\^ines $\CC_*(\BC\TC,\{V_\s\})$ associ\'e \`a $\{V_\s\}$ est acyclique sauf en degr\'e z\'ero, donc il induit une r\'esolution de $H_0(\CC_*(\BC\TC,\{V_\s\})).$ De plus, $V_\s$ est isomorphe \`a $H_0(\CC_*(\BC\TC,\{V_\s\}))^{G^+_\s}$ canoniquement.}

\medskip

En fait, on d\'emontre ce r\'esultat dans un cadre un peu plus g\'en\'eral, pour les syst\`emes de coefficients associ\'es \`a un syst\`eme d'idempotents $(e_s)_{s\in\BC\TC_0}$ dans le langage de Meyer et Solleveld \cite{MS}.

\medskip

D\'ecrivons bri\`evement le contenu des diff\'erentes parties. Au paragraphe 2.1, apr\`es quelques rappels sur les syst\`emes de coefficients, on pr\'ecise dans quel cadre on utilise le langage de \cite{MS}. Ensuite, on d\'emontre le th\'eor\`eme B dans les paragraphes 2.2 - 2.4. Le th\'eor\`eme A est obtenu dans la section 3. Au paragraphe 3.1, on fait un rappel sur les r\'esultats g\'eom\'etriques obtenus dans \cite{Wang-Sigma1}. Ensuite, on donne l'admissibilit\'e et l'absence de torsion pour les coefficients entiers. Celles-ci d\'ecoulent des r\'esultats connus \cite{bonnafe-rouquier-coxeter} sur les vari\'et\'es de Deligne-Lusztig. Aux paragraphes 3.3 et 3.4, on fait des rappels et des compl\'ements sur les repr\'esentations elliptiques de groupes finis et de groupes $p$-adiques, et on calcul la cohomologie de vari\'et\'e de Deligne-Lusztig Coxeter associ\'ee \`a  $\GL_d.$ On \'etudie la partie non-supercuspidale de la cohomologie $\ell$-adique de $\MC_{Dr,1}/\varpi^\ZM$ dans 3.5. Les ingr\'edients cruciales sont le foncteur de l'induction \guillemotleft ~tordue~\guillemotright~ de Lusztig et la version explicite de la correspondance de Jacquet-Langlands d\'ecrite par Bushnell et Henniart \cite{Bushnell-Henniart-level0} .

\medskip
\textbf{Remerciements:} Je remercie profond\'ement mon directeur de th\`ese Jean-Fran\c cois Dat tant pour son aide g\'en\'ereuse que pour ses constants encouragements. Je tiens \`a exprimer ma gratitude \`a C\'edric Bonnaf\'e pour r\'epondre \`a mes questions. Une partie de cet article a \'et\'e accomplie lorsque je visitais Institute for Mathematical Sciences, National University of Singapore en 2013. Je le remercie pour l'excellente condition de travail.

\section{Syst\`emes de coefficients sur l'immeuble de Bruhat-Tits}
Soit $K$ une extension finie de
$\QM_p$ d'anneau des entiers $\OC.$ Notre but de cette section est de montrer le th\'eor\`eme B dans l'introduction. Pour les d\'efinitions et les propri\'et\'es fondamentaux de l'immeuble de Bruhat-Tits semi-simple $\BC\TC$ associ\'e \`a $G:=\GL_d(K)$ ($d\geq 2$), on renvoie les lecteurs \`a \cite{fargues} annexe A. On s'int\'eresse \`a la cat\'egorie des $R[G]$-modules {\em lisses} \`a gauche, o\`u $R$ est une alg\`ebre sur $\ZM[\frac{1}{p}]$ ({\em cf.} \cite{vigneras-invent}).

\subsection{Rappels et compl\'ements sur les syst\`emes de coefficients}\label{h}

\ali \label{Sec2P1} Rappelons que $\BC\TC$ est un complexe simplicial partiellement ordonn\'e de sorte que $\tau<\s$ si $\tau$ est une facette de $\s.$ Un simplexe de dimension $0$ est appel\'e un {\em sommet}, et un simplexe de dimension maximale est appel\'e une {\em chambre}. On sait que l'ensemble des sommets de $\BC\TC$ s'identifie \`a l'ensemble des classes d'homoth\'etie de $\OC$-r\'eseaux dans l'espace vectoriel $K^d$ (un $\OC$-reseau de $K^d$ est un sous-$\OC$-module $M$ de $K^d$ tel que $M\otimes_{\OC}K\cong K^d$). Un ensemble de simplexes  $\s_1,\ldots,\s_k$ est appel\'e {\em adjacent} s'il existe un simplexe $\s$ tel que $\s_i<\s$ $\forall i\in\{1,\ldots,k\}.$ Si $\s$ et $\tau$ deux simplexes adjacents, on notera $[\s,\tau]$ le plus petit simplexe contenant $\s\cup\tau.$ Pour un sous complexe $\Sigma$ de $\BC\TC,$ on d\'esignera $\Sigma^\circ$ l'ensemble des sommets de $\Sig.$ L'action de $G$ sur $\BC\TC$ se factorise par $\PGL_d(K).$ Soient $\tau,\s$ deux simplexes, notons $H(\s,\tau)$ l'enclos de $\tau$ et $\s$: l'intersection de tous les appartements contenant $\s\cup\tau.$ Pour un simplexe $\s\in\BC\TC,$ notons $G_\s$ son fixateur sous l'action de $G$, i.e.
$$
G_\s=\{g\in G~|~ gx=x\text{ pour tout sommet $x$ de $\s$}\}.
$$
Il admet un unique pro-$p$-sous-groupe distingu\'e maximal, appel\'e son {\em pro-$p$-radical} et not\'e $G_\s^+.$ Plus g\'en\'eralement, fixons un entier $r\geq1$ et consid\'erons
$$
U^{(r)}:=\{g\in G~|~g\equiv 1 \pmod{\varpi^r}\}
$$
le sous-groupe congruence principal de niveau $r$ dans $G.$ Pour tout sommet $x$ de $\BC\TC,$ on peut consid\'erer le sous-groupe de $G^+_x$
$$
U^{(r)}_x:=gU^{(r)}g^{-1},\text{ si $x=g([\OC^d])$ avec $g\in G$}.
$$
Pour un simplexe quelconque $\s\in\BC\TC,$ on peut consid\'erer le groupe ({\em cf.} \cite{SS-crelle})
$$
U^{(r)}_\s:=\text{le sous-groupe de $G$ engendr\'e par $U^{(r)}_x,$ avec $x\in\s.$ }
$$
Schneider et Stuhler ont d\'emontr\'e que lorsque $r=1,$ $U^{(1)}_\s=G^+_\s,~\forall\s\in\BC\TC$ ({\em cf.} \cite[\S 6 Lemme 2]{SS-invent}).

\begin{DEf}\label{c}
(1) {\em Un syst\`eme de coefficients} $\G$ sur $\BC\TC$ \`a valeurs dans la cat\'egorie des $R$-modules est un foncteur contravariant de la cat\'egorie $(\BC\TC,<)$ vers la cat\'egorie de $R$-modules. Concr\`etement, c'est la donn\'ee des $R$-modules $(V_\s)_{\s\in\BC\TC}$ et des $R$-morphismes $\varphi^\s _\tau:V_\s\To{}V_\tau$ si $\tau<\s,$ soumis aux conditions: $\varphi^\s_\s=\id_{V_\s}$ et $\varphi^{\tau}_\omega\circ \varphi^\s_\tau=\varphi^\s_\omega$ si $\omega<\tau<\s.$

(2) Un {\em syst\`eme de coefficients $G$-\'equivariant} est un syst\`eme de coefficients $\G$ muni des isomorphismes $\a_g:V_\s\simto V_{g\cdot\s},$ $\forall g\in G,$ $\s\in\BC\TC$ compatibles avec les $\varphi^\s_\tau,$ et $\a_1=\id,$ $\a_g\circ\a_h=\a_{gh}.$

(3) Un syst\`eme de coefficients $G$-\'equivariant est dit {\em de niveau 0}, si pour tout sommet $x$ de $\BC\TC$ et tout simplexe $\s$ contenant $x,$ le morphisme $\varphi^\s_x:V_\s\to V_x$ induit un isomorphisme $V_\s\simto V_x^{G^+_\s}.$ On note $\Coef_0(G,R)$ la cat\'egorie ab\'elienne des syst\`emes de coefficients $G$-\'equivariants de niveau $0,$ dont les morphismes sont des morphismes $G$-\'equivariants de syst\`emes de coefficients.

\end{DEf}

\begin{rem}
On peut consid\'erer les syst\`emes de coefficients $G^\circ$-\'equivariants de niveau 0, o\`u $G^\circ=\Ker(\val_K\circ\det:\GL_d(K)\to K^\times),$ car pour tout $\s\in\BC\TC,$ $G^+_\s\subset G^\circ.$ Ce sont les syst\`emes de coefficients satisfaisant la condition (2) de \ref{c} pour tout $g\in G^\circ$ et la condition (3). Le syst\`eme de coefficients $\s\mapsto H^q_c(\nu^{-1}(|\s|^*),\Lambda)$ que nous allons consid\'erer dans \S 3 est un syst\`eme de coefficients $G^\circ$-\'equivariant de niveau 0.
\end{rem}

\begin{exe}
L'exemple essentiel de syst\`eme de coefficients de $R$-modules introduit par
Schneider et Stuhler est le suivant: partant alors d'un $R G$-module lisse
V, on d\'efinit un syst\`eme de coefficients \`a valeurs dans la cat\'egorie des $R$-modules :
\begin{itemize}
  \item $\s\mapsto V^{U^{(r)}_\s}.$
  \item $(\tau<\s)\mapsto (V^{U^{(r)}_\s}\into V^{U^{(r)}_\tau})$ dont l'existence est assur\'ee par l'inclusion $U^{(r)}_\tau\subset U^{(r)}_\s.$
\end{itemize}
Bien s\^ur, le syst\`eme de coefficients ainsi obtenu est $G$-\'equivariant.
\end{exe}

\ali\label{e} Plus g\'en\'eralement, on peut consid\'erer la cat\'egorie ab\'elienne des syst\`emes de coefficients associ\'es \`a un syst\`eme d'idempotents $e=(e_\s)_{\s\in\BC\TC}$ satisfaisant les conditions de \cite[Def. 2.1]{MS}. Plus pr\'ecis\'ement, pour chaque sommet $x\in\BC\TC,$ $e_x$ est un idempotent dans l'alg\`ebre de Hecke $\HC(G_x,R)\into \HC(G,R).$ On consid\`ere l'ensemble des syst\`emes d'idempotents $e=(e_x)_{x\in\BC\TC^\circ}$ satisfaisant les trois propri\'et\'es suivantes:

\begin{itemize}
  \item (a) $e_x$ et $e_y$ commutent pour deux sommets adjacents $x,y$.
  \item (b) Soient $x,y,z$ trois sommets avec $z\in H(x,y)$ et $z,x$ adjacents, on a $e_xe_ze_y=e_xe_y.$
  \item (c) $e_{gx}=ge_x g^{-1},~\forall g\in G,~x\in\BC\TC^\circ.$
\end{itemize}
Pour un simplexe $\s\in\BC\TC,$ notons
$$
e_\s:=\mathop{\prod_{x\in\BC\TC^\circ}}_{x\in\s}e_x
$$
qui est un idempotent dans $\HC(G,R)$ bien d\'efini d'apr\`es (a). Pour la suite, on demande de plus que le syst\`eme d'idempotents satisfait une condition suppl\'ementaire:
\begin{itemize}
\item (d) Pour tout sommet $x$ et tout simplexe $\s$ contenant $x,$ $e_\s\in\HC(G_x,R).$
\end{itemize}

\begin{rem}
(i) Soient $e=(e_\s)_{\s\in\BC\TC}$ un syst\`eme d'idempotents satisfaisant les conditions (a)-(c) de \ref{e}, et $V$ un $RG$-module lisse. Alors, le foncteur $\s\mapsto e_\s(V)$ d\'efinit un syst\`eme de coefficients $G$-\'equivariant ({\em cf.} \cite{MS}).

(ii) Notons $e_{U^{(r)}_\s}:=\frac{\chi_{U^{(r)}_\s}}{\mu(U^{(r)}_\s)}$ l'idempotent associ\'e au sous-groupe compact ouvert $U^{(r)}_\s.$ Le syst\`eme d'idempotents $(e_{U^{(r)}_\s})_{\s\in\BC\TC}$ satisfait les conditions (a)-(d), voir \cite[Section 2.2]{MS}.

(iii) Supposons que l'on se donne un syst\`eme d'idempotents $(e_\s)_{\s\in\BC\TC}$ satisfaisant les conditions (a)-(c). S'il existe un entier $r\geq1$ tel que pour tout sommet $x$ et tout simplexe $\s$ contenant $x$, on ait
$$
e_\s=e_{U^{(r)}_\s}e_x,
$$
alors $(e_\s)_{\s\in\BC\TC}$ satisfait la condition (d).
\end{rem}

\begin{DEf}
Soit $e=(e_\s)_{\s\in\BC\TC}$ un syst\`eme d'idempotents satisfaisant les conditions (a)-(d) de \ref{e}. {\em Un $e$-syst\`eme de coefficients} est un syst\`eme de coefficients tel que pour tout $x\in\BC\TC^\circ$ et tout simplexe $\s$ contenant $x,$ le morphisme $\varphi^\s_x:V_\s\to V_x$ induit un isomorphisme $V_\s\simto e_\s(V_x).$ Notons que la condition (d) assure que $e_\s\in\HC(G_x,R)$ agit sur $V_x.$ On a comme pr\'ec\'edemment une notion de $e$-syst\`eme de coefficients $G$-\'equivariant, et on note $\Coef_e(G,R)$ la cat\'egorie ab\'elienne de $e$-syst\`emes de coefficients $G$-\'equivariants, dont les morphismes sont des morphismes $G$-\'equivariants de syst\`emes de coefficients.
\end{DEf}

\ali Pour former le complexe de cha\^ine cellulaire associ\'e \`a un sous complexe $\Sig$, on munit chaque simplexe $\s$ d'une orientation ({\em cf.} \cite[(1.1.2)]{MS}) qui induit une orientation sur chaque sous facette de $\s.$ On d\'efinit
$$
\e_{\tau\s}:=\left\{
  \begin{array}{ll}
    1, & \hbox{si $\tau<\s$ et l'orientation sur $\tau$ co\"incide avec celle induite par $\s$ ;} \\
    -1, & \hbox{si $\tau<\s$ et l'orientation sur $\tau$ est oppos\'ee \`a celle induite par $\s$;} \\
    0, & \hbox{si $\tau$ n'est pas une facette de $\s$.}
  \end{array}
\right.
$$


Soient $\G=(V_\s,\varphi^\s_\tau)$ un syst\`eme de coefficients et $\Sig$ un sous complexe de $\BC\TC,$ on les associe un {\em complexe de cha\^ines} sur $\Sig$ \`a coefficients $\G:$
\ini\begin{equation}\label{complexe}
\CC_*(\Sig,\G):=\CC^{or}_c(\Sig_{d-1},\G)\To{\partial}\cdots\To{\partial}\CC^{or}_c(\Sig_{0},\G)
\end{equation}
o\`u
$\CC^{or}_c(\Sig_{q},\G):=\bigoplus_{\s\in\Sig,~\dim\s=q}V_\s,$ et la diff\'erentielle est donn\'ee par
$$
\partial\big((v_\s)_{\s\in\Sig}\big)_\tau:=\sum_{\tau\in\Sig}\e_{\tau\s} \varphi^\s_\tau(v_\s).
$$
Nous noterons $H_*(\Sig,\G)$ les objets d'homologie du complexe $\CC_*(\Sig,\G).$ Notons que $G$ agit sur les facettes orient\'ees de $\BC\TC_{q}$ de sorte que si $\G$ est un syst\`eme de coefficients $G$-\'equivariant, alors pour tout $q\geq 0,$ $\CC^{or}_c(\BC\TC_{q},\G)$ poss\`ede une action de $G.$

\ali Soit $e=(e_\s)_{\s\in\BC\TC}$ un syst\`eme d'idempotents satisfaisant les conditions (a)-(c) de \ref{e}. Notons $\Rep_e(G,R)$ la cat\'egorie de $\HC(G,R)$-modules non-d\'eg\'en\'er\'es $V$ tels que
$$
V=\sum_{x\in\BC\TC^\circ}e_x(V).
$$
Gr\^ace \`a \cite[Thm. 3.1]{MS}, $\Rep_e(G,R)$ est une sous-cat\'egorie de Serre de la cat\'egorie de $R[G]$-modules lisses. Lorsque $e=(e_{U^{(r)}_\s})_{\s\in\BC\TC},$ on dit que c'est la cat\'egorie de $R[G]$-modules lisses de niveau $r,$ not\'ee par $\Rep_r(G,R).$ Rappelons le th\'eor\`eme principal de \cite{MS}:

\begin{theo}\textsl{(\cite[Thm. 2.4]{MS}) }\label{f}
Sous ces hypoth\`eses, soient $V$ un $R$-module tel que $e_\s\in\End(V)~\forall \s\in\BC\TC,$ et $\G(V)$ le syst\`eme de coefficients $(\s\mapsto e_\s(V))$ associ\'e \`a $V.$ Alors le complexe de cha\^ines
\begin{align*}
\CC^{or}_c(\BC\TC_{d-1},V)\To{\partial}\cdots\To{\partial}\CC^{or}_c(\BC\TC_{0},V)\To{}&V\To{}0\\
(v_x)_{x\in\BC\TC^\circ}\longmapsto &\sum_{x\in\BC\TC^\circ}v_x
\end{align*}
est une r\'esolution de $V.$
\end{theo}

\begin{rem}
(i) Lorsque $e=(e_{U^{(r)}_\s})_{\s\in\BC\TC}$ et $V\in \Rep_r(G,R),$ ce th\'eor\`eme est d\'emontr\'e par Schneider et Stuhler \cite{SS-crelle}.

(ii) Ce th\'eor\`eme est valable pour tout groupe r\'eductif $p$-adique, {\em cf.} \cite{SS-ihes}, \cite{MS}.
\end{rem}
Notre but est de d\'emontrer le r\'esultat suivant:

\begin{theo}\label{thm2}
Soient $e=(e_\s)_{\s\in\BC\TC}$ un syst\`eme d'idempotents satisfaisant les conditions (a)-(d) de \ref{e}, et $\G$ un $e$-syst\`eme de coefficients sur $\BC\TC.$ Alors,
\begin{description}
\item [(a)] en posant $\Sig=\BC\TC,$ le complexe de cha\^ines \ref{complexe} est exact sauf en degr\'e $0,$ i.e. c'est une r\'esolution de $H_0(\BC\TC,\G).$

\item [(b)] $\G$ est isomorphe au syst\`eme de coefficients $(\s\mapsto e_\s(H_0(\BC\TC,\G))).$
\end{description}
\end{theo}

Si l'on consid\`ere la cat\'egorie des $e$-syst\`emes de coefficients $G$-\'equivariants, on a le corollaire suivant:
\begin{coro}
Soit $e=(e_\s)_{\s\in\BC\TC}$ un syst\`eme d'idempotents satisfaisant les conditions (a)-(d) de \ref{e}, le foncteur
\begin{align*}
\Rep_e(G,R)&\To{}\Coef_e(G,R)\\
V&\longmapsto \G(V)
\end{align*}
admet un quasi-inverse $\G\mapsto H_0(\BC\TC,\G),$ donc induit une \'equivalence de cat\'egories.
\end{coro}

\begin{preuve}
Si $\G$ est un syst\`eme de coefficient $G$-\'equivariant, on sait que $H_0(\BC\TC,\G)$ est muni une action de $G.$ Gr\^ace au (b) du th\'eor\`eme pr\'ec\'edent, on sait que le foncteur $V\mapsto\G(V)$ est essentiellement surjectif. Donc il suffit de montrer que le foncteur est pleinement fid\`ele. Soient $V,W\in\Rep_e(G,R),$ $\G$ induit un morphisme
\begin{align*}
\G:\Hom_{RG}(V,W)&\To{}\Hom_{\Coef_e}(\G(V),\G(W))\\
f&\longmapsto (f_\s=f|_{e_\s(V)}:V_{\s}\to W_\s).
\end{align*}
Soit $f\in\Hom_{RG}(V,W)$ tel que $\G(f)=0.$ Alors, $\forall x\in\BC\TC^{\circ},$ $f|_{e_x(V)}=f_x=0.$ Par hypoth\`ese que $V\in\Rep_e(G,R),$ donc $V=\sum_{x\in\BC\TC^\circ}e_x(V).$ On en d\'eduit que $f=0.$ Ceci d\'emontre l'injectivit\'e de $\G$. Pour la surjectivit\'e, soit $(g_\s)_{\s\in\BC\TC}\in\Hom_{\Coef_e}(\G(V),\G(W)),$ il induit un morphisme de complexes
$$
g:\CC_*(\BC\TC,\G(V))\To{}\CC_*(\BC\TC,\G(W)).
$$
D'apr\`es le th\'eor\`eme \ref{f}, le complexe
$$
\CC_*(\BC\TC,\G(V))\to V\to 0 \      \ \text{ (resp. $\CC_*(\BC\TC,\G(W))\to W\to 0$)}
$$
est une r\'esolution de $V$ (resp. $W$). Donc $H_0(g)$ induit un morphisme de $RG$-modules
$$
H_0(g):V\To{}W.
$$
Par d\'efinition, pour tout sommet $x\in\BC\TC,$ on a $H_0(g)|_{e_x(V)}=g_x.$ Donc pour tout simplexe $\s$ contenant $x,$ on a $$H_0(g)|_{e_\s(V)}=(H_0(g)|_{e_x(V)})|_{e_\s(V)}=g_x|_{e_\s(V)}=g_\s.$$
Ceci d\'emontre la surjectivit\'e.
\end{preuve}

\subsection{Les applications locales}

D\'esormais, on utilisera les alphabets latins $x,y,z\ldots$ pour des sommets de $\BC\TC,$ et les alphabets grecs $\s,\tau,\omega,\ldots$ pour des simplexes quelconques.

\ali\label{ff} Pour un $e$-syst\`eme de coefficients $\G=(V_\s)_{\s\in\BC\TC},$ on identifie $V_\s$ et $e_\s(V_x)$ via le morphisme canonique $\varphi^{\s}_x$ lorsque $x\in\s.$ Pour deux simplexes $\tau<\s,$ $V_\s=e_\s(V_\tau),$ notons alors $p^\tau_\s$ le projecteur $V_\tau\onto e_\s(V_\tau)=V_\s.$ Nous allons d\'efinir une famille d'applications $\e^\s_x:V_\s\to V_x$ ($x$ n'est pas forc\'ement contenu dans $\s$) telle que $\forall \tau<\s,$ l'application $\e^\s_x$ se factorise par $\e^\tau_x,$ i.e. le diagramme suivant soit commutatif:
$$
\xymatrix{
V_\s \ar[r]^{\e^\s_x} \ar@{^(->}[d]^{\ph^\s_\tau} & V_x\\
V_\tau \ar[ru]_{\e^\tau_x}
}
$$

Tout d'abord, soient $x,\s$ adjacents, on d\'efinit l'application $\e^\s_x$ par la compos\'ee:
$$
\e^\s_x:V_\s\To{p^\s_{[x,\s]}}V_{[x,\s]}=e_{[x,\s]}(V_x)\into V_x
$$
En particulier, si $x\in\s,$ $\e^\s_x$ est l'inclusion $V_\s=e_\s(V_x)\into V_x$ induite par $\varphi^\s_x.$ On peut aussi d\'efinir un idempotent $\wt{\e^\s_x}\in\End(V_\s):$
$$
\wt{\e^\s_x}:V_\s \To{p^\s_{[x,\s]}}  V_{[x,\s]}=e_{[x,\s]}(V_\s) \into  V_\s.
$$

\begin{lemme}\label{lem1}
Soit $y$ un autre sommet tel que $x,y,\s$ soient adjacents, alors $\wt{\e^\s_y}\circ\wt{\e^\s_x}=\wt{\e^\s_x}\circ\wt{\e^\s_y}.$
\end{lemme}
\begin{preuve}
L'application $\wt{\e^\s_y}\circ\wt{\e^\s_x}$ est la compos\'ee
$$
V_\s\onto e_{[y,\s]}(e_{[x,\s]}(V_\s))\into V_\s
$$
et l'application $\wt{\e^\s_x}\circ\wt{\e^\s_y}$ est la compos\'ee
$$
V_\s\onto e_{[x,\s]}(e_{[y,\s]}(V_\s))\into V_\s.
$$
Gr\^ace \`a l'hypoth\`ese de $(e_\s)_{\s\in\BC\TC}$, on a $e_{[x,\s]}e_{[y,\s]}=e_{[y,\s]}e_{[x,\s]},$ d'o\`u l'\'enonc\'e du lemme.
\end{preuve}

Soit maintenant $y$ un sommet quelconque. On peut choisir un appartement $A$ contenant $x$ et $y,$ et une chambre de Weyl ferm\'ee \footnote{$C$ est la cl\^oture d'une chambre de Weyl.} $C$ telle que $y\in x+C.$
\begin{lemme}\label{k}
L'enclos $H(x,y)$ est \'egal \`a $(x+C)\cap (y-C).$
\end{lemme}
\begin{preuve}
Rappelons que $H(x,y)$ est l'intersection de tous les appartements contenant $x,y.$ Il est aussi \'egal \`a la r\'ealisation g\'eom\'etrique de l'enveloppe convexe simpliciale de $\{x,y\},$ {\em cf.} \cite[Annexe A. Prop. 29]{fargues}. D'abord on sait que les ensembles $x+C$ et $y-C$ sont convexes, leur intersection $(x+C)\cap (y-C)$ est \'egalement convexe et contient $x,y.$ Donc $H(x,y)\subset(x+C)\cap (y-C).$ R\'eciproquement, il suffit de montrer que pour tout sommet $z\in(x+C)\cap (y-C),$ $z$ est contenu dans $H(x,y).$

Comme tout sommet de l'immeuble de Bruhat-Tits semi-simple associ\'e \`a $\GL_d(K)$ est un point sp\'ecial \cite{BT}, on peut alors supposer que $x$ est l'origine de la donn\'ee radicielle valu\'ee, i.e. $x$ est repr\'esentable par le r\'eseau $[\OC^d]$ de l'espace vectoriel $K^d.$ Notons $\D=\{\a_1,\ldots,\a_{d-1}\}$ l'ensemble de racines simples associ\'ees \`a $C,$ on a alors $\a_i(x)=0,~\forall i\in\{1,\ldots,d-1\}.$ Par hypoth\`ese, on sait que pour tout $i\in\{1,\ldots,d-1\},$ $\a_i(y)\geq \a_i(z)\geq0.$ Pour montrer que $z\in H(x,y),$ il suffit de montrer que pour toute racine affine $\a$ telle que $\a(x)\geq0$ et $\a(y)\geq0,$ on a $\a(z)\geq0.$ \'Ecrivons $\a=\sum\l_i\a_i+n$ o\`u $n\in\ZM$ et les $\l_i$ sont tous positifs ou tous n\'egatifs. Notons d'abord que $n\geq 0,$ car $\a(x)\geq0$ et $\a_i(x)=0.$ Si $\l_i\geq 0,$ on a $\a(z)=\sum\l_i\a_i(z)+n\geq0.$ Si $\l_i\leq0,$ alors puisque $\a(y)=\sum\l_i\a_i(y)+n\geq0,$ on a aussi
$$
\a(z)=\sum\l_i\a_i(z)+n=\sum\l_i\a_i(y)+n+\sum(-\l_i)(\a_i(y)-\a_i(z))\geq0.
$$
\end{preuve}

\begin{DEf}
On dit qu'une suite finie de sommets diff\'erents $(z_0,z_1,\ldots,z_{m-1},z_m)$ est un {\em chemin tendu} de $y$ vers $x,$ si $z_0=y,$ $z_m=x$ et $z_{i-1},z_i$ sont adjacents de sorte que $z_i\in H(x,z_{i-1})~\forall i\in\{1,\ldots,m\}.$
\end{DEf}

\begin{lemme}\label{m}

\begin{description}
\item [(a)] Si $(z_0,z_1,\ldots,z_{m-1},z_m)$ est un chemin tendu, alors $(z_m,z_{m-1},\ldots,z_0)$ l'est aussi.

\item [(b)] L'entier $m$ dans la d\'efinition ci-dessus ne d\'epend pas du choix de chemin tendu. On l'appelle la distance entre $x$ et $y,$ not\'ee par $\rho(x,y).$

\item [(c)] Soit $(z_0,z_1,\ldots,z_{m-1},z_m)$ un chemin tendu, alors $\forall 0\leq k\leq m,$ $(z_0,\ldots,z_k)$ (resp. $(z_k,\ldots,z_m)$) est un chemin tendu.

\item [(d)] Soient $(z_0,\ldots,z_k)$ et $(z_k,\ldots,z_m)$ deux chemins tendus et $z_k\in H(z_0,z_m),$ alors $(z_0,\ldots,z_k,\ldots,z_m)$ est un chemin tendu.
\end{description}
\end{lemme}
\begin{preuve}
(a) Choisissons un appartement $A$ contenant $z_0$ et $z_m,$ et une chambre de Weyl ferm\'ee $C$ telle que $z_0\in z_m+C.$ Par d\'efinition, on a $z_i\in H(z_{i-1},z_m)=(z_m+C)\cap(z_{i-1}-C)$ {\em cf.} \ref{k}. Ceci entra\^ine que $z_{i-1}\in H(z_0,z_i),$ car $z_{i-1}\in z_0-C.$ Autrement dit, $(z_m,z_{m-1},\ldots,z_0)$ est un chemin tendu.

(b) Le choix d'un appartement $A$ comme ci-dessus nous fournit un tore maximal d\'eploy\'e $T$ dans $G$ tel que la r\'ealisation g\'eom\'etrique $|A|$ de $A$ s'identifie \`a l'espace euclidean $X_*(T)\otimes \RM/(X_*(Z(G))\otimes\RM),$ o\`u $Z(G)$ d\'esigne le centre de $G.$ On peut supposer que $T$ s'identifie au sous-groupe de matrices diagonales. On peut supposer de plus que les racines simples $\D=\{\a_1,\ldots,\a_{d-1}\}$ associ\'ees \`a la chambre de Weyl ferm\'ee $C$ que l'on a choisie sont de la forme $\a_i(t)=t_i/t_{i+1},$ o\`u $t=\diag(t_1,\ldots,t_{d}).$ Notons $\{\omega_1,\ldots,\omega_{d-1}\}$ l'ensemble des copoids fondamentaux associ\'es \`a $\D,$ o\`u $\omega_i(a)=\diag(\underbrace{a,\ldots,a}_i,1\ldots,1),~\forall a\in K^\times.$ Supposons de plus que le sommet $x$ est l'origine de cet espace affine $X_*(T)\otimes \RM/(X_*(Z(G))\otimes\RM)\cong \RM^{d-1}.$ Consid\'erons les cocaract\`eres $\l_k:K^\times \to T$ tels que $(\l_k(a))_{ij}=a$ si $i=j=k$ et $(\l_k(a))_{ij}=1$ si $i=j\neq k.$ Alors $\l_1,\ldots,\l_{d-1}$ induisent une base de $|A|\cong \RM^{d-1},$ sous laquelle $\omega_i$ s'identifie au vecteur $\omega_i:=(\underbrace{1,\ldots,1}_i,0,\ldots,0)\in\RM^{d-1},$ et la chambre vectorielle $x+C$ s'identifie au c\^one $\RM_+\omega_1+\cdots+\RM_+\omega_{d-1}\subset\RM^{d-1}.$ Ainsi le sommet $y$ s'identifie \`a un point $(a_1,\ldots,a_{d-1})\in\ZM^{d-1},$ o\`u $a_1\geq \cdots\geq a_{d-1}\geq 0,$ et la chambre vectorielle $y-C$ s'identifie au c\^one $(a_1,\ldots,a_{d-1})+\RM_{-}\omega_1+\cdots+\RM_{-}\omega_{d-1}\subset\RM^{d-1}.$ Un sommet $z\in H(x,y)$ adjacent \`a $y$ s'identifie alors \`a un point $(b_1,\ldots,b_{d-1})\in\ZM^{d-1}$ satisfaisant $b_1\geq\cdots\geq b_{d-1}\geq0$ et $1\geq a_1-b_1\geq\cdots\geq a_{d-1}-b_{d-1}\geq 0.$  Donc on voit que le nombre $m$ est \'egal \`a $a_1,$ ind\'ependant du choix de chemin tendu.

(c) Par d\'efinition, on sait que $(z_k,z_{k+1},\ldots,z_m)$ est un chemin tendu. D'apr\`es le lemme \ref{k}, $z_k\in H(z_{k-1},z_m)=(z_m+C)\cap (z_{k-1}-C),$ on a $z_{k-1}\in z_k+C.$ On en d\'eduit par r\'ecurrence que $z_i\in z_k+C, ~\forall 0\leq i\leq k-1.$ De plus, notons que $z_i\in H(z_{i-1},z_m)=(z_m+C)\cap (z_{i-1}-C),$ on a alors que $z_i\in z_{i-1}-C.$ Donc $z_i\in (z_k+C)\cap(z_{i-1}-C)=H(z_{i-1},z_k), ~\forall 0\leq i\leq k-1.$ C'est-\`a-dire $(z_0,z_1,\ldots,z_k)$ est un chemin tendu.

(d) Il s'agit de montrer que pour tout $i\in\{1,\ldots,k\},$ $z_i\in H(z_{i-1},z_m).$ Par hypoth\`ese, $z_k\in (z_m+C)\cap (z_0-C),$ donc $z_0\in z_k+C.$ Pour chaque $i,$ $z_i\in H(z_{i-1},z_k)=(z_k+C)\cap(z_{i-1}-C),$ donc $z_i\in z_k+C\subset z_m+C.$ On en d\'eduit que $z_i\in (z_{i-1}-C)\cap (z_m+C)=H(z_{i-1},z_m).$
\end{preuve}


La preuve du lemme pr\'ec\'edent donne \'egalement l'existence d'un chemin tendu de $y$ vers $x,$ on d\'efinit alors $\e^y_x$ comme une composition de $\e^{z_{i-1}}_{z_i}$:
$$
\e^y_x:=\e^{z_{m-1}}_x\circ\cdots\circ \e^y_{z_1}.
$$

\begin{prop}\label{l}
Cette d\'efinition ne d\'epend pas du choix de chemin tendu.
\end{prop}

\begin{preuve}
On d\'emontre cette proposition par r\'ecurrence sur $m=\rho(x,y).$ Lorsque $m=1,$ on n'a qu'un seul chemin tendu reliant $x$ et $y,$ l'\'enonc\'e est trivial. Si $m>1,$ soit $(z'_0=y,\ldots,z'_m=x)$ un autre chemin tendu de $y$ vers $x.$ On peut supposer que $z'_i\neq z_i~\forall i\in\{1,\ldots,m-1\}.$ Sinon, il existe un entier $k\in\{1,\ldots,m-1\}$ tel que $z_k=z'_k.$ Par le lemme \ref{m} (c), le chemin $(y=z_0,z_1,\ldots,z_k)$ (resp. $(y=z'_0,z'_1,\ldots,z'_k=z_k)$) est un chemin tendu de $y$ vers $z_k,$ et le chemin $(z_k,z_{k+1},\ldots,z_m=x)$ (resp. $(z_k=z'_k,z'_{k+1},\ldots,z'_m=x)$) est un chemin tendu de $y$ vers $z_k.$ Par r\'ecurrence, on a
$$
\e^{z_{m-1}}_x\circ\cdots\circ \e^{z_k}_{z_{k+1}}=\e^{z'_{m-1}}_x\circ\cdots\circ \e^{z_k}_{z'_{k+1}}\   \ \text{ et } \   \ \e^{z_{k-1}}_{z_k}\circ \cdots\circ \e^{y}_{z_{1}}=\e^{z'_{k-1}}_{z_k}\circ\cdots\circ \e^{y}_{z'_{1}}.
$$
On en d\'eduit que
\begin{align*}
\e^{z_{m-1}}_x\circ\cdots\circ \e^y_{z_1}&=\e^{z_{m-1}}_x\circ\cdots\circ \e^{z_k}_{z_{k+1}}\circ\e^{z_{k-1}}_{z_k}\circ \cdots\circ \e^{y}_{z_{1}}\\
&=\e^{z'_{m-1}}_x\circ\cdots\circ \e^{z_k}_{z'_{k+1}}\circ\e^{z'_{k-1}}_{z_k}\circ\cdots\circ \e^{y}_{z'_{1}}=\e^{z'_{m-1}}_x\circ\cdots\circ \e^y_{z'_1}.
\end{align*}

Supposons maintenant que $z'_i\neq z_i~\forall i\in\{1,\ldots,m-1\}.$ D'apr\`es \cite[Lemme 2.9]{MS}, $y,z_1,z'_1$ sont adjacents. Si $m=\rho(x,y)=2,$ on a par la condition (b) dans \ref{e} que $e_xe_{z_1}e_y=e_xe_y=e_xe_{z'_1}e_y.$ On sait alors que le morphisme $\e^{z_1}_x\e^y_{z_1}$ est donn\'e par la compos\'ee:
\begin{align*}
V_y\To{p^y_{[y,z_1]}}V_{[y,z_1]}=e_{[y,z_1]}(V_{z_1})&\To{p^{z_1}_{[x,z_1]}}e_{[x,z_1]}e_{[y,z_1]}(V_{z_1})=e_xe_{[z_1,z'_1]}e_{[y,z_1]}(V_{z_1})=e_xe_{[z_1,z'_1]}(V_{[y,z_1,z'_1]})\\
&\into e_xe_{[z_1,z'_1]}(V_{[z_1,z'_1]}) = V_{[x,z_1,z'_1]}\into V_x.
\end{align*}
On voit qu'il s'identifie \`a la compos\'ee:
$$
V_y\To{p^y_{[y,z_1,z'_1]}}V_{[y,z_1,z'_1]}\into V_{[z_1,z'_1]}\To{p^{[z_1,z'_1]}_{[z_1,z'_1,x]}}V_{[x,z_1,z'_1]}\into V_x.
$$
De la m\^eme mani\`ere, on sait que ceci s'identifie \'egalement au morphisme $\e^{z'_1}_x\e^y_{z'_1}.$

Si $m=\rho(x,y)\geq 3,$ gr\^ace au lemme qui suit, il existe un sommet $z\in H(x,z_1)\cap H(x,z'_1)$ diff\'erent de $x$ tel que $z'_1,z_1\in H(z,y).$ Posons un chemin tendu de $z_1$ vers $z$ (resp. $z'_1$ vers $z$), alors $(y,z_1,\ldots,z)$ (resp. $(y,z'_1,\ldots,z)$) est un chemin tendu de $y$ vers $z,$ et $(z_1,\ldots,z,x)$ (resp. $(z'_1,\ldots,z,x)$) est un chemin tendu de $z_1$ (resp. $z'_1$) vers $x$ gr\^ace au lemme \ref{m} (d). En utilisant l'hypoth\`ese de r\'ecurrence, on a
\begin{align*}
\e^{z_{m-1}}_x\circ\cdots\circ \e^y_{z_1}&=\e^{z_1}_x\circ \e^y_{z_1}=\e^z_x\circ\e^{z_1}_z\circ \e^y_{z_1}=\e^z_x\circ \e^y_{z}\\
&=\e^z_x\circ\e^{z'_1}_z\circ \e^y_{z'_1}=\e^{z'_1}_x\circ \e^y_{z'_1}=\e^{z'_{m-1}}_x\circ\cdots\circ \e^y_{z'_1}.
\end{align*}
On conclut par r\'ecurrence.
\end{preuve}

\begin{lemme}\label{lem}
Soient $x,y,z$ trois sommets tels que $y,z$ adjacents et $\min\{\rho(x,y),\rho(x,z)\}\geq 2,$ alors il existe un sommet $w$ diff\'erent de $x$ tel que $w\in H(x,y)\cap H(x,z)$ et $x,w$ soient adjacents.
\end{lemme}
\begin{preuve}
Choisissons un appartement $A$ contenant $x$ et $[y,z],$ et une chambre de Weyl ferm\'ee $C$ tel que $[y,z]\subset x+C.$ Identifions $|A|$ \`a l'espace affine $\RM^{d-1}$ et $x$ au point d'origine comme dans la preuve du lemme \ref{m}. Comme $y,z$ sont adjacents, on peut supposer que $y$ (resp. $z$) s'identifie au point $(a_1,\ldots,a_{d-1})\in\ZM^{d-1}$ (resp. $(b_1,\ldots,b_{d-1})\in\ZM^{d-1}$) tel que $a_1\geq\cdots \geq a_{d-1}\geq0$ (resp. $b_1\geq\cdots \geq b_{d-1}\geq0$) et $|a_i-b_i|\leq 1,~\forall i.$ Notons $r\in\{1,\ldots,d-1\}$ (resp. $s\in\{1,\ldots,d-1\}$) l'entier tel que $a_r\geq2$ et $a_{r+1}\leq 1$ (resp. $b_s\geq2$ et $b_{s+1}\leq 1$). On peut alors prendre $w\in A$ le sommet associ\'e au point $(\underbrace{1,\ldots,1}_{\min\{r,s\}},0\ldots,0)\in \RM^{d-1}.$
\end{preuve}

Le corollaire suivant est un analogue \guillemotleft ~local~\guillemotright~de la condition (b) dans \cite[Def. 2.1]{MS}.
\begin{coro}\label{coro1}
Soient $x,y,z$ trois sommets tels que $z\in H(x,y),$ alors $\e^z_x\circ\e^y_z=\e^y_x.$
\end{coro}
\begin{preuve}
D'apr\`es le lemme \ref{m} (d), si l'on choisit un chemin tendu de $y$ vers $z,$ et un chemin tendu de $z$ vers $x,$ alors le chemin compos\'e est un chemin tendu de $y$ vers $x.$ L'\'enonc\'e du corollaire d\'ecoule du fait que $\e^y_x$ ne d\'epend pas du choix de chemin tendu, {\em cf.} \ref{l}.
\end{preuve}

\begin{lemme}\label{n}
Soient $x,y,z$ trois sommets et $y,z$ adjacents. Alors on a
\begin{description}
\item [(a)] $\e^y_x\circ\e^{[y,z]}_y=\e^z_x\circ\e^{[y,z]}_z.$

\item [(b)] $p^y_{[y,z]}\circ\e^x_y=p^z_{[y,z]}\circ\e^x_z.$

\end{description}
\end{lemme}
\begin{preuve}
On d\'emontre (a) par r\'ecurrence sur $m=\max\{\rho(x,y),\rho(x,z)\}.$ Lorsque $m=1,$ $x,y,z$ sont adjacents, on a par d\'efinition $\e^y_x\circ\e^{[y,z]}_y=\e^{[x,y,z]}_x\circ p^{[y,z]}_{[x,y,z]}=\e^z_x\circ\e^{[y,z]}_z.$ Si $m>1,$ d'apr\`es le lemme \ref{lem}, il existe un sommet $w\in H(x,y)\cap H(x,z)$ tel que $\rho(w,y)=\rho(x,y)-1$ et $\rho(w,z)=\rho(x,z)-1.$ Par r\'ecurrence, on a $\e^y_w\circ\e^{[y,z]}_y=\e^z_w\circ\e^{[y,z]}_z.$ En vertu du corollaire pr\'ec\'edent, on obtient que
$$
\e^y_x\circ\e^{[y,z]}_y=\e^w_x\circ\e^y_w\circ\e^{[y,z]}_y=\e^w_x\circ\e^z_w\circ\e^{[y,z]}_z=\e^z_x\circ\e^{[y,z]}_z.
$$
L'assertion (b) peut \^etre d\'emontr\'ee de la m\^eme mani\`ere, et on laisse le d\'etail au lecteur.
\end{preuve}

\ali Maintenant, pour un simplexe quelconque $\s,$ on d\'efinit l'application $\e^\s_x$ de la mani\`ere suivante: choisissons un sommet $y\in\s$ et posons
$$
\e^\s_x: V_\s\To{\e^\s_y} V_y \To{\e^y_x} V_x.
$$

\begin{lem}
Cette d\'efinition ne d\'epend pas du choix de $y.$ En cons\'equence, soient $\tau<\s$ deux simplexes, on a $\e^\s_x=\e^\tau_x\circ\varphi^\s_\tau$ (voir \ref{ff}).
\end{lem}
\begin{preuve}
Soit $y'$ un autre sommet de $\s.$ D'apr\`es le lemme \ref{n} (a), on a $$\e^y_x\circ\e^\s_y=\e^y_x\circ\e^{[y,y']}_y\circ\varphi^\s_{[y,y']}=\e^{y'}_x\circ\e^{[y,y']}_{y'}\circ\varphi^\s_{[y,y']}=\e^{y'}_x\circ\e^\s_{y'}.$$
\end{preuve}

\subsection{Les projecteurs $u^{\Sig}_{\Sig'}$}\label{j}

Soient $\Sig$ un sous complexe fini convexe de $\BC\TC,$ $\G=(V_\s)_{\s\in\BC\TC}$ un $e$-syst\`eme de coefficients. Fixons un sommet $x\in\Sig,$ et notons $\underline{V_x}:=(\s\mapsto V_x,\varphi^\s_x=\id_{V_x})$ le syst\`eme de coefficients constant \`a valeurs dans $V_x.$ Alors, les applications locales $\{\e^\s_x\}_{\s\in\Sig}$ induisent un morphisme de syst\`emes de coefficients, et donc un morphisme de complexes
$$
\oplus_{\s\in\Sig}\e^\s_x:\CC_*(\Sig,\G)\To{}\CC_*(\Sig,\underline{V_x}),
$$
et un morphisme sur les homologies $p^{\Sig}_x:=H_0(\oplus\e^\s_x):H_0(\Sig,\G)\to H_0(\Sig,\underline{V_x}).$
Soit $\Sig'\subset\Sig$ un sous complexe fini convexe, notons
$$\oplus_{\s\in\Sig'} \id_{V_\s}:\CC_*(\Sig',\G)\To{}\CC_*(\Sig,\G)$$
le morphisme des complexes de cha\^ines induit par $\{\id_{V_\s}\}_{\s\in\Sig'},$
et $$i^\Sig_{\Sig'}:=H_0(\oplus\id_{V_\s}):H_0(\Sig',\G)\to H_0(\Sig,\G).$$

\begin{lemme}
On a $H_0(\Sig,\underline{V_x})=V_x,$ et pour $y\in\Sig^\circ$ la compos\'ee de morphismes d'homologies
$$
H_0(\{y\},\G)=V_y\To{i^\Sig_y}H_0(\Sig,\G)\To{p^{\Sig}_x}H_0(\Sig,\underline{V_x})=V_x
$$
s'identifie \`a $\e^y_x.$ En particulier, $p^\Sig_x\circ i^\Sig_x=\id_{V_x}.$
\end{lemme}
\begin{preuve}
Comme $\Sig$ est fini convexe donc contractile, le complexe de cha\^ines $\CC_*(\Sig,\underline{V_x})$ est une r\'esolution de $V_x.$ Par d\'efinition de $p^\Sig_x,$ on voit que $p^\Sig_x\circ i^\Sig_y=\e^y_x.$ Or, $\e^x_x$ est l'identit\'e sur $V_x.$ On a $p^\Sig_x\circ i^\Sig_x=\id_{V_x}.$
\end{preuve}

Gr\^ace au lemme pr\'ec\'edent, $V_x$ est un sous $R$-module de $H_0(\Sig,\G).$ Posons alors
$$
e^\Sig_x:H_0(\Sig,\G)\To{p^\Sig_x}V_x\To{i^\Sig_x}H_0(\Sig,\G).
$$
C'est un idempotent de $H_0(\Sig,\G)$ d'image $i^\Sig_x(V_x)$ (voir la proposition qui suit).

\begin{prop}\label{prop1}
Nous avons les propri\'et\'es suivantes:
\begin{description}
\item [(a)] $H_0(\Sig,\G)=\sum_{x\in\Sig^\circ}i^\Sig_x(V_x).$

\item [(b)] $e^\Sig_x\circ e^\Sig_x=e^\Sig_x,$ i.e. $e^\Sig_x$ est un idempotent de $H_0(\Sig,\G).$

\item [(c)] Si $x,y\in \Sig^\circ$ sont adjacents, alors $e^\Sig_x\circ e^\Sig_y=e^\Sig_y\circ e^\Sig_x.$

\item [(d)] Si $x,y,z\in\Sig^\circ,$ $z\in H(x,y)$ et $z,x$ sont adjacents, alors $e^\Sig_x \circ e^\Sig _z \circ e^\Sig_y=e^\Sig_x\circ e^\Sig _y.$
\end{description}
\end{prop}

\begin{preuve}
(a) et (b) sont claires. Pour (c) (resp. (d)) il suffit de montrer que $\forall w\in\Sig^\circ,$ $e^\Sig_x\circ e^\Sig_y=e^\Sig_y\circ e^\Sig_x$ (resp. $e^\Sig_x \circ e^\Sig _z \circ e^\Sig_y=e^\Sig_x\circ e^\Sig _y$) sur $i^\Sig_w(V_w).$ (d) d\'ecoule du corollaire \ref{coro1}. En effet, gr\^ace au lemme pr\'ec\'edent, on a
$$
e^\Sig_x \circ e^\Sig _z \circ e^\Sig_y\circ i^\Sig_w=i^\Sig_x\circ p^\Sig_x\circ i^\Sig_z\circ p^\Sig_z\circ i^\Sig_y\circ p^\Sig_y\circ i^\Sig_w=i^\Sig_x\circ \e^z_x\circ \e^y_z\circ\e^w_y
$$
et
$$
e^\Sig_x\circ e^\Sig_y\circ i^\Sig_w= i^\Sig_x\circ p^\Sig_x\circ i^\Sig_y\circ p^\Sig_y\circ i^\Sig_w=i^\Sig_x\circ \e^y_x\circ \e^w_y.
$$
En vertu du corollaire \ref{coro1}, on obtient l'\'egalit\'e $e^\Sig_x \circ e^\Sig _z \circ e^\Sig_y\circ i^\Sig_w=e^\Sig_x\circ e^\Sig_y\circ i^\Sig_w.$

Pour (c), soit $a$ un \'el\'ement de $V_w,$ $$e^\Sig_x\circ e^\Sig_y\circ i^\Sig_w(a)-e^\Sig_y\circ e^\Sig_x\circ i^\Sig_w(a)=i^\Sig_x\circ \e^y_x\circ \e^w_y(a)-i^\Sig_y\circ\e^x_y\circ\e^w_x(a).$$ Notons $i'^\Sig_x:V_x\to V_x\oplus V_y\oplus\bigoplus_{s\in\Sig^\circ\ba\{x,y\}}V_s$ (resp. $i'^\Sig_y:V_y\to V_x\oplus V_y\oplus\bigoplus_{s\in\Sig^\circ\ba\{x,y\}}V_s$) l'immersion naturelle $u\mapsto (u,0,0\ldots)$ (resp. $v\mapsto (0,v,0\ldots)$). On a donc deux diagrammes commutatifs:
$$
\xymatrix{
V_x\ar[r]^{i'^\Sig_x} \ar[rd]_{i^\Sig_x}             &\bigoplus_{s\in\Sig^\circ}V_s\ar@{->>}[d]  &\text{et}  & V_y\ar[r]^{i'^\Sig_y} \ar[rd]_{i^\Sig_y}             &\bigoplus_{s\in\Sig^\circ}V_s\ar@{->>}[d] \\
& H_0(\Sig,\G) &  &  & H_0(\Sig,\G)
}
$$
Il s'agit alors de montrer que $i'^\Sig_x\circ \e^y_x\circ \e^w_y(a)-i'^\Sig_y\circ\e^x_y\circ\e^w_x(a)\in\partial(\oplus_{\s\in\Sig,\dim\s=1}V_\s).$ Consid\'erons l'\'el\'ement $b:=p^y_{[x,y]}\circ\e^w_y(a)=p^x_{[x,y]}\circ\e^w_x(a)\in V_{[x,y]}$ ({\em cf.} Lemme \ref{n} (b)). Par d\'efinition, on a donc $\varphi^{[x,y]}_x(b)=\varphi^{[x,y]}_x\circ p^y_{[x,y]}\circ\e^w_y(a)=\e^y_x\circ \e^w_y(a)$ et $\varphi^{[x,y]}_y(b)=\varphi^{[x,y]}_y\circ p^x_{[x,y]}\circ\e^w_x(a)=\e^x_y\circ\e^w_x(a).$ Alors, la diff\'erence $i'^\Sig_x\circ \e^y_x\circ \e^w_y(a)-i'^\Sig_y\circ\e^x_y\circ\e^w_x(a)$ est \'egale \`a $\partial(b).$ Ceci d\'emontre l'\'enonc\'e (c) du lemme.
\end{preuve}

\begin{DEf}
Soit $\s$ un simplexe de $\Sig.$ Gr\^ace \`a (c) de la proposition pr\'ec\'edente, on peut d\'efinir un idempotent
$$e^\Sig_\s:=\prod_{x<\s,x\in\Sig^\circ} e^\Sig_x\in\End(H_0(\Sig,\G)).$$
\end{DEf}

\begin{coro}\label{gg}
\begin{description}
\item [(a)] Soient $\tau,\s$ deux simplexes adjacents de $\Sig,$ alors $e^\Sig_\s\circ e^\Sig_\tau=e^\Sig_{[\s,\tau]}.$

\item [(b)] Soient $\s,\tau,\omega$ trois simplexes de $\Sig$ tels que $\omega\in H(\s,\tau).$ Alors $e^\Sig_\s\circ e^\Sig_\omega\circ e^\Sig_\tau=e^\Sig_\s \circ e^\Sig_\tau.$
\end{description}
\end{coro}

\begin{preuve}
En vertu des propri\'et\'es donn\'ees par la proposition \ref{prop1}, on peut adapter la strat\'egie de la d\'emonstration de \cite[Prop. 2.2]{MS} dans notre situation. On laisse les d\'etails au lecteur.
\end{preuve}

\medskip

Comme dans \cite[Thm. 2.12]{MS}, on pose la d\'efinition suivante:
\begin{DEf}
Soit $\Sig'\subset\Sig$ un sous complexe fini convexe de $\Sig,$ on d\'efinit un endomorphisme $u^\Sig_{\Sig'}\in\End(H_0(\Sig,\G))$ par la formule:
$$
u^\Sig_{\Sig'}:=\sum_{\s\in\Sig'}(-1)^{\dim \s}e^\Sig_\s.
$$
\end{DEf}

\begin{prop}\label{i}
L'endomorphisme $u^\Sig_{\Sig'}$ est un idempotent tel que
\begin{align*}
u^\Sig_{\Sig'}(H_0(\Sig,\G))&=\sum_{x\in{\Sig'}^\circ } \im e^\Sig_x\\
\Ker u^\Sig_{\Sig'} &= \bigcap_{x\in{\Sig'}^\circ} \Ker e^\Sig_x\\
H_0(\Sig,\G)&= u^\Sig_{\Sig'}(H_0(\Sig,\G))\oplus \Ker u^\Sig_{\Sig'}
\end{align*}
De plus, si $\Sig'=\Sig,$ $\Ker u^\Sig_\Sig=\bigcap_{x\in\Sig^\circ}\Ker e^\Sig_x=0.$
\end{prop}

\begin{preuve}
Gr\^ace au corollaire \ref{gg}, la strat\'egie de la d\'emonstration de \cite[Thm. 2.12]{MS} s'applique. Si $\Sig'=\Sig,$ $H_0(\Sig,\G)=\sum_{x\in\Sig^\circ}\im e^\Sig_x,$ donc d'apr\`es la proposition \ref{prop1} (a), on a $\Ker u^\Sig_\Sig=0.$
\end{preuve}

\begin{coro}\label{coro2}
Soient $\Sig_+,\Sig_-$ deux sous complexes finis convexes de $\Sig$ tels que $\Sig=\Sig_+\cup \Sig_-.$ Notons $\Sig_0:=\Sig_+\cap \Sig_-$ leur intersection qui est encore convexe. Alors $u^\Sig_\Sig=u^\Sig_{\Sig_+}+u^\Sig_{\Sig_-}-u^\Sig_{\Sig_0},$ $u^\Sig_{\Sig_+}u^\Sig_{\Sig_-}=u^\Sig_{\Sig_-}u^\Sig_{\Sig_+}=u^\Sig_{\Sig_0}$ et
\begin{align*}
u^\Sig_{\Sig_+}(H_0(\Sig,\G))\cap u^\Sig_{\Sig_-}(H_0(\Sig,\G))&=u^\Sig_{\Sig_0}(H_0(\Sig,\G))\\
\Ker u^\Sig_{\Sig_+}+\Ker u^\Sig_{\Sig_-} &=\Ker u^\Sig_{\Sig_0}.
\end{align*}
\end{coro}

\begin{preuve}
On peut utiliser la strat\'egie de la preuve de \cite[Corollary 2.13]{MS}.
\end{preuve}

\subsection{D\'emonstration du th\'eor\`eme \ref{thm2}}
D\'emontrons tout d'abord la premi\`ere partie du th\'eor\`eme \ref{thm2}.

\begin{prop}
Soient $\Sig$ un sous complexe fini convexe de $\BC\TC$ fix\'e et $\Sig'\subset\Sig$ un sous complexe convexe, alors
\begin{description}
\item [(a)] $H_n(\Sig',\G)=0$ pour tout $n>0.$

\item [(b)] L'application $i^{\Sig}_{\Sig'}:H_0(\Sig',\G)\To{} H_0(\Sig,\G)$ (voir \ref{j}) est injective d'image $u^\Sig_{\Sig'}(H_0(\Sig,\G)).$
\end{description}
\end{prop}

\begin{preuve}
Tout d'abord on d\'emontre (b). Notons que l'on a un diagramme commutatif:
$$
\xymatrix{
H_0(\Sig',\G)\ar[r]^{i^{\Sig}_{\Sig'}}\ar[d]^{p_x^{\Sig'}}  & H_0(\Sig,\G)\ar[d]^{p_x^\Sig} \\
V_x\ar[r]^{\id}  & V_x.
}
$$
Soit $a\in H_0(\Sig',\G)$ tel que $i^{\Sig}_{\Sig'}(a)=0.$ Donc $p_x^{\Sig'}(a)= p_x^{\Sig}\circ i^{\Sig}_{\Sig'}(a)=0,~\forall x\in\Sig'^{\circ}.$ On en d\'eduit que $e^{\Sig'}_x(a)=i^{\Sig'}_{x}\circ p_x^{\Sig'}(a)=0,~\forall x\in\Sig'^{\circ}.$ Donc $a\in \bigcap_{x\in\Sig'^\circ}\Ker e^{\Sig'}_x=\{0\}$ ({\em cf.} \ref{i}).

D'apr\`es la proposition \ref{i}, on a
\begin{align*}
u^{\Sig}_{\Sig'}(H_0(\Sig,\G))&=\sum_{x\in\Sig'^\circ}\im e^\Sig_x=\sum_{x\in\Sig'^\circ}i^\Sig_x(V_x)\\
&=\sum_{x\in\Sig'^\circ}i^\Sig_{\Sig'}\circ i^{\Sig'}_x(V_x)=i^\Sig_{\Sig'}(\sum_{x\in\Sig'^\circ}i^{\Sig'}_x(V_x))\\
&=i^\Sig_{\Sig'}(u^{\Sig'}_{\Sig'}H_0(\Sig',\G))=i^\Sig_{\Sig'} H_0(\Sig',\G).
\end{align*}
\medskip

On d\'emontre (a) par r\'ecurrence.
\begin{lem}
Supposons que $\Sig$ ne soit pas un simplexe, et $\forall \Sig'\subsetneqq \Sig$ fini convexe la propri\'et\'e $(a)$ de la proposition soit v\'erifi\'ee pour $\Sig'.$ Alors elle est v\'erifi\'ee pour $\Sig.$
\end{lem}

\begin{preuve}
Lorsque $\Sig$ n'est pas un simplexe, gr\^ace \`a \cite[Section 2.5]{MS}, on peut d\'ecomposer $\Sig=\Sig_+\cup\Sig_-$ de sorte que $\Sig_+,\Sig_-$ et leur intersection $\Sig_0$ soient sous complexes finis convexes propres de $\Sig.$ En particulier, (a) est v\'erifi\'ee pour $\Sig_+,\Sig_-$ et $\Sig_0.$ Les complexes de cha\^ines associ\'es forment une suite exacte de complexes
$$
\CC_*(\Sig_0,\G)\rightarrowtail \CC_*(\Sig_+,\G)\oplus \CC_*(\Sig_-,\G)\onto \CC_*(\Sig,\G),
$$
dont la suite exacte longue d'homologie du type Mayer-Vietoris associ\'ee implique que $H_n(\Sig,\G)=0$ pour $n\geq 2.$ De plus, l'injectivit\'e de l'application $H_0(\Sig_0,\G)\To{i^{\Sig_+}_{\Sig_0}}H_0(\Sig_+,\G)$ implique que $H_1(\Sig,\G)=0.$ Donc (a) est v\'erifi\'ee pour $\Sig.$
\end{preuve}

Il nous reste \`a traiter le cas o\`u $\Sig$ est un simplexe. Soit $\Sig$ un simplexe, rappelons que pour $x\in\Sig^\circ$ et $\s\subset\Sig$ un sous simplexe, on a d\'efini un idempotent $\wt{\e^\s_x}\in\End(V_\s)$ satisfaisant $\wt{\e^\s_x}\circ\wt{\e^\s_y}=\wt{\e^\s_y}\circ\wt{\e^\s_x}$ (voir le lemme \ref{lem1}). Posons comme \cite[(2.5)]{MS} des idempotents
$$
\e^{\s,0}_I:=\prod_{x\in I}\wt{\e^\s_x}\prod_{x\not\in I}(1-\wt{\e^\s_x})\in\End(V_\s)
$$
pour un sous-ensemble $I$ de $\Sig^\circ.$ On a comme dans {\em loc. cit.}
\begin{itemize}
\item $\e^{\s,0}_I(V_\s)=\left\{
                           \begin{array}{ll}
                             0, & \hbox{si $\s^\circ\nsubseteq I$;} \\
                             \e^{\s,0}_I(V_x) \text{ pour un sommet } x\in I, & \hbox{si $\s^\circ\subset I$.}
                           \end{array}
                         \right.
$

\item $\id_{V_\s}=\sum_{I\subset \Sig^\circ}\e^{\s,0}_I\in\End(V_\s)$ et $\e^{\s,0}_I\e^{\s,0}_J=0$ si $I\neq J.$
\end{itemize}

Notons $\e^0_I:=\oplus_\s\e^{\s,0}_I\in\End(\CC_*(\Sig,\G))$ commutant avec la diff\'erentielle. On en d\'eduit une d\'ecomposition
$$
\CC_*(\Sig,\G)=\bigoplus_{I\subset\Sig^\circ}\e^0_I(\CC_*(\Sig,\G)).
$$
Si $I=\emptyset,$ $\e^0_I(\CC_*(\Sig,\G))=0.$ Si $I\neq \emptyset,$ le complexe de cha\^ines $\e^0_I(\CC_*(\Sig,\G))$ calcule l'homologie du sous complexe $\Sig_I$ de $\Sig$ engendr\'e par $I$ \`a valeurs constants $\e^0_I(V_x),$  pour un sommet quelconque $x$ contenu dans $I.$ Comme $\Sig_I$ est contractile, $\e^0_I(\CC_*(\Sig,\G))$ est une r\'esolution de $\e^0_I(V_x).$ Ceci entra\^ine que $H_n(\Sig,\G)=0$ pour $n\geq 1.$
\end{preuve}
\medskip

Maintenant, on peut montrer que $\CC_*(\BC\TC,\G)$ est une r\'esolution de $H_0(\BC\TC,\G)$ ({\em cf.} \cite{MS}). Posons $(\Sig_n)_{n\in\NM}$ une suite croissante de sous complexes finis convexes telle que $\BC\TC=\bigcup_n \Sig_n$ et $\CC_*(\BC\TC,\G)=\limi{}\CC_*(\Sig_n,\G).$ D'apr\`es ce qui pr\'ec\`ede, on sait que pour chaque $n,$ le complexe
$$
\CC_*(\Sig_n,\G)\To{} H_0(\Sig_n,\G)\To{}0
$$
est exact. Or, la limite inductive pr\'eserve l'exactitude dans la cat\'egorie de $R$-modules. On a alors que $\CC_*(\BC\TC,\G)$ est une r\'esolution de $\limi{} H_0(\Sig_n,\G).$ C'est-\`a-dire $H_n(\BC\TC,\G)=0$ pour $n\geq 1$ et $H_0(\BC\TC,\G)=\limi{}H_0(\Sig_n,\G).$ Ceci entra\^ine la premi\`ere partie du th\'eor\`eme.

\ali On d\'emontre la seconde partie du th\'eor\`eme \ref{thm2}. Soit $x$ un sommet quelconque. Notons que pour $\Sig'\subset\Sig$ deux complexes finis convexes contenant $x,$ on a un diagramme commutatif:
$$
\xymatrix{
V_x\ar[r]^{i^{\Sig'}_x}\ar[rd]_{i^{\Sig}_x}  & H_0(\Sig',\G)\ar[d]^{i^{\Sig}_{\Sig'}}\\
                         & H_0(\Sig,\G)
}
$$
Ceci induit une injection $i_x:V_x\to H_0(\BC\TC,\G).$

\begin{lem}
Soient $x,y$ deux sommets adjacents, alors on a un diagramme commutatif:
$$
\xymatrix{
V_x \ar[r]^{i_x}\ar[d]_{\e^x_y}       &  H_0(\BC\TC,\G)\ar[d]^{e_y}\\
V_y \ar[r]^{i_y}  & H_0(\BC\TC,\G)
}
$$
\end{lem}
\begin{preuve}
La d\'emonstration suit la ligne de la d\'emonstration de \ref{prop1} (c). Notons $i'_x:V_x\to V_x\oplus V_y\oplus\bigoplus_{s\in\BC\TC^\circ\ba\{x,y\}}V_s$ l'injection naturelle $u\mapsto (u,0,0\ldots),$ et $i'_y:V_y\to V_x\oplus V_y\oplus\bigoplus_{s\in\BC\TC^\circ\ba\{x,y\}}V_s$ l'injection naturelle $v\mapsto (0,v,0\ldots)$ telles que $i_x(u)=f(i'_x(u))$ et $i_y(v)=f(i'_y(v))$ o\`u $f$ est l'application canonique $f:\bigoplus_{s\in\BC\TC^\circ}V_s\onto H_0(\BC\TC,\G).$

Lorsque $x=y,$ on a $e_x(i'_x(a))=i'_x(a),~\forall a\in V_x.$ C'est-\`a-dire, $i_x\circ \e^x_x=e_x\circ i_x.$

Lorsque $x\neq y,$ soit $a$ un \'el\'ement de $V_x,$ il s'agit de montrer que $$e_y(i'_x(a))-i'_y(\e^x_y(a))\in \partial(\bigoplus_{\s\in\BC\TC_1}V_\s).$$ Par d\'efinition, on a
$$
e_y(i'_x(a))-i'_y(\e^x_y(a))=e_ye_x(i'_x(a))-i'_y(\e^x_y(a))=(e_{[x,y]}(a),-\e^x_y(a),0\ldots).
$$
Posons $b:=p^x_{[x,y]}(a)\in V_{[x,y]},$ alors $e_{[x,y]}(a)=\varphi^{[x,y]}_x(b)$ et $\e^x_y(a)=\varphi^{[x,y]}_y(b).$ On en d\'eduit que $(e_{[x,y]}(a),-\e^x_y(a),0\ldots)$ appartient \`a $\partial(\bigoplus_{\s\in\BC\TC_1}V_\s).$
\end{preuve}
D'apr\`es ce qui pr\'ec\`ede, le morphisme $i_x$ se factorise par $e_x(H_0(\BC\TC,\G)).$ Soit $\s$ un simplexe contenant $x,$ $i_\s:=i_x|_{e_\s(V_x)}$ envoie $\varphi^\s_x(V_\s)=e_\s(V_x)$ dans $e_\s(H_0(\BC\TC,\G)),$ et il ne d\'epend pas du choix de sommet contenu dans $\s.$ Donc les morphismes $\{i_\s\}_{\s\in\BC\TC}$ induisent un morphisme de syst\`emes de coefficients
$$\G\To{} (\s\mapsto e_\s(H_0(\BC\TC,\G)))_{\s\in\BC\TC}.$$
En tenant compte du fait que $H_0(\BC\TC,\G)=\sum_{y\in\BC\TC^\circ} i_y(V_y),$ il suffit de montrer que $\forall y\in\BC\TC^\circ,$ $e_\s(i_y(V_y))\subset i_x(V_\s).$ Prenons un chemin tendu $z_0=y,\ldots,z_m=x$ de $y$ vers $x,$ on a alors
$$
e_\s(i_y(V_y))=e_\s(e_y i_y(V_y))=e_\s e_{z_m}e_y(i_y(V_y))=\cdots=e_\s e_{z_m}\cdots e_{z_2}e_{z_1}e_y(i_y(V_y)).
$$
D'apr\`es le lemme pr\'ec\'edent, on a
\begin{align*}
e_\s e_{z_m}\cdots e_{z_2}e_{z_1}e_y(i_y(V_y))&=e_\s e_{z_m}\cdots e_{z_2}i_{z_1}(\e^y_{z_1}(V_y))\subset e_\s e_{z_m}\cdots e_{z_2}i_{z_1}(V_{z_1})\\
&=e_\s e_{z_m}\cdots i_{z_2}(\e^{z_1}_{z_2}(V_{z_1}))\subset e_\s e_{z_m}\cdots i_{z_2}(V_{z_2})\\
&\subset\cdots\subset e_\s i_x(V_x)=i_x(V_\s),
\end{align*}
car $e_\s \circ i_x=i_\s\circ \e^x_\s.$ Ceci termine la d\'emonstration du th\'eor\`eme.

\section{La cohomologie du rev\^etement mod\'er\'e de l'espace de Drinfeld}

Dans cette section, on s'int\'eresse \`a la partie {\em non-supercuspidale} de la cohomologie du rev\^etement mod\'er\'e de l'espace de Drinfeld (la partie {\em supercuspidale} a \'et\'e trait\'ee dans \cite{Wang-Sigma1}). Rappelons tout d'abord les r\'esultats g\'eom\'etriques obtenus dans \cite{Wang-Sigma1}.

\subsection{Rappels sur le rev\^etement mod\'er\'e de l'espace de Drinfeld}

\ali\label{3.1::S1} Soient $K$ une extension finie de
$\QM_p$ d'anneau des entiers $\OC$ et $\varpi$ une uniformisante de
$\OC$ tels que $\OC/\varpi\simeq\FM_q.$ On
fixe $K^{ca}$ une cl\^oture alg\'ebrique de $K$ et $\widehat{K^{ca}}$
son compl\'et\'e $\varpi$-adique. Soient $D$ l'alg\`ebre \`a division centrale sur
$K$ d'invariant $1/d,$ $\OC_D$ l'anneau des entiers de $D$ et $\Pi_D$ une uniformisante. Soient
$K_d$ une extension non-ramifi\'ee de degr\'e $d$ de $K$ contenue
dans $D$ d'anneau des entiers de $\OC_{K_d}$. On note $\breve{K}$ le compl\'et\'e de
l'extension non-ramifi\'ee maximale $K^{nr}\subset K^{ca}$ de $K,$ $\breve{\OC}$ son anneau des entiers. On d\'esigne $W_K$ le groupe de Weil associ\'e \`a $K,$ et $I_K$ le groupe d'inertie. Fixons un rel\`evement $\varphi\in\Gal(K^{ca}/K)$ de $(\Frob_q:x\mapsto x^{-q})\in\Gal(\oFq/\FM_q).$

Soit $d\geq 2$ un entier, notons $\Omega_K^{d-1}$ l'espace sym\'etrique de Drinfeld \cite{drinfeld-ell} de dimension $d-1,$ d\'efini comme le compl\'ementaire de l'ensemble des hyperplans $K$-rationnels dans l'espace projectif $\PM^{d-1}_K.$ Il admet une structure d'espace rigide-analytique au sens de Raynaud-Berkovich. Notons $|\BC\TC|$ la r\'ealisation g\'eom\'etrique de l'immeuble de Bruhat-Tits semi-simple $\BC\TC$ associ\'e \`a $G:=\GL_d(K),$ et rappelons que nous avons une {\em application de r\'eduction} $\tau:\Omega_K^{d-1}\to|\BC\TC|.$ Pour $\s=\{s_0,\ldots,s_k\}$ avec $s_i$ des sommets un $k$-simplexe de $\BC\TC,$ notons $|\s|\subset|\BC\TC|$ sa facette associ\'ee et $|\s|^*:=\bigcup_{\s\subset\s'}|\s'|=\bigcap_i|s_i|^*.$ Donc, les donn\'ees $\{\tau^{-1}(|\s|^*)\}_{\s\in\BC\TC}$ fournissent un recouvrement admissible de $\Omega_K^{d-1}.$

Deligne a construit un mod\`ele semi-stable $\wh{\Omega}^{d-1}_\OC$ de $\Omega_K^{d-1}$ sur $\Spf \OC$ en recollant les mod\`eles locaux $\{\wh{\O}^{d-1}_{\OC,\s}\}_{\s\in\BC\TC},$ tel que sa fibre sp\'eciale g\'eom\'etrique $\o{\Omega}:=\wh{\Omega}^{d-1}_{\OC}\otimes_\OC\o{\FM}_q$ est une r\'eunion des composantes irr\'eductibles param\'etr\'ees par les sommets de $\BC\TC$ ({\em cf.} \cite[2.1.5]{Wang-Sigma1}) i.e. $\o{\O}=\bigcup_{s\in\BC\TC^\circ}\o{\Omega}_s.$ Soit $\s=\{s_0,\ldots,s_k\},$ notons $\o\O_\s:=\bigcap_i\o\O_{s_i},$ $\o{\Omega}_\s^0:=\o{\Omega}_\s\backslash\bigcup_{s'\not\in \s}\o{\Omega}_{s'},$ et $j_\s:\o{\Omega}_\s^0\into\o{\Omega}_\s$ l'inclusion naturelle. Rappelons que Berkovich a d\'efini dans ce cas un morphisme de sp\'ecialisation $$\spe:\Omega^{d-1,ca}_K:=\Omega_K^{d-1}\widehat{\otimes}_K\widehat{K^{ca}}\To{}\o{\Omega}$$ (qui est appel\'e le {\em morphisme de r\'eduction} dans \cite[\S 1]{berk-vanishing}).

\ali Dans \cite{drinfeld-covering}, Drinfeld a d\'emontr\'e que le sch\'ema formel $\wh{\O}^{d-1}_{\OC}\wh{\otimes}_{\OC}\breve{\OC}$ (pro)-repr\'esente l'espace de modules des $\OC_D$-modules formels sp\'eciaux de dimension $d$ et hauteur $d^2$ sur la cat\'egorie de $\breve{\OC}$-alg\`ebres dans lesquelles $\varpi$ est nilpotent. Notons $\XG$ le $\OC_D$-module formel sp\'ecial universel, et posons $\Sigma_n:=\underline{\Isom}_{\OC_D}(\Pi_D^{-n}\OC_D/\OC_D,\XG[\Pi^n_D]^{an}),$ o\`u $\XG[\Pi_D^n]$ d\'esigne les $\Pi_D^n$-torsions. Par construction, $\Sig_n$ est un $(\OC_D/\Pi_D^n\OC_D)^\times$-torseur sur $\Omega^{d-1}_K\wh{\otimes}_K\breve{K}$ muni d'une action de $G^\circ\times \OC_D^\times\times I_K,$ o\`u $$G^\circ:=\Ker(\val_K\circ\det:\GL_d(K)\to K^\times).$$ Les morphismes de transitions $\Sig_n\To{\times\Pi_D} \Sig_{n-1}$ sont $G^\circ\times \OC_D^\times\times I_K$-\'equivariants. Lorsque $n=1,$ l'espace rigide-analytique $\Sig^{ca}_1:=\Sigma_1\widehat{\otimes}_K\widehat{K^{ca}}$ est un $(\OC_D/\Pi_D\OC_D)^\times\simeq\FM_{q^d}^\times$-torseur sur $\Omega_K^{d-1,ca}:=\Omega^{d-1}_K\wh{\otimes}_K\wh{K^{ca}},$ que l'on appelle {\em le rev\^etement mod\'er\'e} de l'espace de Drinfeld de dimension $d-1.$

Notons $p:\Sigma^{ca}_1\to\Omega_K^{d-1,ca}$ la mutiplication par $\Pi_D,$ et $\nu:=\tau\circ p.$ On a un diagramme commutatif dont les fl\`eches sont $G^\circ$-\'equivariantes
$$\xymatrix{
\Sig^{ca}_1 \ar[d]^p \ar[rd]^\nu &\\
\O^{d-1,ca}_K\ar[r]^{\tau}& |\BC\TC|
}
$$
et tel que $\Sigma^{ca}_1$ admet un recouvrement des ouverts admissibles $\{\nu^{-1}(|\s|^*)\}_{\s\in\BC\TC}.$

Fixons $\ell\neq p$ un nombre premier et une cl\^oture alg\'ebrique $\oQl$ de $\QM_\ell.$ Soient $\s'\subset \s$ deux simplexes de $\BC\TC$, l'immersion ouverte $\nu^{-1}(|\s|^*)\into\nu^{-1}(|\s'|^*)$ induit un morphisme canonique:
\[
H^q_c(\nu^{-1}(|\s|^*),\L)\longrightarrow H^q_c(\nu^{-1}(|\s'|^*),\L),~\forall q\geq0
\]
o\`u $\L=\oZl$ ou $\oQl.$ Donc les donn\'ees
$$
\begin{cases}
(\s\in\BC\TC)\mapsto H^q_c(\nu^{-1}(|\s|^*),\L)\\
(\s'\subset \s)\mapsto \big(H^q_c(\nu^{-1}(|\s|^*),\L)\to H^q_c(\nu^{-1}(|\s'|^*),\L)\big)
\end{cases}
$$
d\'efinissent un syst\`eme de coefficients $G^\circ$-\'equivariant (voir \ref{c}) \`a valeurs dans les $\L$-modules que nous noterons simplement $\s\mapsto H^q_c(\nu^{-1}(|\s|^*),\L).$ Ce syst\`eme de coefficients calcule la cohomologie de $\Sig^{ca}_1.$

\begin{fact}(\cite[Prop. 3.2.4]{dat-drinfeld})\label{suite spec}
Il existe une suite spectrale $G^\circ$-\'equivariante
\ini\begin{equation}\label{suite spec1}
E^{pq}_1=\CC^{or}_c(\BC\TC_{(-p)},\s\mapsto H^q_c(\nu^{-1}(|\s|^*),\Lambda))\Longrightarrow H^{p+q}_c(\Sigma^{ca}_1,\Lambda)
\end{equation}
dont la diff\'erentielle $d^{pq}_1$ est celle du complexe de cha\^ines du syst\`eme de coefficients $\s\mapsto H^q_c(\nu^{-1}(|\s|^*),\Lambda).$
\end{fact}

On regroupe ci-dessous les r\'esultats dans \cite[\S 2]{Wang-Sigma1} (o\`u $\Sig_1$ est not\'e $\Sig$).

\begin{theo}\label{3.1::T1} Soient $s$ un sommet de $\BC\TC$ et $\s$ un simplexe contenant $s.$ Posons $\breve{K}^t=\breve{K}[\varpi_t]/(\varpi_t^{q^d-1}-\varpi)$ une
extension mod\'er\'ement ramifi\'ee de degr\'e $q^d-1$ de
$\breve{K}$ d'anneau des entiers $\breve{\OC}^t,$ notons $\Sigma_{1,s}$ l'espace rigide $\Sigma_1\times_{\Omega_{K}^{d-1}}\tau^{-1}(|s|),$ et consid\'erons la normalisation de $\wh\O^{d-1}_{\OC,s}\wh\otimes_{\OC}\breve{\OC}^t$ dans $\Sigma_{1,s}\otimes_{\breve{K}}\breve{K}^t$ que l'on notera $\wh{\Sigma}^0_{1,s}.$ Alors,
\begin{description}
\item[(a)] $\wh{\Sigma}^0_{1,s}$ est un $\FM_{q^d}^\times$-torseur $G_s$-invariant au-dessus de $\wh\O^{d-1}_{\OC,s}\wh\otimes_{\OC}\breve{\OC}^t;$ et si l'on note $\o{\Sig}_{1,s}^0$ sa fibre sp\'eciale, on en d\'eduit un isomorphisme $I_K/I_{K(\varpi_t)}\cong\FM_{q^d}^\times$-\'equivariant
    $$
H^q(\nu^{-1}(|s|),\L)\cong H^q(\o{\Sigma}_{1,s}^0,\L),~\forall q\geq0.
$$
De plus, il existe une donn\'ee de descente \`a la Weil sur $\Sig^{ca}_1$ (qui est induite par $\Pi_D^{-1}\circ\varphi$ sur $\MC_{Dr,1}^{ca}$ {\em cf.} \ref{3.2::S2}) telle qu'elle induit une donn\'ee de descente $\Fr$ sur $\o{\Sigma}_{1,s}^0.$

\item[(b)] Quitte \`a choisir une base d'un r\'eseau qui repr\'esente $s,$ on a un isomorphisme $G_s/G_s^+\cong \GL_d(\FM_q)$-\'equivariant $\o{\Sig}^0_{1,s}\cong \DL^{d-1}_{\oFq}:=\DL^{d-1}\otimes_{\FM_q}\oFq$ o\`u $\DL^{d-1}$ d\'esigne la sous-vari\'et\'e ferm\'ee de l'espace affine $\AM^d_{\FM_q}=\Spec {\FM_q[X_0,\ldots,X_{d-1}]}$ d\'efinie par l'\'equation
\ini\begin{equation}\label{Eq:2}
\det((X_i^{q^j})_{0\leq i,j\leq d-1})^{q-1}=(-1)^{d-1}.
\end{equation}
De plus, cet isomorphisme est compatible aux actions de $\FM_{q^d}^\times,$ o\`u $\FM_{q^d}^\times$ agit sur $\DL^{d-1}_{\oFq}$ par $X_i\mapsto \zeta X_i,~\forall \z\in\FM_{q^d}^\times;$ la structure $\FM_q$-rationnelle induite par $\Fr$ sur $\o{\Sigma}_{1,s}^0$ co\"incide via cet isomorphisme avec celle induite par $\Frob:X_i\mapsto X_i^q$ sur $\DL^{d-1}_{\oFq}.$

\item[(c)] L'immersion ouverte $\nu^{-1}(|s|)\into\nu^{-1}(|s|^*)$ induit un isomorphisme:
\[
H^q(\nu^{-1}(|s|^*),\Lambda)\simto H^q(\nu^{-1}(|s|),\Lambda),~\forall q\geq0.
\]

\item[(d)] Soit $\s$ un simplexe contenant $s.$ Le morphisme canonique $H^q_c(\nu^{-1}(|\s|^*),\Lambda)\to H^q_c(\nu^{-1}(|s|^*),\Lambda)$ provenant de l'immersion ouverte $\nu^{-1}(|\s|^*)\into\nu^{-1}(|s|^*)$ induit un isomorphisme
\begin{equation*}
H^q_c(\nu^{-1}(|\s|^*),\Lambda)\simto H^q_c(\nu^{-1}(|s|^*),\Lambda)^{G_\s^+}.
\end{equation*}
\end{description}
\end{theo}
\begin{preuve}
(a) d\'ecoule de \cite[2.3.8, 2.3.9, 3.1.8]{Wang-Sigma1}. (b) d\'ecoule de \cite[2.4.6, 2.5.3, 3.1.10]{Wang-Sigma1}. L'assertion (c) est \cite[Th\'eor\`eme 2.2.3]{Wang-Sigma1} qui est une cons\'equence de \cite[Lemme 5.6]{zheng}. La d\'emonstration de (d) est donn\'e dans \cite[Th\'eor\`eme 2.5.9]{Wang-Sigma1} qui est une cons\'equence du th\'eor\`eme principal de \cite{wang_DL}.
\end{preuve}

\begin{coro}\label{3.1::C1}
Notons $\G_q(\L)$ le syst\`eme de coefficients $\s\mapsto H^q_c(\nu^{-1}(|\s|^*),\L).$ Alors, la suite spectrale \ref{suite spec1} d\'eg\'en\`ere en $E_1,$ et elle induit un isomorphisme
$$H^q_c(\Sig^{ca}_1,\L)\simto H_0(\BC\TC,\G_q(\L)),~\forall q\geq0.$$
\end{coro}
\begin{preuve}
D'apr\`es (d) du th\'eor\`eme pr\'ec\'edent, on sait que $\G_q(\L)$ est un syst\`eme de coefficients $G^\circ$-\'equivariant de niveau 0 ({\em cf.} \ref{h}). Donc, l'\'enonc\'e du corollaire d\'ecoule du th\'eor\`eme \ref{thm2}.
\end{preuve}

\subsection{La cohomologie \`a coefficients entiers}\label{3.2::S1}

Dans cette partie, on consid\`ere la cohomologie \`a coefficients entiers, i.e. $\L=\oZl.$ On introduit tout d'abord quelques notations. Posons $GDW:=G\times D^\times \times W_K,$ et $v:GDW\to \ZM$ l'homorphisme qui envoie un \'el\'ement $(g,\d,w)$ vers l'entier $\val_K(\det(g^{-1})\Nr(\d)\Art^{-1}(w))\in\ZM,$ o\`u $\Nr:D^\times\to K^\times$ d\'esigne la norme r\'eduite et $\Art^{-1}:W_K\onto W_K^{ab}\to K^\times$ d\'esigne la compos\'ee de l'inverse du morphisme d'Artin qui envoie l'uniformisante $\varpi$ vers le Frobenius g\'eom\'etrique $\varphi$ ({\em cf.} \ref{3.1::S1}). On d\'esigne $(GDW)^0:=v^{-1}(0)$ le noyau de $v.$ Pour $f|d,$ on notera $[GDW]_{f}$ le sous-groupe distingu\'e form\'e des \'el\'ements $(g,\d,w)$ tels que $v(g,\d,w)\in f\ZM.$ On consid\`ere $G$ (resp. $D^\times,$ $W_K$) comme un sous-groupe de $GDW$ via l'inclusion naturelle, et on notera $[G]_{f}$ (resp. $[D]_{f},$ $[W_K]_{f}$) son intersection avec $[GDW]_{f}.$ D\`es que l'on identifie $K^\times$ au centre de $G,$ on a $[GDW]_d=(GDW)^0\varpi^\ZM.$

\ali\label{3.2::S2} On utilise le langage de Rapoport-Zink \cite{RZ} afin de d\'efinir une action du groupe $GDW.$ Rappelons bri\`evement (voir \cite[3.1.4]{Wang-Sigma1} \cite[3.1]{Dat-elliptic} pour plus des d\'etails) que nous avons alors un sch\'ema formel $\wh{\MC}_{Dr,0}$ isomorphe (non-canoniquement) \`a $\wh\O^{d-1}_{\OC}\times\ZM$ de fibre g\'en\'erique g\'eom\'etrique $\MC^{ca}_{Dr,0}\cong \O^{d-1,ca}_{K}\times\ZM.$ $\MC_{Dr,0}$ poss\`ede un rev\^etement \'etale $\MC^{ca}_{Dr,1}\cong\Sig_1^{ca}\times\ZM$ de groupe de Galois $\FM_{q^d}^\times.$ Il existe une action de $GDW$ sur $\MC^{ca}_{Dr,1}$ telle que la composante $\Sig^{ca}_1\times\{0\}$ soit stable sous $(GDW)^0,$ et que $\varpi\in K^\times\subset G$ permute les composantes par $n\mapsto n+d$ sur $\ZM.$

Posons
$$
\Hb^{q}_{c,\oZl}:=H^q_c(\MC_{Dr,1}^{ca}/\varpi^\ZM,\oZl), ~\forall q\geq0.
$$
D'apr\`es la description pr\'ec\'edente,
$$
\Hb^q_{c,\oZl}=\Ind^{GDW}_{[GDW]_d}H^{q}_c(\Sig_1^{ca},\oZl),
$$
en d\'eclarant que $\varpi$ agit trivialement sur $H^{q}_c(\Sig_1^{ca},\oZl).$

En combinant avec les propri\'et\'es connues des vari\'et\'es de Deligne-Lusztig, on d\'eduit les corollaires suivants.

\begin{coro}\label{Cor::2}
Pour tout entier $q\geq 0,$ le $\oZl G$-module $\Hb^q_{c,\oZl}$ est admissible.
\end{coro}

\begin{preuve}
 Il suffit de montrer que $H^q_c(\Sig^{ca}_{1},\oZl)$ est un $\oZl G^\circ$-module admissible. D'apr\`es \ref{3.1::C1}, \ref{thm2} (b) et \ref{3.1::T1} (a) (b), $$H^q_c(\Sig^{ca}_1,\oZl)^{G_x^+}\simto H_0(\BC\TC,\G_q(\oZl))^{G^+_x}\simto H^q_c(\DL^{d-1}_{\oFq},\oZl),$$ pour tout sommet $x$ de $\BC\TC.$ Notons que $\DL^{d-1}_{\oFq}$ est une vari\'et\'e affine sur $\oFq.$ D'o\`u l'\'enonc\'e du corollaire.
\end{preuve}

\begin{coro}\label{Cor::3}
Pour tout entier $q\geq 0,$ le $\oZl$-module $\Hb^q_{c,\oZl}$ est sans torsion et non-divisible.
\end{coro}
\begin{preuve}
Il suffit de montrer les m\^eme \'enonc\'es pour $H^q_c(\Sig^{ca}_{1},\oZl).$ On d\'emontre tout d'abord que le $\oZl$-module $$C:=H^q_c(\DL^{d-1}_{\oFq},\oZl)$$ est sans torsion.
Ceci d\'ecoule de \cite[Lemma 3.9, Corollary 4.2]{bonnafe-rouquier-coxeter}. En effet, notons $\o{G}:=\GL_d(\FM_q),$ et notons $C_{\tors}\subset C$ la partie de torsions. Bonnaf\'e et Rouquier d\'efinissent dans \cite[3.3.4]{bonnafe-rouquier-coxeter} un \'el\'ement $e(\psi)\in \oZl \o G$ tel que $\oZl \o{G} e(\psi)$ soit un prog\'en\'erateur de $\oZl\o G,$ {\em cf. loc. cit.} Corollary 4.2 (voir \cite{Bonnafe} pour une approache plus \'el\'ementaire).  Le $\oZl$-module $\oZl\o{G} e(\psi)C_{\tors}$ est un sous-module de $\oZl \o{G} e(\psi) C$ qui est un $\oZl$-module libre de type fini d'apr\`es \cite[Lemma 3.9]{bonnafe-rouquier-coxeter}. Donc $\oZl\o{G} e(\psi)C_{\tors}=0.$ Ceci entra\^ine que $C_{\tors}=0$ d'apr\`es \cite[Corollary 4.2]{bonnafe-rouquier-coxeter}.

Notons que le $\oZl G^\circ$-module $H^q_c(\Sig_1^{ca},\oZl)$ est engendr\'e par $$H^q_c(\Sig_1^{ca},\oZl)^{G^+_x}\cong H^q_c(\DL^{d-1}_{\oFq},\oZl)$$ pour tout sommet $x\in\BC\TC,$ i.e.
$$
 H^q_c(\Sig_1^{ca},\oZl)=\sum_{x\in\BC\TC^\circ}H^q_c(\Sig^{ca}_1,\oZl)^{G^+_x}.
$$
On en d\'eduit que $H^q_c(\Sig^{ca}_1,\oZl)$ est sans torsion. En vertu de l'admissibilit\'e dans le corollaire pr\'ec\'edent, $H^q_c(\Sig^{ca}_1,\oZl)$ est non-divisible.
\end{preuve}

\subsection{Rappels sur les repr\'esentations elliptiques}
Avant de donner des cons\'equences sur la partie non-supercuspidale de la cohomologie, nous rappelons les propri\'et\'es des repr\'esentions elliptiques de $\GL_d(\FM_q)$ et de $\GL_d(K).$

\ali Soit $k$ un corps fini. Si $\o\pi_i$ est une repr\'esentation de $\GL_{d_i}(k),$ $i=1,\ldots,t,$ on note $\times_i\o\pi_i:=\o\pi_1\times\cdots\times\o\pi_t$ l'induite de la repr\'esentation $\o\pi_1\otimes\cdots\otimes\o\pi_t$ du sous-groupe de Levi diagonal par blocs $\GL_{d_1}(k)\times\cdots\times \GL_{d_t}(k)$ de $\GL_{d_1+\cdots+d_t}(k),$ le long du parabolique triangulaire par blocs sup\'erieur correspondant $\o P_{(d_1,\ldots,d_t)}.$ Notons $R(\GL_d(k))$ le groupe de Grothendieck des $\oQl$-repr\'esentations de longueur finie de $\GL_d(k).$

\begin{Def}
Soient $k=\FM_q$ et $f|d.$ On dit qu'un caract\`ere $\th':\FM_{q^f}^\times\to\oQl^\times$ est {\em $f$-primitif}, si $\th'$ ne se factorise pas par la norme $N_{\FM_{q^f}/\FM_{q^{f'}}}:\FM_{q^{f}}^\times\onto\FM_{q^{f'}}^\times,~\forall f'|f$ et $f'\neq f.$
\end{Def}

Un caract\`ere $f$-primitif $\th'$ correspond \`a une repr\'esentation irr\'eductible cuspidale $\o\pi_f(\th')$ de $\GL_f(\FM_q)$ par la correspondance de Green. Notons $e:=d/f.$ Les sous-quotients irr\'eductibles de $\times^e \o\pi_f(\th')$ (le produit $e$ fois) sont en bijection avec les $\End_{\GL_d(\FM_q)}(\times^e \o\pi_f(\th'))$-modules \`a droites simples. Howlett et Lehrer (voir \cite[3.6]{Dipper}) ont d\'efini un isomorphisme d'alg\`ebres
\[
\End_{\GL_d(\FM_q)}(\times^e\o\pi_f(\th'))\cong \HC_{q^f}(\SG_e),
\]
o\`u $\HC_{q^f}(\SG_e)$ d\'esigne l'alg\`ebre de Hecke du groupe sym\'etrique $\SG_e.$ Rappelons que $\HC_{q^f}(\SG_e)$ poss\`ede une base $\{T_w~|~w\in\SG_e\}$ telle que
$$
\begin{cases}
T_v\cdot T_w=T_{vw}, \text{ si $l(vw)>l(w)$};\\
T_vT_w=q^f T_{vw}+(q^f-1)T_w, \text{ sinon}.
\end{cases}
$$
L'ensemble des caract\`eres irr\'eductibles de $\HC_{q^f}(\SG_e)$ est param\'etr\'e par l'ensemble des partitions de $e:~(\l\vdash e)\mapsto S_\l(q^f),$ o\`u $S_\l(q^f)$ est le module de Specht associ\'e \`a $\l$ \cite{Dipper}. Si $q^f=1,$ $S_\l(1)$ est le module de Specht du groupe sym\'etrique $\SG_e.$
On en d\'eduit une bijection $$\l\mapsto \o\pi^\l_{\th'}:=S_\l(q^f)\otimes_{\HC_{q^f}(\SG_e)}(\times^e\o\pi_f(\th'))$$ entre l'ensemble des partitions de $e$ et l'ensemble des sous-quotients irr\'eductibles de $\times^e\o\pi_f(\th').$

\begin{Def}
Une repr\'esentation irr\'eductible $\o\pi$ de $\GL_d(\FM_q)$ sera dite {\em elliptique} si son image dans $R(\GL_d(\FM_q))$ n'est pas contenue dans le sous-groupe engendr\'e par les induites paraboliques de repr\'esentations de sous-groupes de Levi propres.
\end{Def}

\begin{lemme}\label{3.3::L1}
Soit $\o\pi$ une repr\'esentation irr\'eductible elliptique de $\GL_d(\FM_q).$ Alors, il existe un unique triple $(f,\th',i)$ avec $f|d,$ $\th'$ un caract\`ere $f$-primitif de $\FM_{q^f}^\times$ et $i\in\{0,\ldots,e-1\}$ avec $e:=d/f$ tel que $\o\pi\cong \o\pi^i_{\th'}:=\o\pi^{(i+1,1^{e-1-i})}_{\th'}.$
\end{lemme}
\begin{preuve}
D'apr\`es la classification d\^ue \`a Green des repr\'esentations irr\'eductibles de $\GL_d(\FM_q)$, le support cuspidal d'une repr\'esentation elliptique $\o\pi$ est de la forme $( (\GL_f(\FM_q))^{d/f},\otimes^{d/f}\o\pi_f(\th'))$ pour un entier $f|d$ et un caract\`ere $f$-primitif $\th'.$ Donc, on a $\o\pi\cong \o\pi^\l_{\th'}$ pour $\l\vdash e:=d/f.$

Pour $n\in\NM,$ notons $\o\BG^{\GL_{nf}(\FM_q)}_{f,\th'}$ la sous-cat\'egorie pleine des $\oQl$-repr\'esentations de $\GL_{nf}(\FM_q)$ de longueur finie engendr\'ee par les sous-quotients de $\times^n\o\pi_f(\th').$ On a donc une famille  d'\'equivalences de cat\'egories
\begin{align*}
\o\a^n_{\FM_q}:\o\BG^{\GL_{nf}(\FM_q)}_{f,\th'}&\To{}\mathop{\hbox {Mod-$\HC_{q^f}(\SG_n)$}}\\
V&\longmapsto \Hom_{\GL_{nf}(\FM_q)}(\times^n\o\pi_f(\th'),V)\\
M\otimes_{\HC_{q^f}(\SG_n)}\big(\times^n\o\pi_f(\th') \big)&\longmapsfrom M
\end{align*}
o\`u \hbox {Mod-$\HC_{q^f}(\SG_n)$} d\'esigne la cat\'egorie des $\HC_{q^f}(\SG_n)$-modules \`a droites de longueur finie.


\begin{fait}
Soient $\mu:=(\mu_1\geq\cdots\geq\mu_r)$ une partition de $e$ et $\SG_\mu:=\SG_{\mu_1}\times\cdots\times\SG_{\mu_r}$ le sous-groupe de $\SG_e$ associ\'e \`a $\mu.$ Posons
$$
x_\mu:=\sum_{w\in \SG_\mu}T_w\in\HC_{q^f}(\SG_e)\cong\End_{\GL_d(\FM_q)}(\times^e\o\pi_f(\th')).
$$
Alors, on a
\begin{description}
\item[(a)] $x_\mu\HC_{q^f}(\SG_e)=\Ind_{\HC_{q^f}(\SG_\mu)}^{\HC_{q^f}(\SG_e)}1;$

\item[(b)] $x_\mu(\times^e\o\pi_f(\th'))=\Ind_{\o P_{\mu\cdot f}}^{\GL_d(\FM_q)}\o\s,$ o\`u $\o P_{\mu\cdot f}$ d\'esigne le \para standard $\o P_{(\mu_1f,\ldots,\mu_rf)}$ de $\GL_d(\FM_q),$ et $\o\s$ est une repr\'esentation cuspidale du quotient r\'eductif de $\o P_{\mu\cdot f}.$
\end{description}
\end{fait}
\begin{preuve}
(a) d\'ecoule de \cite[4.2]{Dipper}. (b) d\'ecoule de \cite[6.2]{James-PLMS}. Notons qu'en fait on conna\^it la forme pr\'ecise de $\o\s.$
\end{preuve}

Rappelons que $S_\l(q^f)$ appara\^it avec multiplicit\'e un dans $\Ind_{\HC_{q^f}(\SG_\l)}^{\HC_{q^f}(\SG_e)}1,$ et que dans le groupe de Grothendieck de \hbox {Mod-$\HC_{q^f}(\SG_e)$}, $[\Ind_{\HC_{q^f}(\SG_\l)}^{\HC_{q^f}(\SG_e)}1]=[S_\l(q^f)]+\sum_{\mu>\l}m_{\l,\mu}S_\mu(q^f).$ On en d\'eduit que
\begin{align*}
[S_\l(q^f)]&=\sum_{\mu\geq\l}a_{\l,\mu}[\Ind_{\HC_{q^f}(\SG_\mu)}^{\HC_{q^f}(\SG_e)}1]\\
&=\sum_{\mu\geq\l}a_{\l,\mu} x_\mu(\HC_{q^f}(\SG_e)).
\end{align*}
Il s'ensuit que
\[
\o\pi^\l_{\th'}=\sum_{\mu\geq\l}a_{\l,\mu} x_\mu(\times^e\o\pi_f(\th')).
\]
Par d\'efinition, $\o\pi^{\l}_{\th'}$ est elliptique \ssi $a_{\l,(e)}\neq 0.$ Par ailleurs, il exist un isomorphisme d'alg\`ebres entre $\HC_{q^f}(\SG_e)$ et l'alg\`ebre du groupe $\CM[\SG_e]$ tel que $S_\mu(q^f)$ correspond \`a $S_\mu(1)$ et que $\Ind_{\HC_{q^f}(\SG_\mu)}^{\HC_{q^f}(\SG_e)}1$ correspond \`a $\CM[\SG_e/\SG_\l]$ ({\em cf.} \cite[4.3]{Dipper}). D'apr\`es \cite[2.3.17]{James}, $a_{\l,(e)}\neq 0$ \ssi $\l$ est de la forme $(i+1,1^{e-1-i})$ avec $i\in\{0,\ldots,e-1\}$.

\end{preuve}

\ali Soit $K$ une extension finie de $\QM_p.$ Notons $R(\GL_d(K))$ le groupe de Grothendieck des $\oQl$-repr\'esentations de longueur finie de $\GL_d(K).$ Comme d'habitude, pour $\pi_i$ une repr\'esentation de $\GL_{d_i}(K),$ $i=1,\ldots,t,$ on note $\times_i\pi_i:=\pi_1\times\cdots\times\pi_t$ l'induite normalis\'ee de la
repr\'esentation $\pi_1\otimes\cdots\otimes\pi_t$ du sous-groupe de Levi diagonal par blocs $\GL_{d_1}(K)\times\cdots\times \GL_{d_t}(K)$ de $\GL_{d_1+\cdots+d_t}(K),$ le long du parabolique triangulaire par blocs sup\'erieur $P_{(d_1,\ldots,d_t)}.$

Rappelons bri\`evement la classification des repr\'esentations elliptiques \cite[\S 2]{Dat-elliptic}.

\begin{Def}
Une repr\'esentation irr\'eductible $\pi$ de $\GL_d(K)$ est {\em elliptique} si son image dans $R(\GL_d(K))$ n'est pas contenue dans le sous-groupe engendr\'e par les induites paraboliques de repr\'esentations de \levis propres.
\end{Def}

Pour une repr\'esentation irr\'eductible $\rho$ de $D^\times,$ il existe un unique diviseur $f\in\NM$ de $d$ et une unique repr\'esentation supercuspidale irr\'eductible de $\GL_{f}(K)$ tels que $\JL(\rho)$ apparaisse dans l'induite parabolique standard normalis\'ee $|\det|^{\frac{1-e}{2}}\pi_\rho\times\cdots\times|\det|^{\frac{e-1}{2}}\pi_\rho,$ o\`u $e:=d/f$ et $\JL$ d\'esigne la correspondance de Jacquet-Langlands induisant une bijection entre les repr\'esentations irr\'eductibles de $D^\times$ et les s\'eries discr\`etes de $\GL_d(K).$ Notons $M_\rho$ le \levi standard $(\GL_{f}(K))^{e}$ et $\overrightarrow{\pi_\rho}:=|\det|^{\frac{1-e}{2}}\pi_\rho\otimes\cdots\otimes|\det|^{\frac{e-1}{2}}\pi_\rho,$ alors la paire $(M_\rho,\overrightarrow{\pi_\rho})$ est un repr\'esentant du support cuspidal de $\JL(\rho).$ Posons $S_\rho$ l'ensemble des racines simples du centre de $M_\rho$ dans l'alg\`ebre de Lie du parabolique triangulaire sup\'erieur $P_\rho$ de Levi $M_\rho.$ On peut identifier $S_\rho$ \`a l'ensemble $\{1,\ldots,e-1\}.$

\begin{fact}
Soit $\pi$ une repr\'esentation elliptique de $\GL_d(K),$ alors il existe une unique repr\'esentation irr\'eductible $\rho$ de $D^\times$ et un unique sous-ensemble $I$ de $S_\rho$ tels que $\pi$ peut s'\'ecrire de la mani\`ere unique sous la forme $$\pi_{\rho}^I:=\Cosoc\big(i^{\GL_d(K)}_{M_{\rho,I}}\big(\Soc\big(i^{M_{\rho,I}}_{M_\rho}(\overrightarrow{\pi_\rho})\big)\big)\big),$$
o\`u $M_{\rho,I}$ est le centralisateur du noyau commun des $\a\in I$ et le symbole $i$ d\'esigne l'induite normalis\'ee. Dans ce cas, $\pi_\rho^I$ a le m\^eme support cuspidal que $\JL(\rho).$ On a l'\'egalit\'e $\LJ(\pi^I_\rho)=(-1)^{|I|}[\rho]$ dans $R(D^\times)$ le groupe de Grothendieck des $\oQl$-repr\'esentations de longueur finie de $D^\times,$ o\`u $\LJ:R(\GL_d(K))\to R(D^\times)$ d\'esigne le transfert de Langlands-Jacquet \cite[2.1]{Dat-elliptic}.
\end{fact}

\ali\label{4.3Sec:2} La repr\'esentation elliptique $\pi_\rho^I$ de $\GL_d(K)$ est de niveau z\'ero, i.e. $(\pi_\rho^I)^{1+\varpi M_d(\OC)}\neq 0,$ \ssi $\rho$ est de niveau z\'ero, i.e. $\rho^{1+\Pi_D\OC_D}\neq 0.$ Bushnell et Henniart \cite{Bushnell-Henniart-level0} d\'ecrivent dans ce cas, la classification explicite des repr\'esentations irr\'eductibles de niveau z\'ero de $D^\times$ ainsi que la correspondance de Jacquet-Langlands associ\'ee.

Pour $n\in \NM$ un entier, notons $K_n$ l'extension non-ramifi\'ee de $K$ de degr\'e $n.$ Un caract\`ere mod\'er\'e \footnote{C'est-\`a-dire $\wt\th|_{1+\varpi \OC_{K_d}}$ est trivial.} $\wt\th:K_d^\times\to\oQl^\times$ sera dit {\em $f$-primitif} si $f$ est l'entier positif minimal tel qu'il existe un caract\`ere $\wt\th':K_f^\times\to \oQl^\times$ avec $\wt\th=\wt\th'\circ N_{K_d/K_f}.$ Notons dans ce cas $e:=d/f.$ Le couple $(K_f/K,\wt\th')$ est une {\em paire admissible mod\'er\'ee}, {\em cf.} \cite{Bushnell-Henniart-level0}.

Il existe une bijection $\wt\th\mapsto \rho(\wt\th)$ entre l'ensemble des classes d'\'equivalence des caract\`eres $\wt\th:K_d^\times\to \oQl^\times$ et l'ensemble des classes d'isomorphie des repr\'esentations irr\'eductibles de niveau z\'ero de $D^\times.$ En effet, si $\wt\th$ est $f$-primitif, on fixe une $K$-injection $K_f\into D,$ unique \`a conjugaison par un \'el\'ement de $D^\times$ pr\`es, et on identifie $K_f$ \`a une $K$-sous-alg\`ebre de $D^\times.$ Notons $B$ le centralisateur de $K_f$ dans $D.$ Alors, $B$ est une $K_f$-alg\`ebre \`a division centrale de dimension $e^2.$ On d\'esigne $\Nr_{B/K_f}:B^\times\to K_f^\times$ la norme r\'eduite et $U^1_D$ le sous-groupe $1+\Pi_D\OC_D$ de $D^\times.$ Donc $\wt\th$ induit un caract\`ere $\Theta$ du groupe $[D]_f=B^\times U^1_D$ ({\em cf.} \ref{3.2::S1}) par
$$
\Theta(bu)=\wt\th'(\Nr_{B/K_f}(b)),~b\in B^\times,~u\in U^1_D,
$$
o\`u $\wt\th'$ est le caract\`ere de $K_f^\times$ tel que $\wt\th=\wt\th'\circ N_{K_d/K_f}$. Alors, on d\'efinit la repr\'esentation irr\'eductible $\rho(\wt\th)$ de $D^\times$ comme suit
$$
\rho(\wt\th):=\Ind^{D^\times}_{[D]_{f}}\Theta,
$$
et on obtient de telle mani\`ere toutes les repr\'esentations de niveau z\'ero de $D^\times.$

\begin{theo}\label{4.3Prop:1}
Soit $\pi$ une repr\'esentation de niveau z\'ero de $\GL_d(K)$ telle que $\pi^{1+\varpi M_d(\OC)}$ soit irr\'eductible et elliptique en tant que repr\'esentation de $\GL_d(\FM_q),$ donc isomorphe \`a une $\o\pi^i_{\th'}$ avec $f|d$ et $\th'$ un caract\`ere $f$-primitif de $\FM_{q^f}^\times$ (voir \ref{3.3::L1}). Alors, $\pi$ est irr\'eductible elliptique de la forme $\pi^I_{\rho(\wt\th)}$ o\`u $\wt\th$ est un caract\`ere mod\'er\'e tel que $\wt\th|_{\OC^\times_{K_d}/1+\varpi \OC_{K_d}}\cong\th:=\th'\circ N_{\FM_{q^d}/\FM_{q^f}}$ et $I=\{1,\ldots,i\}$ ou $\{e-i,\ldots,e-1\}$ avec $e:=d/f,$ $0\leq i\leq e-1$ sous la convention $I=\emptyset$ si $i=0.$
\end{theo}

\begin{preuve}
Rappelons que le foncteur $M\mapsto M^{1+\varpi M_d(\OC)}$ induit une \'equivalence de cat\'egories entre la cat\'egorie des repr\'esentations lisses de niveau z\'ero de $\GL_d(K)$ et la cat\'egorie des modules sur l'alg\`ebre de Hecke $\HC(\GL_d(K),1+\varpi M_d(\OC)).$ Comme $\pi^{1+\varpi M_d(\OC)}$ est irr\'eductible et non nulle, $\pi$ est irr\'eductible.

Supposons que $\pi$ ne soit pas elliptique. Alors, on peut \'ecrire
$$
[\pi]=\sum_{P\subsetneqq \GL_d(K)}a_P[\Ind^{\GL_d(K)}_P \s_{M_P}],~a_P\in\ZM,
$$
o\`u $\s_{M_P}$ d\'esigne une repr\'esentation de niveau z\'ero du \levi $M_P$ de $P.$ Donc, en prenant les invariants sous $1+\varpi M_d(\OC),$ on a ({\em cf.} \cite[III. 3.14]{vigneras-book})
$$
\pi^{1+\varpi M_d(\OC)}=\o\pi^i_{\th'}=\sum_{\o P\subsetneqq \GL_d(\FM_q)}a_P\Ind^{\GL_d(\FM_q)}_{\o P} (\s_{M_{P}}^{M\cap(1+\varpi M_d(\OC))}),
$$
o\`u $\o P:= P(\OC)\pmod \varpi.$ On obtient alors une contradiction car $\o\pi^i_{\th'}$ est elliptique.

Comme $\pi^{1+\varpi M_d(\OC)}=\o\pi^i_{\th'},$ le couple $((\GL_f(\FM_q))^e,(\o\pi_f(\th'))^{\otimes e})$ est un repr\'esentant du support cuspidal de $\pi^{1+\varpi M_d(\OC)}.$ Donc, il existe un caract\`ere mod\'er\'e $\wt \th':K_f^\times\to \oQl^\times$ prolongeant $\th'$ tel que le support cuspidal de $\pi$ est \'egal \`a $((\GL_f(K))^e, \overrightarrow{\pi_f(\wt\th')}),$ o\`u $\pi_f(\wt\th')$ est la repr\'esentation supercuspidale de $\GL_f(K)$ associ\'ee \`a $\wt\th',$ {\em cf.} \cite[Prop. 8]{Bushnell-Henniart-level0}. Notons $\wt\th:=\wt\th'\circ N_{K_d/K_f},$ alors $\LJ(\pi)=\pm[\rho(\wt\th)].$ Donc on peut \'ecrire $\pi$ sous la forme $\pi_{\rho(\wt\th)}^I$ pour un sous-ensemble $I$ de $\{1,\ldots,e-1\}.$

Consid\'erons $\BG^{\GL_{nf}(K)}_{f,\th'}$ ($\forall n\in\NM$) la sous-cat\'egorie pleine des $\oQl$-repr\'esentations de longueur finie de $\GL_{nf}(K)$ form\'ee des objets dont tous les sous-quotients irr\'eductibles contiennent $((\GL_f(K))^n, \psi\cdot\overrightarrow{\pi_f(\wt\th')})$ dans leur support cuspidal, pour un certain caract\`ere non-ramifi\'e $\psi$ de $(\GL_f(K))^n.$ On sait que c'est une sous-cat\'egorie de Serre et que l'on a des \'equivalences
\begin{align*}
\BG^{\GL_{nf}(K)}_{f,\th'}&\To{}\mathop{\hbox {Mod-$\End_{\GL_{nf}(K)}(\times^n\pi_f(\wt\th'))$}}\\
V&\mapsto\Hom_{\GL_{nf}(K)}(\times^n\pi_f(\wt\th'),V),
\end{align*}
o\`u $\mathop{\hbox {Mod-$\End_{\GL_{nf}(K)}(\times^n\pi_f(\wt\th'))$}}$ d\'esigne la cat\'egorie des $\End_{\GL_{nf}(K)}(\times^n\pi_f(\wt\th'))$-modules \`a droites de longueur finie. Bushnell et Kutzko \cite[5.6.6]{Bushnell-Kutzko-book} ont d\'efini un isomorphisme d'alg\`ebres $$\End_{\GL_{nf}(K)}(\times^n\pi_f(\wt\th'))\cong \HC(n,q^f),$$ o\`u $\HC(n,q^f)$ est l'alg\`ebre de Iwahori-Hecke. Nous avons donc une \'equivalence de cat\'egories (bijective sur les objets simples)
\[
\a^n_K:\BG^{\GL_{nf}(K)}_{f,\th'}\To{}\mathop{\hbox {Mod-$\HC(n,q^f)$}}.
\]
Cette \'equivalence est compatible \`a l'\'equivalence $\o\a^n_{\FM_q},$ les induites paraboliques normalis\'ees et les foncteurs de Jacquet normalis\'es (voir \cite{Dat-elliptic}). Notons que $\pi=\pi_{\rho(\wt\th)}^I\in\BG^{\GL_{d}(K)}_{f,\th'}.$ Nous avons donc un diagramme commutatif
$$
\xymatrix{
\BG^{\GL_d(K)}_{f,\th'} \ar [d]_{\Hom_{1+\varpi M_d(\OC)}(1,-)}  \ar[r]^{\a^e_K} & \mathop{\hbox {Mod-$\HC(e,q^f)$}} \ar[d]^{\res.} \ar[r]^{(\a^e_{K_f})^{-1}} & \BG^{\GL_e(K_f)}_{1,\th'}\cong\BG^{\GL_e(K_f)}_{1,1} \ar[d]^{\Hom_{1+\varpi M_e(\OC_{K_f})}(1,-)} \\
\o\BG^{\GL_d(\FM_q)}_{f,\th'}  \ar[r]^{\o\a^e_{\FM_q}} & \mathop{\hbox {Mod-$\HC_{q^f}(\SG_e)$}} \ar[r]^{(\o\a^e_{\FM_{q^f}})^{-1}} & \o\BG^{\GL_e(\FM_{q^f})}_{1,\th'}\cong\o\BG^{\GL_e(\FM_{q^f})}_{1,1}
}
$$

En d\'ecorant de signes $'$ les objets relatifs au corps $K_f$ et $\FM_{q^f}$ selon le contexte, l'image de $\pi^I_{\rho(\wt\th)}$ (resp. $\o\pi^\l_{\th'}$) sous $(\a^e_{K_f})^{-1}\circ\a^e_K$ (resp. $(\o\a^e_{\FM_{q^f}})^{-1}\circ\o\a^e_{\FM_q} $) est de la forme $\pi'^I_1$ (resp. $\o\pi'^\l_1$) ({\em cf.} \cite[2.1.13]{Dat-elliptic}).

\begin{lemme}
On a $\l_I=(i+1,1^{(e-1-i)}).$
\end{lemme}
\begin{preuve}
Gr\^ace aux \'equivalences ci-dessus, on se ram\`ene \`a $\th'=1$ et $f=1.$ La repr\'esentation induite $\Ind^{\GL_e(\FM_{q^f})}_{\o P'_I}1$ contient la repr\'esentation irr\'eductible $\o\pi'^{\l_I}_1$ avec multiplicit\'e un, et on a de plus
$$
\Hom_{\GL_e(\FM_{q^f})}(\o\pi'^{\l_I}_1,\Ind^{\GL_e(\FM_{q^f})}_{\o P'_\mu}1)=0, \quad \text{si $\l_I<\mu$}.
$$

Dans le groupe de Grothendieck $R(\GL_e(K_f))$, on a l'\'egalit\'e
$$
[\pi'^I_{1}]=[\Ind^{\GL_e(K_f)}_{P'_I}1]+\sum_{J\supsetneqq I}(-1)^{|J\ba I|}[\Ind^{\GL_e(K_f)}_{P'_J}1].
$$
Ceci entra\^ine que dans $R(\GL_e(\FM_{q^f})),$ on a
$$
[\o\pi'^I_1]=[\Ind^{\GL_e(\FM_{q^f})}_{\o P'_I}1]+\sum_{J\supsetneqq I}(-1)^{|J\ba I|}[\Ind^{\GL_e(\FM_{q^f})}_{\o P'_J}1].
$$
Alors, on a
\begin{align*}
\dim\Hom_{\GL_e(\FM_{q^f})}(\o\pi'^{\l_I}_1,\o\pi'^i_1)
&=\dim\Hom_{\GL_e(\FM_{q^f})}(\o\pi'^{\l_I}_1,(\pi'^I_1)^{1+\varpi M_e(\OC_{K_f})})\\
&=\dim\Hom_{\GL_e(\FM_{q^f})}([\o\pi'^{\l_I}_1],[\o\pi'^{\l_I}_1]+\sum_{\mu>\l_I}a_\mu [\Ind^{\GL_e(\FM_{q^f})}_{\o P'_\mu}1])\\
&=1.
\end{align*}
Il s'ensuit que $\l_I=(i+1,1^{e-1-i}).$
\end{preuve}

Revenons \`a la preuve du th\'eor\`eme \ref{4.3Prop:1}, il nous reste \`a montrer que $I=\{1,\ldots,i\}$ ou $\{e-i,\ldots,e-1\}.$ Notons $r_{P'_{(1,\ldots,1)}}$ le foncteur de Jacquet normalis\'e associ\'e au sous-groupe de Borel $P'_{(1,\ldots,1)}$ de $\GL_e(K_f),$ et $\o U'_{(1,\ldots,1)}$ le radical unipotent de $\o P'_{(1,\ldots,1)}.$ On a
\begin{align*}
\dim_{\oQl}r_{P'_{(1,\ldots,1)}}(\pi'^I_{1})&=\dim_{\oQl}\Hom_{(1+\varpi \OC_{K'})^e}(1,r_{P'_{(1,\ldots,1)}}(\pi'^I_{1}))\\
&=\dim_{\oQl}(\o\pi'^i_1)^{\o U'_{(1,\ldots,1)}}\\
&=\binom{e-1}{i},
\end{align*}
car $\dim (\o\pi'^i_1)^{\o U'_{(1,\ldots,1)}}$ est la dimension du module de Specht associ\'e \`a la partition $(i+1,1^{e-1-i}).$

Par ailleurs, d'apr\`es \cite[(2.1.1) Fact]{Dat-Nagoya}, on a
$$
\dim_{\oQl}r_{P'_{(1,\ldots,1)}}(\pi'^I_{1})=\#\{w\in\SG_e~|~I=\{i\in\{0,\ldots,e-1\}~|~w(i-1)<w(i)\}\}.
$$
Par un calcul direct de ce sous-ensemble de $\SG_e,$ on obtient que
$$
\#\{w\in\SG_e~|~I=\{i\in\{0,\ldots,e-1\}~|~w(i-1)<w(i)\}\}=\binom{e-1}{i}
$$
\ssi $I=\{1,\ldots,i\}$ ou $\{e-i,\ldots,e-1\}.$ Ceci termine la preuve du th\'eor\`eme \ref{4.3Prop:1}.
\end{preuve}

\subsection{Cohomologie de la vari\'et\'e de Deligne-Lusztig}\label{S3.4}

On relie les repr\'esentations elliptiques $\o\pi^i_{\th'}$ dans 3.3 aux cohomologies de vari\'et\'e de Deligne-Lusztig $\DL^{d-1}_{\oFq}$ ({\em cf.} \ref{3.1::T1} (b)). Les lettres \'epaisses $\Gb,\Pb,\Lb,\Ub,\Vb...$ seront utilis\'ees pour les groupes alg\'ebriques lin\'eaires sur $\oFq.$ $\Gb$ sera toujours le groupe lin\'eaire $\GL_d$ sur $\oFq,$ et on fixe un morphisme de Frobenius $F:(a_{i,j})\mapsto(a_{i,j}^q).$ Soit $\Vb$ un sous-groupe unipotent de $\Gb,$ la vari\'et\'e de Deligne-Lusztig associ\'ee est d\'efinie par (voir \cite{bonnafe-rouquier-ihes})
$$
Y^\Gb_\Vb:=\{g\Vb\in\Gb/\Vb~|~g^{-1}Fg\in \Vb F(\Vb)\}.
$$

\ali\label{3.4::Para1} Fixons une injection d'alg\`ebres
$$\iota:\FM_{q^d}\into M_d(\FM_q)=\End_{\FM_q}(\FM_q^d).$$
Le centralisateur de $\iota(\FM_{q^d}^\times)$ dans $\Gb$ est un tore maximal $F$-stable $\Tb$ de $\Gb,$ d\'eploy\'e sur $\FM_{q^d}.$ $\Tb$ est un tore Coxeter de $\Gb,$ et $\Tb^F\cong \FM_{q^d}^\times.$ Notons $V:=\FM_{q}^d$ et $\o{V}:=\oFq\otimes_{\FM_q}V.$ Alors, on a
$$
\o{V}=\bigoplus^{d-1}_{i=0}\o{V}_i, \text{ avec } \o{V}_i=\{v\in\o{V}~|~\iota(\l)(v)=\l^{q^i}v,\forall \l\in\FM_{q^d}^\times\}.
$$
Posons $\Bb$ le sous-groupe de Borel associ\'e au drapeau complet
$$
0\subset\o{V}_0\subset\o{V}_0\oplus\o{V}_1\subset\cdots\subset\bigoplus_{i<d-1}\o{V}_i\subset\o{V},
$$
et $\Ub$ son radical unipotent.

Soit $d=ef.$ Notons $\Lb$ le centralisateur de $\iota(\FM_{q^f}^\times)$ dans $\Gb.$ $\Lb$ est un sous-groupe de Levi $F$-stable de $\Gb$ d\'efini sur $\FM_{q^f};$ $\Lb$ contient $\Tb,$ et $\Lb^F\cong\GL_e(\FM_{q^f}).$ Si l'on note
$$
\o{V}^j=\{v\in\o{V}~|~\iota(\l)(v)=\l^{q^j}v,\forall \l\in\FM_{q^f}^\times\}=\bigoplus_{i\equiv j \mod f}\o{V}_i,
$$
$\Lb$ est associ\'e \`a la d\'ecomposition $\o{V}=\bigoplus^{f-1}_{j=0}\o{V}^j.$ Posons $\Pb$ le \para associ\'e au drapeau
$$
0\subset\o{V}^0\subset\o{V}^0\oplus\o{V}^1\subset\cdots\subset\bigoplus_{j<f-1}
\o{V}^j\subset\o{V},
$$
et $\Vb$ son radical unipotent. Notons $\Bb_{\Lb}:=\Bb\cap\Lb=\Tb\ltimes\Ub_{\Lb}$ un \borel de $\Lb$ de radical unipotent $\Ub_\Lb.$

Soit $\th'$ un caract\`ere $f$-primitif de $\FM_{q^f}^\times.$ En composant avec la norme de $\FM_{q^d}$ sur $\FM_{q^f},$ on d\'efinit un caract\`ere $\th:=\th'\circ N_{\FM_{q^d}/\FM_{q^f}}$ de $\FM_{q^d}^\times.$ La classe de conjugaison Frobenius de $\th$ correspond \`a une classe de conjugaison d'\'el\'ements semi-simples $\{s\}$ dans le groupe dual $\Gb^*$ de $\Gb.$ On identifie $\Lb$ au groupe dual du centralisateur de $s$ dans $\Gb^*.$ Dans \cite[Thm. A']{bonnafe-rouquier-ihes}, Bonnaf\'e et Rouquier ont associ\'e \`a $s$ un idempotent central $e^{\Lb^F}_s\in\oQl\Lb^F.$ Soit $t\in\{s\}\cap \Tb^*$ correspondant \`a $\th;$ l'idempotent $e^{\Tb^F}_t$ est alors l'idempotent associ\'e au caract\`ere $\th.$

Posons $R^\Gb_{\Lb\subset\Pb}$ l'induction de Lusztig d\'efinie par somme altern\'ee \cite[\S 1]{lusztig-finiteness}. Le morphisme $R^\Gb_{\Lb\subset\Pb}$ induit une bijection
$$
R^\Gb_{\Lb\subset\Pb}:\o\BG^{\GL_{e}(\FM_{q^f})}_{1,\th'}\To{}\o\BG^{\GL_{d}(\FM_q)}_{f,\th'}.
$$

\begin{fact}\label{3.4F::1}
On a $R^\Gb_{\Lb\subset\Pb}((\th'\circ\det)\otimes\o\pi'^i_1)=(-1)^{d+e}\o\pi^i_{\th'},$ o\`u $\o\pi'^i_1$ signifie la repr\'esentation elliptique de $\GL_e(\FM_{q^f})$ associ\'ee au caract\`ere trivial et la partition $(i+1,1^{e-1-i}).$
\end{fact}
\begin{preuve}
D'apr\`es \cite[4.5]{Dipper} et \cite[Page 116]{FS}, $\o\pi^i_{\th'}=\e_{\Gb^F}\e_{\Lb^F}R^\Gb_{\Lb\subset\Pb}((\th'\circ\det)\otimes\o\pi'^i_1)$ avec $\e_{\Gb^F}=(-1)^d$ et $\e_{\Lb^F}=(-1)^{e}.$
\end{preuve}

\begin{prop}\label{4.3Prop::1}
Rappelons bri\`evement que la vari\'et\'e $\DL^{d-1}_{\oFq}$ est munie une action de $\GL_d(\FM_q),$ et qu'elle est un $\FM_{q^d}^\times$-torseur au-dessus de l'espace de Drinfeld $\O^{d-1}_{\oFq}$ sur $\oFq$ qui est le compl\'ementaire de tous les hyperplans $\FM_q$-rationnels dans $\PM^{d-1}_{\oFq}$ ({\em cf.} \cite[2.5.1]{Wang-Sigma1}). Notons $M(\th)$ la partie $\th$-isotypique d'un $\FM_{q^d}^\times$-module $M.$ Alors, en tant que repr\'esentations de $\GL_d(\FM_q),$ on a $$H^{d-1+i}_c(\DL^{d-1}_{\oFq},\oQl)(\th)\cong\o\pi^i_{\th'},~\forall i\in\{0,\ldots,e-1\}.$$
\end{prop}

\begin{preuve}
On utilise le r\'esultat de l'ind\'ependance de paraboliques prouv\'e dans \cite{Dat-parab}.  Nous avons des isomorphismes de vari\'et\'es
\begin{align*}
Y^\Gb_{\Ub}&\cong \DL^{d-1}_{\oFq}\\
Y^\Lb_{\Ub_\Lb}&\cong \DL^{e-1}_{\oFq},
\end{align*}
ainsi qu'un isomorphisme gr\^ace \`a \cite[Lemma 3]{lusztig-finiteness},
\ini\begin{equation}\label{4.3Equ:3}
Y^\Gb_\Vb\times_{\Lb^F}Y^\Lb_{\Ub_\Lb}\cong Y^\Gb_{\Vb\Ub_\Lb}.
\end{equation}
Notons que $\Vb\Ub_{\Lb}$ est le radical unipotent d'un \borel $\Bb':=(\Bb\cap\Lb)\Vb$ de $\Gb.$ D'apr\`es, \cite[Thm. B']{bonnafe-rouquier-ihes}, $Y^\Gb_\Vb$ (resp. $Y^\Lb_{\Ub_\Lb}$) est de dimension $e(d-e)$ (resp. $e-1$), et le complexe $R\G(Y^\Gb_\Vb)e^{\Lb^F}_s$ est concentr\'ee en degr\'e $e(d-e).$ En vertu de \cite[Prop. 5.3]{Dat-parab}, les parties $\th$-isotypiques des complexes de cohomologies de $Y^\Gb_\Ub$ et de $Y^\Gb_{\Vb\Ub_\Lb}$ sont isomorphes, \`a un d\'ecalage de la diff\'erence de leurs dimensions et un twist \`a la Tate pr\`es. Tenant compte du fait que $R\G(Y^\Lb_{\Ub_{\Lb}})e^{\Tb^F}_t=e^{\Lb^F}_sR\G(Y^\Lb_{\Ub_{\Lb}})e^{\Tb^F}_t$ (\cite[Thm. 11.4]{bonnafe-rouquier-ihes}), on a
\begin{align*}
H^{d-1+i}_c(\DL^{d-1}_{\oFq},\oQl)(\th)&=H^{e(d-e)+e-1+i}_c(Y^\Gb_\Vb\times_{\Lb^F}Y^\Lb_{\Ub_\Lb},\oQl)e^{\Tb^F}_t\\
&=H^{e(d-e)}(Y^\Gb_{\Vb})e^{\Lb^F}_s\otimes_{\oQl[\Lb^F]}H^{e-1+i}_c(Y^\Lb_{\Ub_{\Lb}})(\th)\\
&=(-1)^{e(d-e)}R^\Gb_{\Lb\subset\Pb}((\th'\circ\det)\otimes \o\pi'^i_1)\\
&=(-1)^{e(e-1)(f-1)}\o\pi^i_{\th'}=\o\pi^i_{\th'}.
\end{align*}
o\`u la derni\`ere \'egalit\'e d\'ecoule de \ref{3.4F::1}. D'o\`u l'\'enonc\'e.
\end{preuve}

On peut calculer les valeurs propres de Frobenius sur les cohomologies de $Y^\Gb_\Ub\cong\DL^{d-1}_{\oFq}.$

\begin{prop}\label{3.4::Prop1}
Le Frobenius $F^f$ agit sur $H^{d-1+i}_c(Y^\Gb_\Ub,\oQl)(\th)$ par le scalaire
\[
(-1)^{e(f-1)}\th'((-1)^{e-1})\cdot (q^f)^{\frac{d-e}{2}+i}.
\]
\end{prop}
\begin{preuve} En vertu de \cite[Prop. 5.3]{Dat-parab}, on a un isomorphisme $F$-\'equivariant
$$
H^{d-1+i}_c(Y^\Gb_\Ub)(\th)(\frac{(e-1)(d-e)}{2})\cong H^{e(d-e)+e-1+i}_c(Y^\Gb_{\Vb\Ub_\Lb})(\th).
$$
Comme $\Vb$ est $F^f$-stable, l'isomorphisme \ref{4.3Equ:3}:
$$
Y^\Gb_\Vb\times_{\Lb^F}Y^\Lb_{\Ub_\Lb}\cong Y^\Gb_{\Vb\Ub_\Lb}
$$
est $F^f$-\'equivariant. Donc il suffit de calculer la valeur propre de $F^f$ sur
$$
H^{e(d-e)+e-1+i}_c(Y^\Gb_\Vb\times_{\Lb^F}Y^\Lb_{\Ub_\Lb})(\th)(-\frac{(e-1)(d-e)}{2}).
$$
Comme pr\'ec\'edemment, on a l'expression suivante
$$
H^{e(d-e)+e-1+i}_c(Y^\Gb_\Vb\times_{\Lb^F}Y^\Lb_{\Ub_\Lb})(\th)= H^{e(d-e)}_c(Y^\Gb_\Vb)e^{\Lb^F}_s\otimes_{\oQl[\Lb^F]}H^{e-1+i}_c(Y^\Lb_{\Ub_\Lb})e^{\Tb^F}_t.
$$
Alors, le valeur propre de $F^f$ est \'egal \`a $\mu_1\cdot\mu_2,$ o\`u $\mu_1$ (resp. $\mu_2$) d\'esigne la valeur propre de $F^f$ sur $H^{e-1+i}_c(Y^\Lb_{\Ub_\Lb})e^{\Tb^F}_t$ (resp. $H^{e(d-e)}_c(Y^\Gb_\Vb)(-\frac{(e-1)(d-e)}{2})e^{\Lb^F}_s$).

Commen\c cons par le calcul de $\mu_1.$
\begin{lemme}
$\mu_1=\th'((-1)^{e-1})q^{fi}.$
\end{lemme}
\begin{preuve}
Ceci est essentiellement donn\'e par Digne et Michel \cite{Digne-Michel}. Rappelons que $Y^\Lb_{\Ub_{\Lb}}$ est isomorphe \`a la vari\'et\'e $\DL^{e-1}_{\FM_{q^f}}$ d\'efinie par l'\'equation:
$$
\det((X_i^{q^{fj}})_{0\leq i,j\leq e-1})^{q^f-1}=(-1)^{e-1}.
$$
 Notons que dans \cite{Digne-Michel}, les auteurs ont consid\'er\'e la vari\'et\'e d\'efinie par l'\'equation $$\det((X_i^{q^{fj}})_{0\leq i,j\leq e-1})^{q^f-1}=1.$$ Soit $g\in\GL_e(\FM_{q^f}),$ notons, suivant eux \cite[Page 12]{Digne-Michel} (not\'e $N^1_{\wdo}(\th)(g)$ dans {\em loc. cit.}),
\begin{align*}
N_{\th}(g)&:=\Trace(gF^f|H^{*}_c(Y^\Lb_{\Ub_\Lb})e^{\Tb^F}_t)\\
&:=\sum_k(-1)^k\Trace(gF^f|H^{k}_c(Y^\Lb_{\Ub_\Lb})e^{\Tb^F}_t).
\end{align*}
Rappelons que $s$ se correspond \`a $\th,$ et $\th=\th'\circ N_{\FM_{q^d}/\FM_{q^f}}.$

\begin{fait}

\begin{description}
    \item[(a)] En tant que repr\'esentations de $\GL_e(\FM_{q^f}),$ on a
  $$
  H^{k}_c(Y^\Lb_{\Ub_\Lb})e^{\Tb^F}_t=H^k_c(\O^{e-1}_{\o\FM_{q^f}},\oQl) \otimes \th'^{-1}\circ\det,~\forall k.
  $$
  \item[(b)] Pour tout $g\in\GL_e(\FM_{q^f}),$ on a
  \ini\begin{equation}\label{4.4Eq:1}
  N_\th(g)=\th'((-1)^{e-1})\cdot N_1(g)\cdot \th'^{-1}(\det(g));
  \end{equation}
\end{description}
\end{fait}
\begin{preuve}
(a) d\'ecoule du fait que l'\'el\'ement $(g,\z)\in \GL_e(\FM_{q^f})\times \FM_{q^d}^\times$ avec $\det(g)\cdot N_{\FM_{q^d}/\FM_{q^f}}(\z)=1$ stabilise une composante connexe g\'eom\'etrique fix\'ee de $\DL^{e-1}_{\FM_{q^f}}.$
 (b) d\'ecoule de \cite[V. Corollaire 3.3]{Digne-Michel}, avec le facteur suppl\'ementaire $\th'((-1)^{e-1})$ qui vient du fait que $\DL^{e-1}_{\FM_{q^f}}$ est d\'efinie par l'\'equation
$$
\det((X_i^{q^{fj}})_{0\leq i,j\leq e-1})^{q^f-1}=(-1)^{e-1}.
$$
On revoie les lectures \`a \cite[Thm. 3.1.12]{Wang-Sigma1} pour les d\'etails.
\end{preuve}

Posons $g$ un \'el\'ement r\'egulier elliptique de $\GL_e(\FM_{q^f}).$ D'apr\`es \cite[Prop. 3.3.9]{DOR}, on a
\begin{align*}
N_1(g)&=\sum_k(-1)^k\Trace(gF^f|H^k_c(\O^{e-1}_{\o\FM_{q^f}},\oQl))\\
&=\sum_{k=e-1}^{2e-2}(-1)^k q^{f(k-e+1)}\cdot(-1)^{2e-2-k}\Trace(g|1)\\
&=\sum_{k=e-1}^{2e-2}q^{f(k-e+1)}.
\end{align*}

Notons $\a_k$ la valeur propre de $F^f$ sur $H^k_c(\DL^{e-1}_{\o\FM_{q^f}},\oQl)e^{\Tb^F}_t.$ Alors, on a
\begin{align*}
N_\th(g)&=\sum_k(-1)^k\Trace(gF^f|H^k_c(\DL^{e-1}_{\o\FM_{q^f}},\oQl)e^{\Tb^F}_t)\\
&=\sum_k(-1)^k\a_k\Trace(g|H^k_c(\DL^{e-1}_{\o\FM_{q^f}},\oQl)e^{\Tb^F}_t)\\
&=\sum_k(-1)^k\a_k\th'^{-1}(\det(g))\Trace(g|H^k_c(\O^{e-1}_{\o\FM_{q^f}},\oQl))\\
&=\sum_{k=e-1}^{2e-2}(-1)^k\a_k\th'^{-1}(\det(g))\cdot(-1)^{2e-2-k}\Trace(g|1)\\
&=\th'^{-1}(\det(g))\cdot\sum_{k=e-1}^{2e-2}\a_k.
\end{align*}
Comme les $\a_k/q^{f(k-e+1)}$ sont des nombre de Weil de poids $0,$ en vertu de l'\'egalit\'e \ref{4.4Eq:1}, on a $\a_k=\th'((-1)^{e-1})q^{f(k-e+1)}.$ En particulier, on a $\mu_1=\th'((-1)^{e-1})q^{fi}.$
\end{preuve}

\begin{lemme}
$\mu_2=(-1)^{e(f-1)}q^{\frac{ef(f-1)}{2}}.$
\end{lemme}
\begin{preuve}
La strat\'egie de la d\'emonstration est la suivante. Notons $R_u(\Hb)$ le radical unipotent d'un groupe r\'eductif connexe $\Hb.$ On va choisir un \borel $\Bb$ dans $\Pb:=\Lb\Vb$ tel que l'on a un isomorphisme $F^f$-\'equivariant:
$$
Y^\Gb_{R_u(\Bb)}\cong Y^\Gb_{\Vb}\times_{\Lb^F}Y^\Lb_{R_u(\Bb\cap\Lb)}
$$
et que l'on conna\^it la valeur propre de $F^f$ sur la partie $\th$-isotypique de la cohomologie de $Y^\Lb_{R_u(\Bb\cap\Lb)}.$ Ensuite, on trouve un autre \borel $\Bb'$ tel que $\Bb^*\cap\Lb^*=\Bb'^*\cap\Lb^*,$ et que l'on conna\^it la valeur propre de $F^f$ sur la partie $\th$-isotypique de sa cohomologie. En vertu du r\'esultat de \cite{Dat-parab}, les valeurs propres de $F^f$ sur la partie $\th$-isotypique de la cohomologie de $Y^\Gb_{R_u(\Bb)}$ et de celle de $Y^\Gb_{R_u(\Bb')}$ sont les m\^emes \`a un twist de Tate explicit pr\`es. Donc, on en d\'eduit $\mu_2.$

Notons $V_1=V_2=\cdots=V_e:=\FM_q^f,$ et $\o{V}_i:=\oFq\otimes_{\FM_q}V_i.$ Fixons des injections d'alg\`ebres $\iota_i:\FM_{q^f}\into \End_{\FM_q}V_i,~\forall 1\leq i\leq e.$ Posons $V:=\bigoplus_{i=1}^e V_i$ et $\o{V}:=\oFq\otimes_{\FM_q}V.$ On identifie $\Gb$ au groupe de $\FM_q$-automorphisme de $\o{V}.$ Les $\iota_i$ induisent une injection
$$
\iota:=\iota_1\times\cdots\times \iota_e:(\FM_{q^f}^\times)^e\into \Aut_{\FM_q}(V_1)\times\cdots\times\Aut_{\FM_q}(V_e)\subset \Aut_{\FM_q}(V_1\oplus\cdots\oplus V_e)\into \Gb.
$$
Le centralisateur de $\iota((\FM_{q^f}^\times)^e)$ dans $\Gb$ est un tore maximal $F$-stable $\Tb$ dont la structure rationelle est donn\'ee par la restriction de scalaires $\Res_{\FM_{q^f}/\FM_q}(\GM_m)^e.$ De plus, $\Tb^F=\iota((\FM_{q^f}^\times)^e).$ L'action de $\FM_{q^f}^\times$ sur $\o{V}_i$ induite par $\iota_i$ fait une d\'ecomposition en sous-espaces de dimension 1
$$
\o{V}_i=\bigoplus_{j=0}^{f-1}\o{V}_{i,j},\text{ o\`u } \o{V}_{i,j}:=\{v\in\o{V}_i~|~\iota_i(\l)(v)=\l^{q^j}v,\forall\l\in\FM_{q^f}^\times\}.
$$
On note $\Bb'$ le \borel de $\Gb$ associ\'e au drapeau complet d\'efini par les sous-espaces vectoriels successifs
$$
0,\o{V}_{1,0},\o{V}_{1,1},\ldots,\o{V}_{1,f-1},\o{V}_{2,0},\o{V}_{2,1},\ldots,\o{V}_{2,f-1},\ldots,\o{V}_{e,0}\ldots\o{V}_{e,f-1},
$$
et $R_u(\Bb')$ son radical unipotent. Notons $\Lb'$ le \levi $F$-stable associ\'e \`a la d\'ecomposition
$$
\o{V}=\bigoplus_{i=1}^e\o{V}_i,
$$
et $\Pb'$ le \para $F$-stable associ\'e au drapeau
$$
0\subset\o{V}_1\subset\o{V}_1\oplus\o{V}_2\subset\cdots\subset\bigoplus_{j<e}\o{V}_j\subset\o{V}
$$
de radical unipotent $R_u(\Pb').$ Le \borel $\Bb'\cap\Lb'$ de $\Lb'$ est de radical unipotent $R_u(\Bb'\cap\Lb').$ On a $\Lb'^F\cong (\GL_f(\FM_q))^e.$ La vari\'et\'e de Deligne-Lusztig $Y^{\Lb'}_{R_u(\Bb'\cap\Lb')}$ (resp. $Y^{\Gb}_{R_u(\Pb')}$) s'identifie \`a $(\DL^{f-1}_{\FM_{q^f}})^e$ (resp. $\Gb^F/R_u(\Pb')^F$).

Notons que $\iota_i$ induit une structure de $\FM_{q^f}$-espace vectoriel sur $V_i.$ Donc l'injection $\iota=\iota_1\times\cdots\times\iota_e$ se factorise par
$$
(\FM_{q^f}^\times)^e\into\Aut_{\FM_{q^f}}(V_1)\times\cdots\times\Aut_{\FM_{q^f}}(V_e)\subset\Aut_{\FM_{q^f}}(\bigoplus_{i=1}^e V_i)\subset \Aut_{\FM_{q}}(\bigoplus_{i=1}^e V_i)\into\Gb.
$$
L'immersion diagonale $\D:\FM_{q^f}^\times\to(\FM_{q^f}^\times)^e$ induit une d\'ecomposition
$$
\o{V}=\bigoplus_{j=0}^{f-1}\o{V}^j, \text{ o\`u } \o{V}^j:=\{v\in\o{V}~|~\iota\circ\D(\l)(v)=\l^{q^j}v,\forall \l\in\FM_{q^f}^\times\}.
$$
En fait, $\o{V}^j=\bigoplus_{i=1}^{e}\o{V}_{i,j}.$ Notons $\Lb$ le \levi $F$-stable associ\'e \`a cette d\'ecomposition $\o{V}=\bigoplus_{j=0}^{f-1}\o{V}^j,$ et $\Pb$ le \para de $\Gb$ correspondant au drapeau
$$
0\subset\o{V}^0\subset\o{V}^0\oplus\o{V}^1\subset\cdots\subset\bigoplus_{j<f-1}\o{V}^j\subset\o{V},
$$
et $R_u(\Pb)$ le radical unipotent de $\Pb.$ On note $\Bb$ le \borel $F$-stable de $\Gb$ associ\'e au drapeau complet d\'efini par les sous-espaces vectoriels successifs
$$
0,\o{V}_{1,0},\o{V}_{2,0},\ldots,\o{V}_{e,0},\o{V}_{1,1},\o{V}_{2,1},\ldots,\o{V}_{e,1},\ldots,\o{V}_{1,f-1}\ldots\o{V}_{e,f-1},
$$
et $R_u(\Bb)$ son radical unipotent. Le tore maximal $\Tb$ est contenu dans $\Lb,$ et $\Bb\cap\Lb\supset\Tb$ est un \borel de $\Lb$ de radical unipotent $R_u(\Bb\cap\Lb).$ On a $\Lb^F\cong\GL_e(\FM_{q^f}),$ la vari\'et\'e de Deligne-Lusztig $Y^\Gb_{R_u(\Pb)}$ s'identifie \`a la vari\'et\'e $Y^\Gb_{\Vb}$ introduite dans \ref{3.4::Para1}.

Consid\'erons les deux vari\'et\'es $Y^\Gb_{R_u(\Bb)}$ et $Y^\Gb_{R_u(\Bb')}.$ D'apr\`es \cite[Lemma 3]{lusztig-finiteness}, on a des isomorphismes
\begin{align*}
Y^\Gb_{R_u(\Bb)}&\cong Y^\Gb_{R_u(\Pb)}\times_{\Lb^F}Y^\Lb_{R_u(\Bb\cap\Lb)}\\
Y^\Gb_{R_u(\Bb')}&\cong Y^\Gb_{R_u(\Pb')}\times_{\Lb'^F}Y^{\Lb'}_{R_u(\Bb'\cap\Lb')}.
\end{align*}
Par d\'efinition, $\Bb^*\cap\Lb^*=\Bb'^*\cap\Lb^*$ et $\Lb^*=C_{\Gb^*}(s).$ Donc on peut appliquer le r\'esultat de \cite{Dat-parab} sur l'ind\'ependance du \para . On obtient

\begin{align*}
H^{e(d-e)}_c(Y^\Gb_{R_u(\Pb)}\times_{\Lb^F}Y^\Lb_{R_u(\Bb\cap\Lb)})(\th)&=(-1)^{e(d-e)}\cdot(-1)^{d-e}\cdot H^{d-e}_c(Y^\Gb_{R_u(\Pb')}\times_{\Lb'^F}Y^{\Lb'}_{R_u(\Bb'\cap\Lb')})(\th)(\frac{(e-1)(d-e)}{2})\\
&=H^{d-e}_c(Y^\Gb_{R_u(\Pb')}\times_{\Lb'^F}Y^{\Lb'}_{R_u(\Bb'\cap\Lb')})(\th)(\frac{(e-1)(d-e)}{2}).
\end{align*}

Notons que
$$
H^{e(d-e)}_c(Y^\Gb_{R_u(\Pb)}\times_{\Lb^F}Y^\Lb_{R_u(\Bb\cap\Lb)})(\th)=H^{e(d-e)}_c(Y^\Gb_{R_u(\Pb)})e^{\Lb^F}_s\otimes_{\oQl[\Lb^F]}(\th'\circ\det|_{\Lb^F}\otimes\Ind_{(\Bb\cap\Lb)^F}^{\Lb^F}1),
$$
et que
$$
H^{d-e}_c(Y^\Gb_{R_u(\Pb')}\times_{\Lb'^F}Y^{\Lb'}_{R_u(\Bb'\cap\Lb')})(\th)=\Ind_{\Pb'^F}^{\Gb^F}\big((H^{f-1}_c(\O^{f-1}_{\o\FM_q},\LC_{\th'}))^{\otimes e}\big),
$$
o\`u $\LC_{\th'}$ d\'esigne le syst\`eme local sur $\O^{f-1}_{\o\FM_q}$ associ\'e \`a $\th'.$
On en d\'eduit l'\'egalit\'e suivante:
\[
H^{e(d-e)}_c(Y^\Gb_{R_u(\Pb)})(-\frac{(e-1)(d-e)}{2})e^{\Lb^F}_s\otimes_{\oQl[\Lb^F]}\th'\circ\det|_{\Lb^F}\otimes\Ind_{(\Bb\cap\Lb)^F}^{\Lb^F}1=\Ind_{\Pb'^F}^{\Gb^F}\big((H^{f-1}_c(\O^{f-1}_{\o\FM_q},\LC_{\th'}))^{\otimes e}\big).
\]
Comme $\th'$ est $f$-primitif, $H^{f-1}_c(\O^{f-1}_{\o\FM_q},\LC_{\th'})$ est cuspidale en tant que repr\'esentation de $\GL_f(\FM_q).$ D'apr\`es \cite[Thm. 3.1.12]{Wang-Sigma1}, la valeur propre de $F^f$ sur $H^{f-1}_c(\O^{f-1}_{\o\FM_q},\LC_{\th'})$ est \'egale \`a $(-1)^{f-1}q^{\frac{f(f-1)}{2}}.$ Comme $F^f$ agit trivialement sur les vari\'et\'es $\Lb^F/(\Bb\cap\Lb)^F$ et $\Gb^F/\Pb'^F,$ sa valeur propre sur $\th'\circ\det|_{\Lb^F}\otimes\Ind_{(\Bb\cap\Lb)^F}^{\Lb^F}1$ est \'egale \`a $1$ et sa valeur propre sur $\Ind_{\Pb'^F}^{\Gb^F}\big((H^{f-1}_c(\O^{f-1}_{\o\FM_q},\LC_{\th'}))^{\otimes e}\big)$ est \'egale \`a $(-1)^{e(f-1)}q^{\frac{ef(f-1)}{2}}.$ Rappelons que $Y^\Gb_{\Vb}\cong Y^{\Gb}_{R_u(\Pb)}.$ On en d\'eduit que $\mu_2$ est \'egale \`a $(-1)^{e(f-1)}q^{\frac{ef(f-1)}{2}}.$
\end{preuve}

Gr\^ace aux deux lemmes pr\'ec\'edents, le Frobenius $F^f$ agit sur $H^{d-1+i}_c(Y^\Gb_\Ub,\oQl)(\th)$ par le scalaire $$\mu_1\cdot\mu_2=(-1)^{e(f-1)}\th'((-1)^{e-1})\cdot (q^f)^{\frac{d-e}{2}+i}.$$
\end{preuve}

\subsection{La cohomologie \`a coefficients $\ell$-adiques}\label{Sec:9}

On \'etudie la partie non-supercuspidale de la cohomologie $\ell$-adique:
$$\Hb^{d-1+i}_{c,\oQl}:=\Hb^{d-1+i}_{c,\oZl}\otimes_{\oZl}\oQl=\Ind^{GDW}_{[GDW]_d}H^{d-1+i}_c(\Sig^{ca}_1,\oQl),~0\leq i\leq d-1.$$

\ali\label{Sec:7} Soient $f|d$ et $\th':\FM_{q^f}^\times\to\oQl^\times$ un caract\`ere $f$-primitif. Notons $\th:=\th'\circ N_{\FM_{q^d}/\FM_{q^{f}}},$ et $e=d/f$ comme pr\'ec\'edemment.

Notons $K_f\subset K^{ca}$ l'extension non-ramifi\'ee de $K$ de degr\'e $f.$ Le caract\`ere $\th$ (ou plut\^ot $\th'$) nous fournit un caract\`ere mod\'er\'ement ramifi\'e $\wt\th'$ de $K_f,$ i.e. il est trivial sur le sous-groupe $1+\varpi\OC_{K_f},$ d\'efini par la compos\'ee:
\begin{align*}
\wt\th': K_f\onto (\OC^\times_{K_f}/1+\varpi\OC_{K_f}) \times\varpi^{\ZM}\cong \FM_{q^f}^\times\times\varpi^{\ZM}&\To{}\oQl^\times\\
(\z,\varpi)&\mapsto \th'(\z)\cdot \th'((-1)^{e-1})
\end{align*}
Alors, la paire $(K_f/K,\wt\th')$ est une paire admissible mod\'er\'ee \cite[1.1]{Bushnell-Henniart-level0}.

Notons $\wt\th:=\wt\th'\circ N_{K_d/K_f},$ et $\rho(\wt\th):=\Ind_{[D]_f}^{D^\times}\Theta$ la repr\'esentation irr\'eductible de niveau z\'ero de $D^\times$ associ\'ee \`a $\wt\th,$ o\`u $\Theta=\wt\th'\circ\Nr_{B/K_f},$ {\em cf.} \ref{4.3Sec:2}. En particulier, on a $$\Th(\Pi_D^f)=\wt\th'(\Nr_{B/K_f}(\Pi_D^f))=\wt\th'((-1)^{e-1}\varpi)=1.$$

\begin{lemme}
En tant que repr\'esentation de $GW,$ on a $$\Hom_{D^\times}(\rho(\wt\th),\Hb^{d-1+i}_{c,\oQl})=\Ind^{GW}_{[GW]_{f}}\big(\Hom_{\FM_{q^d}^\times}(\th,H^{d-1+i}_c(\Sig^{ca}_1,\oQl))\big),$$ o\`u $[GW]_{f}=\{(g,w)\in GW~|~v(g,1,w)\in f\ZM\},$ et l'action d'un \'el\'ement $(g,w)\in [GW]_f$ sur $\Hom_{\FM_{q^d}^\times}(\th,H^{d-1+i}_c(\Sig^{ca}_1,\oQl))$ vient de celle de $(g,\Pi_D^{-v(g,1,w)},w)\in (GDW)^0$ sur $H^{d-1+i}_c(\Sig^{ca}_1,\oQl)$.
\end{lemme}
\begin{preuve}
D'apr\`es la r\'eciprocit\'e de Frobenius, on a
\begin{align*}
\Hom_{D^\times}(\rho(\wt\th),\Hb^{d-1+i}_{c,\oQl})&\cong \Hom_{[D]_f}(\Theta,\Ind^{GDW}_{[GDW]_d}H^{d-1+i}_c(\Sig_1^{ca},\oQl))\\
&\cong\Hom_{\FM_{q^d}^\times}\big(\th,(\Ind^{GDW}_{(GDW)^0\varpi^{\ZM}}H^{d-1+i}_c(\Sig_1^{ca},\oQl))^{\Pi_D^f}\big).
\end{align*}
Comme $[GDW]_f=[GDW]_d\Pi_D^{f\ZM},$ on peut prolonger l'action naturelle de $[GDW]_d$ sur $\Sig^{ca}_1$ en une action (non naturelle) de $[GDW]_f$ en faisant agir $\Pi_D^f$ trivialement. On a alors
\begin{align*}
\Hom_{D^\times}(\rho(\wt\th),\Hb^{d-1+i}_{c,\oQl})&\cong \Hom_{\FM_{q^d}^\times}(\th,\Ind^{GDW}_{[GDW]_f}H^{d-1+i}_c(\Sig_1^{ca},\oQl))\\
&\cong\Ind^{GW}_{[GW]_f}\Hom_{\FM_{q^d}^\times}(\th,H^{d-1+i}_c(\Sig_1^{ca},\oQl)).
\end{align*}
\end{preuve}

\begin{prop}\label{4.3Prop:2}
Il existe une repr\'esentation irr\'eductible elliptique de niveau z\'ero $\pi_i(\th)$ de $G$ telle que
$$
\Hom_{D^\times}(\rho(\wt\th),\Hb^{d-1+i}_{c,\oQl})\underset{GW}{\cong}\pi_i(\th)\otimes \Ind_{[W_K]_f}^{W_K}V.
$$
o\`u $V$ est un $\oQl$-espace vectoriel de dimension $1$ sur lequel $I_K$ agit via $I_K\onto I_K/I_{K(\varpi_t)}\cong\FM_{q^d}^\times\To{\th}\oQl^\times,$ avec $\varpi_t$ la racine $(q^d-1)$-i\`eme de $\varpi$ fix\'ee dans \ref{3.1::T1}, et $\varphi^f$ agit sur $V$ par un scalaire $\l_i$ qui est de la forme $(-1)^{e(f-1)}\th'((-1)^{e-1})\cdot (q^f)^{\frac{d-e}{2}+i}.$ 

De plus, $\Ind^{W_K}_{[W_K]_f}V$ est irr\'eductible de dimension $f,$ et $\pi_i(\th)$ est de la forme $\pi^I_{\rho(\wt\th_i)}$ pour un caract\`ere mod\'er\'e $\wt\th_i$ de $K_d^\times$ prolongeant $\th$ et un sous-ensemble $I$ de $\{1,\ldots,e\}$ satisfaisant $I=\{1,\ldots,i\}$ ou $\{e-i,\ldots,e-1\}.$
\end{prop}

\begin{preuve}
Notons $x=[\OC^d]$ le sommet standard de $\BC\TC$ et $G_x^+=1+\varpi M_d(\OC).$ Consid\'erons les $G_x^+$-invariants de $\Hom_{D^\times}(\rho(\wt\th),\Hb^{d-1+i}_{c,\oQl}).$ En vertu du lemme \ref{4.3Prop::1}, nous avons
\begin{align*}
\Hom_{D^\times}(\rho(\wt\th),\Hb^{d-1+i}_{c,\oQl})^{G^+_x}&\underset{G_x\times W_K}{\cong}\Ind^{G_x\times W_K}_{G_x\times [W_K]_f}\Hom_{\FM_{q^d}^\times}(\th,H^{d-1+i}_c(\Sig^{ca}_1,\oQl)^{G_x^+})\\
&\underset{G_x\times W_K}{\cong}\Ind^{G_x\times W_K}_{G_x\times [W_K]_f}\big(H^{d-1+i}_c(\DL^{d-1}_{\oFq},\oQl)(\th)\big)\\
&\underset{G_x\times W_K}{\cong}\Ind^{G_x\times W_K}_{G_x\times [W_K]_f} \o\pi^i_{\th'}.
\end{align*}
Gr\^ace \`a \ref{3.1::T1} (a) et (b), l'action de $I_K$ sur $H^{d-1+i}_c(\DL^{d-1}_{\oFq},\oQl)(\th)$ se factorise par $I_K/I_{K(\varpi_t)}$ via le caract\`ere $\th,$ et $\varphi^f$ y agit comme l'endomorphisme $\Frob^f.$ Or, l'action de $\Frob^f$ commute avec celle de $\GL_d(\FM_q)$ et $\o\pi^i_{\th'}$ est une repr\'esentation irr\'eductible de $\GL_d(\FM_q).$ $\Frob^f$ y agit par un scalaire $\l_i.$ D'apr\`es Prop. \ref{3.4::Prop1}, $\l_i$ est de la forme pr\'evue dans l'\'enonc\'e. Donc, on a
$$
\Hom_{D^\times}(\rho(\wt\th),\Hb^{d-1+i}_{c,\oQl})^{G^+_x}\underset{G_x\times W_K}{\cong}\o\pi^i_{\th'}\otimes\Ind_{[W_K]_f}^{W_K}V,
$$
avec $V$ un $\oQl$-espace vectoriel de dimension $1$ comme dans l'\'enonc\'e. Comme l'action de $\FM_{q^d}^\times$ sur $V$ se factorise par le caract\`ere $f$-primitif $\th',$ le $W_K$-module $\Ind^{W_K}_{[W_K]_f}V$ est irr\'eductible.


Posons
$$\pi_i(\th):=\Hom_{W_K}\big(\Ind_{[W_K]_f}^{W_K}V, \Hom_{D^\times}(\rho(\wt\th),\Hb^{d-1+i}_{c,\oQl})\big).$$
Alors, on a
\begin{align*}
(\pi_i(\th))^{G_x^+}&=\Hom_{W_K}\big(\Ind_{[W_K]_f}^{W_K}V,\Hom_{D^\times}(\rho(\wt\th),\Hb^{d-1+i}_{c,\oQl})^{G^+_x}\big)\\
&\cong \o\pi^i_{\th'}.
\end{align*}
D'apr\`es le th\'eor\`eme \ref{4.3Prop:1}, $\pi_i(\th)$ est de la forme $\pi^I_{\rho(\wt\th_i)}$ comme dans l'\'enonc\'e, et on a
$$
\Hom_{D^\times}(\rho(\wt\th),\Hb^{d-1+i}_{c,\oQl})\underset{GW}{\simeq}\pi^I_{\rho(\wt\th_i)}\otimes\Ind_{[W_K]_f}^{W_K}V.
$$
\end{preuve}

\begin{coro}\label{4.4Cor}
Si $i=0,$ la repr\'esentation $\pi_0(\th)=\pi_{\rho(\wt\th_0)}^\emptyset$ est g\'en\'erique, donc elle est ``de la s\'erie discr\`ete''.
\end{coro}
\begin{preuve}
Ceci d\'ecoule de \cite[2.1.11]{Dat-elliptic}.
\end{preuve}

\begin{coro}\label{4.3Cor1}
Le $W_K$-module $\ind^{W_K}_{[W_K]_f}V$ est isomorphe \`a $\s_{\rho(\wt\th)}(-i)$ comme dans le th\'eor\`eme A. (b) dans l'introduction.
\end{coro}

\begin{preuve}
C'est une cons\'equence directe du th\'eor\`eme principal de \cite{Bushnell-Henniart-level0}. En effet, comme d\'ecrit par Bushnell et Henniart, le $W_K$-module $\s_{\rho(\wt\th)}(-i)$ s'identifie \`a l'induction $\ind^{W_K}_{[W_{K_f}]}(\eta^{e(f-1)}_{K_f}\wt\th')$ normalis\'ee \`a la Hecke, o\`u $\eta_{K_f}$ est le caract\`ere non-ramifi\'e quadratique de $K^\times_f.$ Bien s\^ur, on voit $\eta^{e(f-1)}_{K_f}\wt\th'$ comme un caract\`ere de $W_{K_f}$ via l'inverse de morphisme d'Artin $\Art^{-1}_{K_f}:W_{K_f}\onto W_{K_f}^{ab}\simto K_f^\times.$ Donc la restriction \`a $I_{K_f}=I_K$ du caract\`ere $\eta^{e(f-1)}_{K_f}\wt\th'$ se factorise par le quotient $I_{K_f}\onto \FM_{q^f}^\times,$ et la valeur du caract\`ere $\eta^{e(f-1)}_{K_f}\wt\th'$ en $\varphi^f$ est \'egale \`a $$(-1)^{e(f-1)}\wt\th'(\Art^{-1}(\varphi^f))=(-1)^{e(f-1)}\wt\th'(\pi)=(-1)^{e(f-1)}\th'((-1)^{e-1}).$$ D'o\`u le corollaire.
\end{preuve}

Dans le reste de ce paragraphe, on d\'emontre le th\'eor\`eme suivant, qui nous dit en particulier que la repr\'esentation $\pi_i(\th)=\pi_{\rho(\wt\th_i)}^I$ est isomorphe \`a $\pi^I_{\rho(\wt\th)}.$

\begin{theo}\label{Prop::1}
Dans le groupe de Grothendieck $R(D^\times),$ on a l'\'egalit\'e $\LJ(\pi_i(\th))=(-1)^{i}[\rho(\wt\th)],$ i.e. $\pi_i(\th)$ est elliptique de type $\rho(\wt\th).$
\end{theo}



\begin{preuve}
Soit $K_f/K$ l'extension non-ramifi\'ee de degr\'e $f.$ Posons $g'\in\GL_e(K_f)$ la matrice donn\'ee par $x_{i,i+1}=1$ pour $1\leq i\leq e-1,$ $x_{e1}=\varpi,$ et les autres $x_{ij}=0.$ D\`es que l'on fait un choix de base de $K_f$ sur $K,$ pour tout \'el\'ement de $\GL_e(K_f),$ on le voit naturellement comme un \'el\'ement de $\GL_d(K)$ dont la classe de conjugaison ne d\'epend pas du choix de base.

\begin{lem}
Il existe un $\a\in\OC_{K_f}^\times$ tel que $g'\cdot \diag(1,\ldots,1,\a)\in \GL_e(K_f)$ soit un \'el\'ement elliptique dans $\GL_d(K),$ et que $\sum_{i=0}^{f-1}\wt\th'(\varphi^i(\a))\neq 0.$
\end{lem}
\begin{preuve}
Soit $\a$ un \'el\'ement de $\OC^\times_{K_f}$ qui engendre $K_f$ sur $K.$ Alors, son polyn\^ome caract\'eristique $P_\a(X)$ sur $K$ a des racines distinctes $a_1,\ldots,a_f$ dans $K^{ca}.$ Le polyn\^ome caract\'eristique de $g'\cdot \diag(1,\ldots,1,\a)\in \GL_e(K_f)$ sur $K$ est de la forme
$$
P(X)=(X^e-a_1\varpi)(X^e-a_2\varpi)\cdots (X^e-a_f\varpi)\in K[X].
$$
Il s'agit alors de montrer que $P(X)$ est irr\'eductible sur $K.$ Soit $f(X)$ un polyn\^ome unitaire de $K[X]$ divisant $P(X),$ comme les polyn\^omes $X^e-a_i\varpi$ sont irr\'eductibles dans $K_f[X],$ on a $f(X)=\prod_{i\in I}(X^e-a_i\varpi),$ o\`u $I$ est un sous-ensemble de $\{1,\ldots,f\}.$ Posons $g(X):=\prod_{i\in I}(X-a_i\varpi),$ alors $g(X)$ divise $P_\a(X)$ et $f(X)=g(X^e).$ Comme $f(X)\in K[X],$ $g(X)$ l'est aussi. En vertu de la irr\'eductibilit\'e de $P_\a(X)$ sur $K,$ on sait que $I=\emptyset$ ou $\{1,\ldots,f\}.$ Donc $P(X)$ est irr\'eductible sur $K.$


Comme $\th'$ est un caract\`ere $f$-primitif de $\FM_{q^f}^\times,$ les caract\`eres $\th'\circ\Frob_q^i,~0\leq i\leq f-1,$ sont distincts. D'apr\`es Artin, ils sont lin\'eairement ind\'ependants. En vertu du fait que
$$\wt\th'(\varphi^i(\a))=\th'(\varphi^i(\a) \pmod \varpi)=\th'(\Frob_q^i(\a \pmod\varpi)),$$
on a $\sum_{i=0}^{f-1}\wt\th'(\varphi^i(\a))\neq 0.$
\end{preuve}

Choisissons un $\a\in\OC_{K_f}^\times$ satisfaisant le lemme pr\'ec\'edent, et notons $$g_G:=g'\cdot\diag(1,\ldots,1,\a)$$ vu comme un \'el\'ement de $\GL_e(K_f)\into\GL_d(K).$ Fixons un $K$-injection $K_f\into D$ de sorte qu'il soit normalis\'e par $\Pi_D,$ et notons $B$ le centralisateur de $K_f$ dans $D.$ Alors, $g_G$ correspond \`a une classe de conjugaison d'\'el\'ements r\'eguliers elliptiques $\{g_D\}$ de $B^{\times}\subset D^\times.$ On va calculer la trace de $g_G$ (resp. $g_D$) sur $\pi_i(\th)$ (resp. $\rho(\wt\th)$) dans la proposition suivante.

\begin{prop}\label{Prop::4}
On a l'\'egalit\'e de caract\`eres suivante:
$$
\chi_{\pi_i(\th)}(g_G)=(-1)^{d-1+i}\chi_{\rho(\wt\th)}(g_D)\neq 0.
$$
\end{prop}
\begin{preuve}
Rappelons que $\rho(\wt\th)=\Ind_{[D]_f}^{D^\times}\Theta,$ o\`u $\Theta$ est le caract\`ere de $[D]_f$ d\'efini dans \ref{4.3Sec:2}. Comme $g_D\in B^\times\subset[D]_f$ correspond \`a $g_G,$ $(g_D)^e$ est conjugu\'e \`a $\a\varpi.$ Donc, on a $\Nr_{B/K_f}(\Pi_D^{i}g_D\Pi_D^{-i})=\Pi_D^i (-1)^{e-1} \a\varpi\Pi_D^{-i}=(-1)^{e-1}\varphi^i(\a)\varpi.$ Alors, par d\'efinition,
\begin{align*}
\chi_{\rho(\wt\th)}(g_D)&=\chi_{\Ind_{[D]_f}^{D^\times}\Theta}(g_D)=\sum_{i=0}^{f-1}\Theta(\Pi_D^{i}g_D\Pi_D^{-i})\\
&=\sum^{f-1}_{i=0}\wt\th'((-1)^{e-1}\varphi^i(\a)\varpi)=\sum_{i=0}^{f-1}\wt\th'\big(\varphi^i(\a)\big).
\end{align*}

Notons $\s:=\{s_1,\ldots,s_e\}$ la facette stable sous $g_G$ sur laquelle $g_G$ agit par la permutation $\{1,\ldots,e\}\in \SG_e.$ D'apr\`es \cite[Thm. III 4.16 et Lemma III 4.10]{SS-ihes}, on a
$$
\chi_{\pi_i(\th)}(g_G)=\Trace(g_G|\pi_i(\th)^{G^+_\s}).$$
En vertu du th\'eor\`eme \ref{thm2} et la dualit\'e de Poincar\'e, on a
$$
\pi_i(\th)^{G^+_\s}=H^{d-1+i}_c(\tau^{-1}(|\s|^*),\LC_{\th})=H^{d-1-i}(\tau^{-1}(|\s|^*),\LC_{\th^{-1}})^\vee,
$$
o\`u $\tau:\O^{d-1,ca}_K\to|\BC\TC|$ signifie le morphisme de r\'eduction (voir \ref{3.1::S1}), et $\LC_{\th}$ (resp. $\LC_{\th^{-1}}$) d\'esigne le syst\`eme local sur $\O^{d-1,ca}_K$ associ\'e \`a $\th$ (resp. $\th^{-1}$).

Notons que l'action de $g_G$ sur $H^{d-1+i}_c(\Sig^{ca}_1,\oQl)(\th)$ est induite par l'action naturelle de l'\'el\'ement $(g_G,\Pi_D^{-f})\in (GD)^0:=\{(g,\d)\in G\times D^\times~|~v(g,\d,1)=0\}$ sur $\Sig_1^{ca}$ ({\em cf.} \cite[3.1.5]{Dat-elliptic}), gr\^ace au fait que $\Th(\Pi_D^f)=1$ (voir \ref{Sec:7}). On a alors
\begin{align*}
\chi_{\pi_i(\th)}(g_G)&=\Trace\big((g_G,\Pi_D^{-f})|H^{d-1+i}_c(\tau^{-1}(|\s|^*),\LC_{\th})\big)\\
&=\Trace\big((g_G^{-1},\Pi_D^f)|H^{d-1-i}(\tau^{-1}(|\s|^*),\LC_{\th^{-1}})\big).
\end{align*}

\begin{lem}
Le morphisme de restriction $$H^{d-1-i}(\tau^{-1}(|\s|^*),\LC_{\th^{-1}})\To{} H^{d-1-i}(\tau^{-1}(|\s|),\LC_{\th^{-1}})$$ est un isomorphisme.
\end{lem}
\begin{preuve}
Notons $R\Psi_\eta$ le foncteur de cycles \'evanescents formel \`a la Berkovich \cite{berk-vanishing}. D'apr\`es Berkovich, il s'agit de montrer que le morphisme de restriction
$$
R\G(\o\O_\s,R\Psi_\eta(\LC_{\th^{-1}}))\To{}R\G(\o\O_\s^0,R\Psi_\eta(\LC_{\th^{-1}}))
$$
induit un isomorphisme. Choisissons un sommet $s\in\s,$ et notons les inclusions
$$
j_s:\o\O_s^0\into \o\O_s,~ i_\s:\o\O_\s\into \o\O_s~ \text{ et } j_\s:\o\O^0_\s\into\o\O_\s.
$$
En vertu des \'egalit\'es \cite[(2.1), (2.4)]{Wang-Sigma1}, on a
\begin{align*}
R\Psi_\eta(\LC_{\th^{-1}})|_{\o\O_\s}&=i_\s^*Rj_{s*}\big(R\Psi_\eta(\LC_{\th^{-1}})|_{\o\O_s^0}\big)\\
&=i_\s^*Rj_{s*}(\o p_{s*}\oQl)(\th^{-1}).
\end{align*}
Finalement, \cite[Lemme 3.2.2]{wang_DL} nous dit que
\begin{align*}
i_\s^*Rj_{s*}(\o p_{s*}\oQl)(\th^{-1})&=Rj_{\s*}j_\s^*i_\s^*Rj_{s*}(\o p_{s*}\oQl)(\th^{-1})\\
&=Rj_{\s*}(Rj_{s*}(\o p_{s*}\oQl)(\th^{-1})|_{\o\O^0_\s}).
\end{align*}
On en d\'eduit l'\'enonc\'e du lemme.
\end{preuve}

Gr\^ace au lemme pr\'ec\'edent, on a
$$
\chi_{\pi_i(\th)}(g_G)=\Trace\big((g_G^{-1},\Pi_D^f)|H^{d-1-i}(\tau^{-1}(|\s|),\LC_{\th^{-1}})\big).
$$
On va consid\'erer un certain sous-sch\'ema formel de $\wh\O:=\wh\O^{d-1}_{\OC}\wh\otimes_{\OC}\wh{\OC^{ca}}$ afin de calculer cette trace.

Soit $\s=\{\varpi
\Lambda_k\subsetneq\Lambda_0\subsetneq\cdots\subsetneq\Lambda_{k}\}$
un $k$-simplexe. On dit que $\s$ est de {\em type} $(e_0,\ldots,e_k)$ si
$e_0=\dim_{\FM_q}\Lambda_0/\varpi\Lambda_k$ et
$e_i=\dim_{\FM_q}\Lambda_i/\Lambda_{i-1}$ pour $1\leq i\leq k.$ On dit de plus que $\s$ est une {\em $f$-facette}, si $f|e_i,~\forall i.$ Lorsque $e_i=f~\forall i$ et $k=e-1,$ on dit que $\s$ est une {\em $f$-facette maximale}. On note $F(f)$ l'ensemble de $f$-facettes, et $F(f)^c:=\BC\TC\ba F(f).$

Soit $\s=\{s_0,\ldots,s_k\}\in F(f),$ on note
\begin{align*}
\o\O^f_\s:&=\o\O_\s\ba \bigcup_{\s\subset\omega, \omega\in F(f)^c}\o\O_\omega\\
&=\bigcup_{\s\subset\omega,\omega\in F(f)}\o\O_\omega^0.
\end{align*}
Il est \'evident que $\o\O^f_\s=\o\O_{s_0}^f\cap\cdots\cap\o\O^f_{s_k},$ et $\o\O^f_\s=\o\O_\s^0$ si $\s$ est une $f$-facette maximale. Pour $s$ un sommet de $\BC\TC,$ $\o\O^f_s=\o\O_s\ba\bigcup_{s\in\omega, \omega\in F(f)^c}\o\O_\omega$ est une sous-vari\'et\'e ouverte de $\o\O_s,$ donc elle est lisse.

Notons $\o\O(f):=\bigcup_{s\in\BC\TC^\circ}\o\O^f_s$ une vari\'et\'e ouverte de la fibre sp\'eciale $\o\O$ de $\wh\O,$ et $\O(f):=\spe^{-1}(\o\O(f))$ o\`u $\spe$ d\'esigne le morphisme de sp\'ecialisation (voir \ref{3.1::S1}). Posons $\wh\O(f)$ le compl\'et\'e de $\wh\O$ le long de $\o\O(f),$ alors $\O(f)$ s'identifie \`a la fibre g\'en\'erique de $\wh\O(f),$ {\em cf.} \cite[Prop. 1.3]{berk-vanishing}.

Calculons le groupe de cohomologie $H^{d-1-i}(\tau^{-1}(|\s|),\LC_{\th^{-1}})$ en utilisant ce mod\`ele entier.
$$
\xymatrix{
  \O(f) \ar@{^(->}[r]^{j} & \wh\O(f) \ar@{<-^)}^{i}[r] & \o\O(f)
  }
$$

Notons tout d'abord que le faisceau $\LC_{\th^{-1}}$ se prolonge \`a $\wh\O(f),$ i.e. il existe un syst\`eme local de rang un $\wh\LC_{\th^{-1}}$ sur $\wh\O(f)$ tel que $j^*\wh\LC_{\th^{-1}}\cong\LC_{\th^{-1}}.$ En effet, d'apr\`es \ref{3.1::T1}, pour chaque sommet $s\in\BC\TC^\circ,$ $\LC_{\th^{-1}}$ se prolonge en un faisceau sur $\o\O_s^0,$ qui est isomorphe au syst\`eme local $\FC^{\wb}_{\th^{-1}}$ dans \cite[Thm. 7.7]{bonnafe-rouquier-ihes}, ici $\wb=(1,\ldots,d)\in\SG_d.$ Alors, d'apr\`es {\em loc. cit.}, si $\omega$ est une $f$-facette contenant $s,$ on sait que $\FC^{\wb}_{\th^{-1}}$ se prolonge sur la strate $\o\O_\omega^0.$ C'est-\`a-dire le syst\`eme local $\LC_{\th^{-1}}$ se prolonge sur toutes les $\o\O_\omega^0,$ pour $\omega\in F(f),$ donc il se prolonge sur $\o\O(f),$ car $\o\O(f)=\bigcup_{\omega\in F(f)}\o\O^0_\omega.$ Comme $\wh\O(f)$ est $\varpi$-adiquement complet, ce faisceau se rel\`eve uniquement comme un faisceau $\wh\LC_{\th^{-1}}$ sur $\wh\O(f)$ prolongeant $\LC_{\th^{-1}}.$

Notons $R\Psi_\eta^f$ le foncteur de cycles \'evanescents formel \`a la Berkovich \cite{berk-vanishing} associ\'e au mod\`ele entier $\wh\O(f)$ de $\O(f).$ D'apr\`es Berkovich, on a
\begin{align*}
H^{d-1-i}(\tau^{-1}(|\s|),\LC_{\th^{-1}})&=\HM^{d-1-i}(\o\O_\s^f,R\Psi_\eta^f\LC_{\th^{-1}})\\
&=\HM^{d-1-i}(\o\O_\s^0,R\Psi_\eta^f\LC_{\th^{-1}}),
\end{align*}
car $\s$ est $f$-maximale. En vertu de la formule de projection, on obtient que
$$
R\Psi_\eta^f\LC_{\th^{-1}}|_{\o\O^0_\s}=\wh\LC_{\th^{-1}}|_{\o\O^0_\s}\otimes^L R\Psi_\eta^f\oQl|_{\o\O^0_\s}.
$$
Rappelons que dans ce cas (o\`u $\wh{\O}(f)$ est un mod\`ele semi-stable de $\O(f)$), on conna\^it les faisceaux de cycles \'evanescents ({\em cf.} \cite[Thm. 3.2]{Illusie-monodromie})
\ini\begin{align}\label{Eq::6}
R^{0}\Psi_\eta^f\oQl|_{\o\O_\s^0}&=\oQl;\\
R^{1}\Psi_\eta^f\oQl|_{\o\O_\s^0}&=(\oplus_{i=1}^{e} (\oQl)_{i} / \oQl \text{ diagonal})(-1);\notag \\
R^{q}\Psi_\eta^f\oQl|_{\o\O_\s^0}&=\bigwedge^q R^{1}\Psi^f_\eta\oQl|_{\o\O_\s^0}. \notag
\end{align}

Gr\^ace \`a \cite[Lemma 3]{lusztig-finiteness}, $\o\O^0_\s$ est isomorphe \`a $(\O^{f-1}_{\oFq})^e.$ Donc, on a
\ini\begin{equation}\label{Eq::7}
H^{q}_c(\o\O^0_\s,\wh\LC_{\th}|_{\o\O^0_\s})=\left\{
                                               \begin{array}{ll}
                                                 H^{f-1}_c(\O^{f-1}_{\oFq},\LC_{\th'})^{\otimes e}, & \hbox{si $q=e(f-1)$;} \\
                                                 0, & \hbox{si $q\neq e(f-1)$.}
                                               \end{array}
                                             \right.
\end{equation}
Ici, $\LC_{\th'}$ d\'esigne le syst\`eme local sur $\O^{f-1}_{\oFq}$ associ\'e au caract\`ere $\th':\FM_{q^f}^\times\to\oQl^\times$ qui est $f$-primitif. En particulier, en tant que repr\'esentation de $\GL_f(\FM_q),$ $H^{f-1}_c(\O^{f-1}_{\oFq},\LC_{\th'})$ est la repr\'esentation cuspidale $\o\pi_f(\th').$

Maintenant, on peut terminer le calcul de la trace de $g_G.$
\begin{align*}
\chi_{\pi_i(\th)}(g_G)&=\Trace\big((g^{-1}_G,\Pi_D^f)|\HM^{d-1-i}(\o\O_\s^0,R\Psi_\eta^f\LC_{\th^{-1}})\big)\\
&=\Trace\big((g_G^{-1},\Pi_D^f)|H^{e(f-1)}(\o\O^0_\s,\wh\LC_{\th^{-1}}|_{\o\O^0_\s}\otimes R^{e-1-i}\Psi_\eta^f\oQl|_{\o\O^0_\s})\big)\\
&=\Trace\big((g_G,\Pi_D^{-f})|H_c^{f-1}(\O^{f-1}_{\oFq},\LC_{\th'})^{\otimes e}\big)\cdot \Trace\big((g_G^{-1},\Pi_D^{f})|R^{e-1-i}\Psi_\eta^f\oQl|_{\o\O^0_\s}\big)
\end{align*}

Comme le morphisme de p\'eriode $\xi_{Dr}:\MC^{ca}_{Dr,0}\to \O^{d-1,ca}_K$ est $G\times D^\times$-\'equivariant si l'on munit $\O^{d-1,ca}_K$ de l'action naturelle de $G$ et de l'action triviale de $D^\times$ (voir \cite[3.1.1 (ii)]{Dat-elliptic}), $(g_G^{-1},\Pi^{f}_D)$ permute les sommets de $\s=\{s_1,\ldots,s_{e}\}$ comme la permutation $(1,\ldots,e)\in\SG_e.$ Donc il agit sur $R^{1}\Psi^f_\eta\oQl|_{\o\O^0_\s}$ via cette permutation. Il s'ensuit que
$$
\Trace((g_G^{-1},\Pi_D^f)|R^{e-1-i}\Psi_\eta^f\oQl|_{\o\O^0_\s})=(-1)^{e-1-i}.
$$

D'apr\`es la r\`egle de Koszul \cite[Sommes trig. 2.4*, Cycle 1.3]{SGA4.5}, on sait que $(g',\Pi_D^{-f})$ agit sur $H_c^{f-1}(\O^{f-1}_{\oFq},\LC_{\th'})^{\otimes e}$ par la formule
$$
v_1\otimes v_2\otimes\cdots\otimes v_e\mapsto \wt\th'(\Nr_{B/K_f}(\Pi_D^{-f}))\cdot(-1)^{(f-1)^2(e-1)}v_2\otimes\cdots\otimes v_e\otimes v_1.
$$
En vertu de la formule de caract\`ere pour l'induite tensorielle, on a
\begin{eqnarray*}
\lefteqn{\Trace\big((g_G,\Pi_D^{-f})|H_c^{f-1}(\O^{f-1}_{\oFq},\LC_{\th'})^{\otimes e}\big)}\\
&=&\wt\th'(\Nr_{B/K_f}(\Pi_D^{-f}))\cdot(-1)^{(f-1)^2(e-1)}\Trace(\a|H_c^{f-1}(\O^{f-1}_{\oFq},\LC_{\th'}))\\
&=&(-1)^{(f-1)^2(e-1)}\Trace(\a|H_c^{f-1}(\O^{f-1}_{\oFq},\LC_{\th'})).
\end{eqnarray*}
La derni\`ere trace a \'et\'e calcul\'ee par Deligne et Lusztig \cite[Cor. 7.2]{deligne-lusztig}, elle est de la forme
$$
(-1)^{f-1}\sum^{f-1}_{i=0}\th'(\Frob_q^i(\a\pmod \varpi)).
$$
Donc, on a
\begin{align*}
\chi_{\pi_i(\th)}(g_G)&=(-1)^{e-1-i}\cdot(-1)^{(f-1)^2(e-1)}\cdot(-1)^{f-1}\sum^{f-1}_{i=0}\th'(\Frob_q^i(\a\pmod \varpi))\\
&=(-1)^{d-1-i}\cdot\sum^{f-1}_{i=0}\wt\th'(\varphi^i(\a))\\
&=(-1)^{d-1-i}\cdot \chi_{\rho(\wt\th)}(g_D).
\end{align*}
\end{preuve}

\begin{prop}\label{Prop::5}
On a $\pi_0(\th)\cong \JL(\rho(\wt\th)),$ o\`u $\JL$ d\'esigne la correspondance de Jacquet-Langlands. 
\end{prop}
\begin{preuve}
Gr\^ace au corollaire \ref{4.4Cor}, $\pi_0(\th)$ est une s\'erie discr\`ete de $G.$ D'apr\`es la proposition \ref{Prop::4}, on a les \'egalit\'es suivantes:
\ini\begin{align}\label{4.4Eq}
\chi_{\pi_0(\th)}(g_G)&=(-1)^{d-1}\chi_{\rho(\wt\th)}(g_D)\notag\\
&=\chi_{\JL(\rho(\wt\th))}(g_G).
\end{align}

D'apr\`es \ref{Eq::7}, les deux s\'eries discr\`etes $\pi_0(\th)$ et $\JL(\rho(\wt\th))$ contiennent le m\^eme type simple $(J_\th,\l_\th),$ o\`u $J_\th$ est isomorphe \`a $G^+_\s$ pour $\s$ une facette $f$-maximale, et $\l_\th$ est la repr\'esentation d'inflation de $J_\th,$ via le morphisme $J_\th\onto \GL_f(\FM_q)^e,$ de la repr\'esentation $(\o\pi_f(\th'))^{\otimes e},$ {\em cf.} \cite[5.2]{Bushnell-Henniart-level0}. Donc, leurs types simples \'etendus \cite{Bushnell-Henniart-level0} sont conjugu\'es \`a un caract\`ere non-ramifi\'e pr\`es. Plus pr\'ecis\'ement, un type simple \'etendu qui prolonge $(J_\th,\l_\th)$ est un couple $(\Jb_\th,\Lamb),$ o\`u $\Jb_\th$ s'identifie au stabilisateur de $\s$ dans $G,$ et $\Lamb$ est une repr\'esentation de $\Jb_\th$ telle que $\Lamb|_J\cong \l_\th,$ {\em cf.} \cite[4.1]{Bushnell-Henniart-level0}. D'apr\`es \cite[Prop. 7]{Bushnell-Henniart-level0}, on a une bijection entre les s\'eries discr\`etes contenant $(J_\th,\l_\th)$ et les types simples \'etendus prolongeant $(J_\th,\l_\th).$   Notons $(\Jb_\th,\Lamb_{\pi_0(\th)})$ (resp. $(\Jb_\th,\Lamb_{\JL (\rho(\wt\th))})$) le type simple \'etendu associ\'e \`a $\pi_0(\th)$ (resp. $\JL(\rho(\wt\th))$). Gr\^ace \`a \cite[Lemma 10]{Bushnell-Henniart-level0}, il existe un caract\`ere non-ramifi\'e $\chi$ de $K^\times$ tel que $\Lamb_{\pi_0(\th)}\cong\Lamb_{\JL (\rho(\wt\th))}\otimes \chi_{\Jb_\th},$ o\`u $\chi_{\Jb_\th}=\chi\circ\det|_{\Jb_\th}.$

Comme le type simple \'etendu $(\Jb_\th,\Lamb_{\pi_0(\th)})$ (resp. $(\Jb_\th,\Lamb_{\JL (\rho(\wt\th))})$) est de multiplicit\'e un dans $\pi_0(\th)$ (resp. $\JL(\rho(\wt\th))$), {\em cf.} \cite{Bushnell-Henniart-level0}, on sait que $\pi_0(\th)^{G^+_\s}\cong\Lamb_{\pi_0(\th)}$ (resp. $\JL(\rho(\wt\th))^{G^+_\s}\cong \Lamb_{\JL (\rho(\wt\th))})$).

Notons que $g_G\in \Jb_\th.$ D'apr\`es \cite[Thm. III 4.16 et Lemma III 4.10]{SS-ihes}, on a
$$
\chi_{\pi_0(\th)}(g_G)=\Trace(g_G|\Lamb_{\pi_0(\th)}),
$$
et idem pour $\JL(\rho(\wt\th)).$ Donc, en vertu de \ref{4.4Eq}, on a
\begin{align*}
\Trace(g_G|\Lamb_{\JL (\rho(\wt\th))})&=\Trace(g_G|\Lamb_{\pi_0(\th)})\\
&=\Trace(g_G|\Lamb_{\JL (\rho(\wt\th))}\otimes \chi_{\Jb_\th})\\
&=\Trace(g_G|\Lamb_{\JL (\rho(\wt\th))})\cdot \chi_{\Jb_\th}(g_G).
\end{align*}
Donc $\chi_{\Jb_\th}(g_G)=1.$ Il s'ensuit que $\chi_{\Jb_\th}=1.$
\end{preuve}

L'\'enonc\'e du th\'eor\`eme \ref{Prop::1} d\'ecoule de la proposition suivante.
\begin{prop}
Nous avons $\LJ(\pi_i(\th))=(-1)^{i} [\rho(\wt\th)],~\forall i\in\{1,\ldots,e-1\}.$
\end{prop}
\begin{preuve}
Gr\^ace \`a la proposition \ref{4.3Prop:2}, $[\rho(\wt\th_i)]=(-1)^i\LJ(\pi_i(\th))$ pour un caract\`ere mod\'er\'e $\wt\th_i$ prolongeant $\th.$ Comme $\wt\th_i$ prolonge $\th,$ le support cuspidal de $\pi_i(\th)$ contient le type simple $((\GL_f(\OC))^e,(\o\pi_f(\th'))^{\otimes e})$ du Levi standard. Comme les repr\'esentations elliptiques $\pi_{\rho(\wt\th_i)}^I(\cong\pi_i(\th))$ et $\pi_{\rho(\wt\th_i)}^\emptyset$ ont le m\^eme support cuspidal, on obtient que $\pi_{\rho(\wt\th_i)}^\emptyset$ contient le type simple $(J_\th,\l_\th).$

D'apr\`es \cite[2.1.14]{Dat-elliptic}, on a l'\'egalit\'e suivante
$$
\chi_{\pi_{\rho(\wt\th_i)}^\emptyset}(g_G)=(-1)^{|I|}\chi_{\pi_{\rho(\wt\th_i)}^I}(g_G).
$$
En vertu des propositions \ref{Prop::4} et \ref{Prop::5}, nous avons
\begin{align*}
\chi_{\pi_{\rho(\wt\th_i)}^\emptyset}(g_G)&=\chi_{\pi_0(\th)}(g_G)\\
&=\chi_{\JL(\rho(\wt\th))}(g_G).
\end{align*}
En r\'esum\'e, les deux s\'eries discr\`etes $\pi_{\rho(\wt\th_i)}^\emptyset$ et $\JL(\rho(\wt\th))$ contiennent le m\^eme type simple $(J_\th,\l_\th),$ et leurs caract\`eres ont la m\^eme valeur en l'\'el\'ement $g_G.$ Comme dans la preuve de la proposition \ref{Prop::5}, ils sont isomorphes, i.e. $\pi_{\rho(\wt\th_i)}^\emptyset\cong \JL(\rho(\wt\th)).$ Donc $\rho(\wt\th_i)\cong\rho(\wt\th).$
\end{preuve}
\end{preuve}

\def\cprime{$'$} \def\cprime{$'$} \def\cprime{$'$} \def\cprime{$'$}
  \def\cprime{$'$} \def\cprime{$'$} \def\cprime{$'$}

\noindent
\textsc{Haoran Wang}, Max-Planck-Institut f\"ur Mathematik, Vivatsgasse 7, 53111 Bonn, Germany

\noindent\texttt{haoran.wang@mpim-bonn.mpg.de}

\end{document}